\documentclass[reqno,10pt]{amsart}
\usepackage{amscd,amssymb,verbatim}
\numberwithin{equation}{section} \theoremstyle{plain}

\newcommand\alp{\alpha}         
\newcommand\bet{\beta}
\newcommand\gam{\gamma}         \newcommand\Gam{\Gamma}
\newcommand\del{\delta}         
\newcommand\eps{\varepsilon}
\newcommand\zet{\zeta}
\newcommand\tet{\theta}

\newcommand\lam{\lambda}                \newcommand\Lam{\Lambda}
\newcommand\sig{\sigma}         \newcommand\Sig{\Sigma}

\newcommand\ome{\omega}         \newcommand\Ome{\Omega}

\newcommand\calC{{\mathcal{C}}}

\newcommand\calE{{\mathcal{E}}}

\newcommand\calI{{\mathcal{I}}}

\newcommand\calM{{\mathcal{M}}}

\newcommand\calO{{\mathcal{O}}}

\newcommand\calU{{\mathcal{U}}}

\newcommand\calZ{{\mathcal{Z}}}

\newcommand\RR{\mathbb{R}}

\newcommand\ZZ{\mathbb{Z}}

\newcommand\CC{\mathbb{C}}

\newcommand\NN{\mathbb{N}}

 \newcommand\grg{{\mathfrak{g}}}

 \newcommand\grm{{\mathfrak{m}}}
 
 \newcommand\gro{{\mathfrak{o}}}

\newcommand\nek{,\ldots,}
\newcommand\sdp{\times \hskip -0.3em {\raise 0.3ex
\hbox{$\scriptscriptstyle |$}}} % semidirect product

\newcommand\Det{\operatorname{Det}}

\newcommand\End{\operatorname{End\,}}

\newcommand\RANK{\operatorname{rank}}

\newcommand\GL{\operatorname{GL}}

\newcommand\Id{\operatorname {Id}}

\newcommand\IM{\operatorname{Im}}

\newcommand\MOD{\operatorname{mod}}

\newcommand\Ker{\operatorname{Ker}}

\newcommand\rank{\operatorname{rank}}
\newcommand\rk{\operatorname{rk}}

\newcommand\RE{\operatorname{Re}}

\newcommand{\sign}{\operatorname{sign}}

\newcommand\spec{{\rm spec\,}}

\newcommand\tr{\operatorname{tr}}
\newcommand\Tr{\operatorname{Tr}}

\newcommand\on{{\overline{n}}}

\newcommand\os{{\bar{s}}}

\newcommand\oz{{\overline{z}}}

\newcommand\oeta{{\overline{\eta}}}

\newcommand\olam{{\overline{\lambda}}}

\newcommand\hatB{{\widehat{B}}}

\newcommand\hatL{{\widehat{L}}}

\newcommand\tilB{{\widetilde{B}}}

\newcommand\tilM{{\widetilde{M}}}

\newcommand\tiln{{\widetilde{n}}}

\newcommand\tilDel{{\widetilde{\Delta}}}

\newcommand\tilome{{\widetilde{\ome}}}

\renewcommand\tilome{{\widetilde{\ome}}}

\renewcommand{\>}{\rangle}
\newcommand{\<}{\langle}

\theoremstyle{plain}
\newtheorem{Thm}[subsection]{Theorem}
\newtheorem{Cor}[subsection]{Corollary}
\newtheorem{Lem}[subsection]{Lemma}
\newtheorem{Prop}[subsection]{Proposition}
\newtheorem{Conjec}[subsection]{Conjecture}

\newtheorem{Def}[subsection]{Definition}%\renewcommand{\theDef}{\thesection.\arabic{Def}}

\theoremstyle{remark}

\newtheorem{Rem}[subsection]{Remark}%\renewcommand{\theRem}{\thesection.\arabic{Rem}}

\def\TeXref#1{%
        \leavevmode\vadjust{\setbox0=\hbox{{\tt
                \  {\tiny \textrm #1}}}%
        \theight=\ht0
        \advance\theight by \lineskip
        \kern -\theight \vbox to
        \theight{\rightline{\rlap{\box0}}%
        \vss}%
        }}%
\newif\ifShowLabels
\ShowLabelstrue
\newdimen\theight
\def\TeXrefEq#1{%
        \leavevmode\vadjust{\setbox0=\hbox{{\tt
                \  {\tiny \textrm #1}}}%
        \theight=\ht1
        \advance\theight by \lineskip
        \kern -\theight \vbox to
        \theight{\rightline{\rlap{\box0}}%
        \vss}%
        }}%
\newcommand{\refs}[1]{Section ~\ref{S:#1}}
\newcommand{\refss}[1]{Subsection ~\ref{SS:#1}}

\newcommand{\reft}[1]{Theorem ~\ref{T:#1}}
\newcommand{\refl}[1]{Lemma ~\ref{L:#1}}
\newcommand{\refp}[1]{Proposition ~\ref{P:#1}}
\newcommand{\refc}[1]{Corollary ~\ref{C:#1}}

\newcommand{\refd}[1]{Definition ~\ref{D:#1}}
\newcommand{\refr}[1]{Remark ~\ref{R:#1}}
\newcommand{\refe}[1]{\eqref{E:#1}}

\newenvironment{thm}[1]%
        { \begin{Thm} \label{T:#1}  \ifShowLabels \TeXref{T:#1} \fi }%
        { \end{Thm} }

\renewcommand{\th}[1]{\begin{thm}{#1}  }
\renewcommand{\eth}{\end{thm} }

\newenvironment{lemma}[1]%
        { \begin{Lem} \label{L:#1}  \ifShowLabels \TeXref{L:#1} \fi }%
        { \end{Lem} }

\newcommand{\lem}[1]{\begin{lemma}{#1} }
\newcommand{\elem}{\end{lemma}}

\newenvironment{propos}[1]%
        { \begin{Prop} \label{P:#1}  \ifShowLabels \TeXref{P:#1} \fi }%
        { \end{Prop} }

\newcommand{\prop}[1]{\begin{propos}{#1} }
\newcommand{\eprop}{\end{propos}}

\newenvironment{corol}[1]%
        { \begin{Cor} \label{C:#1}  \ifShowLabels \TeXref{C:#1} \fi }%
        { \end{Cor} }
\newcommand{\cor}[1]{\begin{corol}{#1}  }
\newcommand{\ecor}{\end{corol}}

\newenvironment{conjec}[1]%
        { \begin{Conjec} \label{Conj:#1}  \ifShowLabels \TeXref{C:#1} \fi }%
        { \end{Conjec} }
\newcommand{\conj}[1]{\begin{conjec}{#1}  }
\newcommand{\econj}{\end{conjec}}

\newenvironment{defeni}[1]%
        { \begin{Def} \label{D:#1}  \ifShowLabels \TeXref{D:#1} \fi }%
        { \end{Def} }
\newcommand{\defe}[1]{\begin{defeni}{#1}  }
\newcommand{\edefe}{\end{defeni}}

\newenvironment{remark}[1]%
        { \begin{Rem} \label{R:#1}  \ifShowLabels \TeXref{R:#1} \fi }%
        { \end{Rem} }
\newcommand{\rem}[1]{\begin{remark}{#1}}
\newcommand{\erem}{\end{remark}}

\newcommand{\eq}[1]%
        { \ifShowLabels \TeXrefEq{E:#1} \fi
           \begin{equation} \label{E:#1} }
\newcommand{\eeq}{\end{equation}}
\newcommand{\meq}[1]%
        { \ifShowLabels \TeXrefEq{E:#1} \fi
           \begin{multline} \label{E:#1} }
\newcommand{\emeq}{\end{multline}}

\newcommand{\prf}{ \begin{proof} }
\newcommand{\eprf}{ \end{proof} }
\newcommand{\Label}[1]{\label{#1}  \ifShowLabels \TeXref{#1} \fi }

\ShowLabelsfalse% comment this out if labels should be printed

\newcommand{\n}{\nabla}
\newcommand{\na}{\nabla_\alp}
\renewcommand{\tiln}{\tilde{\n}}
\renewcommand{\b}{\bullet}
\newcommand{\pa}{\text{\( \partial\)}}
\newcommand{\odd}{{\operatorname{odd}}}
\newcommand{\even}{{\operatorname{even}}}
\newcommand{\Fp}{\operatorname{F.p._{s=0}}}

\newcommand{\Arg}{\operatorname{\mathbf{Arg}}}
\newcommand{\Ph}{\operatorname{\mathbf{Ph}}}
\newcommand{\Mon}{\operatorname{Mon}}
\newcommand{\Diff}{\operatorname{Diff}}
\newcommand{\Ell}{\operatorname{Ell}}

\newcommand{\TRS}{T^{\operatorname{RS}}}
\newcommand{\TTur}{T^{\operatorname{comb}}}
\newcommand{\Tabs}{T^{\operatorname{abs}}}
\newcommand{\Detgr}{\Det_{\operatorname{gr}}}
\newcommand{\Detgrtet}{\Det_{{\operatorname{gr},\tet}}}
\newcommand{\Detgrteto}{\Det_{{\operatorname{gr},\tet_0}}}\newcommand{\Detgrtetone}{\Det_{{\operatorname{gr},\tet_1}}}

\newcommand{\symb}{\sig}
\newcommand{\tilPi}{\tilde{\Pi}}

\newcommand{\p}{\pi_1(M)}
\newcommand{\Rep}{\operatorname{Rep}(\p,\CC^n)}
\newcommand{\Repo}{\operatorname{Rep_0}(\p,\CC^n)}
\newcommand{\Reph}{\operatorname{Rep}^u(\p,\CC^n)}
\newcommand{\Repho}{\operatorname{Rep}^u_0(\p,\CC^n)}

\newcommand{\LD}{\operatorname{LDet}}

\newcommand{\LDgrtet}{\operatorname{LDet}_{\operatorname{gr},\tet}}

\newcommand{\Flat}{\operatorname{Flat}}
\newcommand{\ccomp}{C}
\begin{document}
%-------------------------------------------------------------
%-------------------------------------------------------------
\begin{flushright}
%\underline{Version: \today}
\end{flushright}

\title{Refined Analytic Torsion}
\author[Maxim Braverman]{Maxim Braverman$^\dag$}
\address{Department of Mathematics\\
        Northeastern University   \\
        Boston, MA 02115 \\
        USA
         }
\email{maximbraverman@neu.edu}
\author[Thomas Kappeler]{Thomas Kappeler$^\ddag$}
\address{Institut fur Mathematik\\
         Universitat Z\"urich\\
         Winterthurerstrasse 190\\
         CH-8057 Z\"urich\\
         Switzerland
         }
\email{tk@math.unizh.ch}
\thanks{${}^\dag$Supported in part by the NSF grant DMS-0204421.\\
\indent${}^\ddag$Supported in part by the Swiss National Science foundation, the programme SPECT, and the European Community through the FP6 Marie Curie RTN ENIGMA
(MRTN-CT-2004-5652)}
%\date{\today}
\begin{abstract}
Given an acyclic representation $\alpha$ of the fundamental group of a compact oriented odd-dimensional manifold, which is close enough to an
acyclic unitary representation, we define a refinement $T_\alpha$ of the Ray-Singer torsion associated to $\alpha$,  which can be viewed as the
analytic counterpart of the refined combinatorial torsion introduced by Turaev. $T_\alpha$ is equal to the graded determinant of the odd
signature operator up to a correction term, the {\em metric anomaly}, needed to make it independent of the choice of the Riemannian metric.

$T_\alpha$ is a holomorphic function on the space of such representations of the fundamental group. When $\alpha$ is a unitary representation, the absolute value of
$T_\alpha$ is equal to the Ray-Singer torsion and the phase of $T_\alpha$ is proportional to  the $\eta$-invariant of the odd signature operator. The fact that the
Ray-Singer torsion and the $\eta$-invariant can be combined into one holomorphic function allows to use methods of complex analysis to study both invariants. In
particular, using these methods we compute the quotient of the refined analytic torsion and Turaev's refinement of the combinatorial torsion generalizing in this way
the classical Cheeger-M\"uller theorem. As an application, we extend and improve a result of Farber about the relationship between the Farber-Turaev absolute torsion
and the $\eta$-invariant.

As part of our construction of $T_\alpha$ we prove several new results about determinants and $\eta$-invariants of non self-adjoint elliptic
operators.
\end{abstract}
\maketitle \tableofcontents

%------------------------------------------------------
%------------------------------------------------------
\section{Introduction}\Label{S:introd}

In this paper we refine the analytic torsion which has been introduced by Ray and Singer \cite{RaSi1}. In our set-up we are given a complex flat
vector bundle $E\to M$ over a closed oriented odd-dimensional manifold $M$ and we denote by $\n$ the flat connection on $M$. Whereas the
Ray-Singer torsion $\TRS(\n)$ is a positive real number, the proposed {\em refined analytic torsion} $T= T(\n)$  will be, in general, a complex
number, hence will have a nontrivial phase. The refined analytic torsion can be viewed as an analytic analogue of the refined combinatorial
torsion, introduced by Turaev \cite{Turaev86,Turaev90} and further developed by Farber and Turaev \cite{FarberTuraev99,FarberTuraev00}. Though
$T$ is {\em not} equal to the Turaev torsion in general, the two torsions are very closely related, as described in \refs{comb}.

\subsubsection*{Definition}\Label{SS:Idef} In this paper the refined analytic torsion is defined for an open set of acyclic flat connections, which contains all
acyclic Hermitian connections (see \cite{BrKappelerRATdetline}, where we extend this definition to arbitrary flat connections). If $\dim{M}\equiv 1\ (\MOD 4)$ or if
the rank of the bundle $E$ is divisible by 4, then the refined analytic torsion $T=T(\n)$ is independent of any choices. If $\dim{M}\equiv 3\ (\MOD 4)$ then $T(\n)$
depends, in addition, on the choice of a compact oriented manifold $N$, whose oriented boundary is diffeomorphic to two disjoint copies of $M$, but only up to a
factor $i^{k\cdot\rk{E}}$ \ ($k\in \ZZ$).

\subsubsection*{Relation to the $\eta$-invariant and the Ray-Singer torsion}\Label{SS:IetaRS} If the connection $\n$ is Hermitian, i.e., if there
exists a Hermitian metric on $E$ which is preserved by $\n$, then the refined analytic torsion $T$ is a complex number whose absolute value is
equal to the Ray-Singer torsion and whose phase is determined by the $\eta$-invariant of the odd signature operator. When $\n$ is not Hermitian,
the relationship between the refined analytic torsion, the Ray-Singer torsion, and the $\eta$-invariant is slightly more complicated, cf.
\refs{deponalp}.

\subsubsection*{Analytic property}\Label{SS:Ianalytic} One of the most important properties of the refined analytic torsion is that it depends, in an
appropriate sense, {\em holomorphically} on the connection $\n$. The fact that the Ray-Singer torsion and the $\eta$-invariant can be combined
into one holomorphic function allows to use methods of complex analysis to study both invariants. In particular, using these methods we
establish  a relationship between the refined analytic torsion and Turaev's refinement of the combinatorial torsion which generalizes the
classical Cheeger-M\"uller theorem about the equality between the Ray-Singer and the combinatorial torsion \cite{Cheeger79,Mu1ller78}. As an
application, we generalize and improve a result of Farber about the comparison between the sign of the Farber-Turaev absolute torsion and the
$\eta$-invariant, \cite{Farber00AT}. In fact, we compare the phase of the Turaev torsion and the $\eta$-invariant in a more general set-up.

\subsubsection*{Regularized determinant}\Label{SS:Irefdet} Our construction of the refined analytic torsion uses determinants of non self-adjoint
elliptic differential operators. In \refs{etainv} and Appendix~\ref{S:det-eta-sa} we prove several new results about these determinants which
generalize well known facts about determinants of self-adjoint differential operators. In particular, we express the determinant of a (not
necessarily self-adjoint) operator $D$ in terms of the determinant of $D^2$, the value at 0 of the $\zeta$-function of $D^2$, and the
$\eta$-invariant of $D$. Note that the $\eta$-invariant of a non-self-adjoint operator was defined and studied by Gilkey \cite{Gilkey84}. In
this paper we use a {\em sign refined} version of Gilkey's construction, cf. \refd{eta}.

%-----------------
\subsubsection*{Related works} In \cite{Turaev86,Turaev90}, Turaev constructed a refined version  of the combinatorial torsion and posed the
problem of constructing its analytic analogue, see also \cite[\S10.3]{FarberTuraev00}. More precisely, one can ask if it can be defined in terms of regularized
determinants of elliptic differential operators and, if so, whether the phase is related to the $\eta$-invariant of these differential operators. In the present paper
we show that on the open neighborhood of the set of acyclic Hermitian connection, where the proposed refined analytic torsion $T(\n)$ is defined, $T(\n)$ solves this
problem. In \cite{BrKappelerRATdetline} we extend the notion of refined  analytic torsion to the set of all flat connections.

In addition to the works of Turaev \cite{Turaev86,Turaev90} and Farber-Turaev \cite{FarberTuraev99,FarberTuraev00} on their refined
combinatorial torsion and the relation of its absolute value to the Ray-Singer torsion \cite{Turaev90,FarberTuraev00}  as well as the study of
its phase \cite{Farber00AT}, we would like to mention a recent paper of Burghelea and Haller, \cite{BurgheleaHaller_Euler}. In that paper, among
many other topics, the authors address the question if {the Ray-Singer torsion $\TRS(\n)$  can be viewed as the absolute value of a (in an
appropriate sense) holomorphic function $f(\n)$ on the space of acyclic connection $\n$}.

Burghelea and Haller gave an affirmative answer to this question and showed that
\eq{BurgheleaHaller}
   \TRS(\n) \ =\ |f_1(\n)\cdot{}f_2(\n)|,
\end{equation}
where $f_1(\n)$ is Turaev's refinement of the combinatorial torsion and $f_2(\n)$ is an explicitly calculated holomorphic function.%
\footnote{Note that, in the case when the dimension of the manifold is odd, \refe{BurgheleaHaller} is similar to our \reft{TTur-TRS}.} The
result of Burghelea and Haller is valid for manifolds of arbitrary dimension. If the dimension of the manifold is odd, the refined analytic
torsion proposed in this paper allows to obtain an analogue of \refe{BurgheleaHaller}. In contrast to \cite{BurgheleaHaller_Euler}, the
holomorphic function on the right hand side of equality \refe{BurgheleaHaller} is constructed in this paper in purely analytic terms, cf.
\reft{RaySinger}. The quotient between the Ray-Singer torsion and the absolute value of Turaev's refinement of the combinatorial torsion is
discussed in \refs{comb}.

In response to a first version of our paper, Dan Burghelea kindly brought to our attention his ongoing project with Stefan Haller \cite{BurgheleaHaller_function}
where they consider, among other things, Laplace-type operators  acting on forms obtained by replacing a Hermitian scalar product on a given complex vector bundle by
a non-degenerate symmetric bilinear form. These operators  are non self-adjoint and have complex-valued zeta-regularized determinants. Burghelea and Haller then
express the square of the Turaev torsion in terms of these determinants and some additional ingredients.

\medskip
The results of this paper were announced in \cite{BrKappelerRATshort}.

%-----------------------------------------------------------------------------------------------------
%-----------------------------------------------------------------------------------------------------

\section{Summary of the Main Results}\Label{S:summary}

Throughout this section $M$ is a closed oriented manifold of odd dimension $\dim{}M=d=2r-1$ and $E$ is a complex vector bundle over $M$ endowed with a flat connection
$\n$.
%--------------------------------------------
\subsection{The odd signature operator}\Label{SS:Ioddsign}
The refined analytic torsion is defined in terms of the odd signature operator, hence, let us begin by recalling the definition of this
operator.

Let $\Ome^\b(M,E)$ denote the space of smooth differential forms on $M$ with values in $E$ and set
\[
   \Ome^{\even}(M,E)\ = \ \bigoplus_{p=0}^{r-1}\Ome^{2p}(M,E),
\]
where $r=\frac{\dim M+1}2$. Fix a Riemannian metric $g^M$ on $M$ and let $*:\Ome^\b(M,E)\to \Ome^{d-\b}(M,E)$ denote the Hodge $*$-operator. The {\em chirality}
operator $\Gam:\Ome^\b(M,E)\to \Ome^{d-\b}(M,E)$ is then given by the formula, cf. \cite[\S3.2]{BeGeVe},
\[
    \Gam\,\ome \ := \ i^r\,(-1)^{k(k+1)/2}\,*\,\ome, \qquad \ome\in \Ome^k(M,E).
\]
The {\em odd signature operator} $B=B(\n,g^M):\Ome^\b(M,E)\to \Ome^{\b}(M,E)$ is defined by
\begin{equation} \Label{E:Ioddsign}\notag
    B \ := \ \Gam\,\n\ + \ \n\,\Gam.
\end{equation}
It leaves $\Ome^\even(M,E)$ invariant. Denote its restriction to $\Ome^\even(M,E)$ by $B_\even$. Then, for $\ome\in \Ome^{2p}(M,E)$, one has
\[
    B\,\ome \ = \ i^r\,(-1)^{p+1}\,\big(\,*\,\n-\n\,*\,\big)\,\ome\ \in\ \Ome^{d-2p-1}(M,E)\,\oplus\,\Ome^{d-2p+1}(M,E).
\]

The odd signature operator was introduced by Atiyah, Patodi, and Singer, \cite[p.~44]{APS1}, \cite[p.~405]{APS2}, and, in the more general
setting used here, by Gilkey, \cite[p.~64--65]{Gilkey84}.

The operator $B_{\even}$ is an elliptic differential operator of order one, whose leading symbol is symmetric with respect to any Hermitian metric $h^E$ on $E$.

In this paper we define the refined analytic torsion in the case when the pair $(\n,g^M)$ satisfies the following simplifying assumptions. The
general case will be addressed elsewhere.

\subsection*{Assumption~I} The connection $\n$ is acyclic, i.e.,
\[
    \IM\big(\n|_{\Ome^{k-1}(M,E)}\big) \ =\  \Ker\big(\n|_{\Ome^k(M,E)}\big),\qquad \text{for every}\quad k=0\nek d.
\]
\subsection*{Assumption~II}  $B_{\even}=B_{\even}(\n,g^M)$ is bijective.

\bigskip
Note that if $\n$ is a Hermitian connection then Assumption~I implies Assumption~II, cf. \refss{openset}. Hence, all acyclic Hermitian
connections satisfy Assumptions~I and II. By a simple continuity argument, cf. \refp{neigofh}, these two assumptions are then satisfied for all
flat connections in an open neighborhood (in $C^0$-topology, cf. \refss{openset1}) of the set of acyclic Hermitian connections.
%We denote this open neighborhood by $\Flat'(E,g^M)$.

\subsection{Graded determinant}\Label{SS:Igrdet}
Set
\eq{Igrading}
  \Ome^k_+(M,E) \ := \ \Ker\,(\n\,\Gam)\,\cap\,\Ome^k(M,E), \qquad
  \Ome^k_-(M,E) \ := \ \Ker\,(\Gam\,\n)\,\cap\,\Ome^k(M,E).
\end{equation}
Assumption~II implies that $\Ome^k(M,E)\ = \ \Ome^k_+(M,E)\oplus \Ome^k_-(M,E)$, cf. \refss{decompos}. Hence, \refe{Igrading} defines a {\em grading} on
$\Ome^k(M,E)$.

Define $\Ome^{\even}_\pm(M,E)= \bigoplus_{p=0}^{r-1}\Ome^{2p}_{\pm}(M,E)$ and let $B_{\even}^\pm$ denote the restriction of $B_{\even}$ to
$\Ome^{\even}_\pm(M,E)$. It is easy to see that $B_{\even}$ leaves the subspaces $\Ome^{\even}_\pm(M,E)$ invariant and it follows from
Assumption~II that the operators $B^\pm_{\even}:\Ome^{\even}_\pm(M,E)\to \Ome^{\even}_\pm(M,E)$ are bijective.

One of the central objects of this paper is the {\em graded determinant} of the operator $B_{\even}$. To construct it we need to choose a {\em
spectral cut} along a ray $R_{\tet}= \big\{\rho{}e^{i\tet}:0\le \rho<\infty\big\}$, where $\tet\in [-\pi,\pi)$ is an Agmon angle for
$B_{\even}$, cf. \refd{Agmon}. Since the leading symbol of $B_{\even}$ is symmetric, $B_{\even}$ admits an Agmon angle $\tet\in (-\pi,0)$. Given
such an angle $\tet$, observe that it is an Agmon angle for $B_{\even}^\pm$ as well. The graded determinant of $B_{\even}$ is the non-zero
complex number defined by the formula
\eq{Igrdeterminant}
   \Detgrtet(B_{\even}) \ := \ \frac{\Det_\tet(B_{\even}^+)}{\Det_\tet(-B_{\even}^-)}.
\end{equation}
By standard arguments, cf. \refss{det-tet}, $\Detgrtet(B_{\even})$ is independent of the choice of the Agmon angle $\tet\in (-\pi,0)$.

%-------------------------------------
\subsection{A convenient choice of the Agmon angle}\Label{SS:IAgmon}
For $\calI\subset \RR$ we denote by $L_{\calI}$ the solid angle
\eq{ILab}\notag
    L_{\calI} \ = \  \big\{\, \rho e^{i\tet}:\, 0 < \rho<\infty,\, \tet\in \calI\,  \big\}.
\end{equation}

Though many of our results are valid for any Agmon angle $\tet\in (-\pi,0)$, some of them are easier formulated if the following conditions are
satisfied:

\begin{enumerate}
\item[\textbf{(AG1)}] \ $\tet\in (-\pi/2,0)$, and
\item[\textbf{(AG2)}] \ there are no eigenvalues of the operator $B_{\even}$ in the solid angles $L_{(-\pi/2,\tet]}$ and $L_{(\pi/2,\tet+\pi]}$.
\end{enumerate}

For the sake of simplicity of exposition, we will assume that $\tet$ is chosen so that these conditions are satisfied throughout the Introduction. Since the leading
symbol of $B_{\even}$ is symmetric (with respect to an arbitrary Hermitian metric on $E$), such a choice of $\tet$ is always possible.

%------------------------------------
\subsection{Relationship with the Ray-Singer torsion and the $\eta$-invariant}\Label{SS:Ietainv}
For a pair $(\n,g^M)$ satisfying Assumptions~I and II set
\eq{IcalB}
    \xi \ = \  \xi(\n,g^M,\tet) \ := \ \frac12\,\sum_{k=0}^{d-1}\,(-1)^k\,\zet_{2\tet}' \big(\,0,\,
    {(\Gam\,\n)^2}{\big|_{\Ome^k_+(M,E)}}\,\big),
\end{equation}
where $\zet_{2\tet}'\big(\,s,\,{(\Gam\,\n)^2}{\big|_{\Ome^k_+(M,E)}}\,\big)$ is the derivative with respect to $s$ of the $\zet$-function of the operator
${(\Gam\,\n)^2}{\big|_{\Ome^k_+(M,E)}}$ corresponding to the spectral cut along the ray $R_{2\tet}$, cf. \refss{zet-det}, and $\tet$ is an Agmon angle satisfying
(AG1)-(AG2).

Let $\eta=\eta(\n,g^M)$ denote the (sign refined) $\eta$-invariant of the operator $B_{\even}(\n,g^M)$, cf. \refd{eta}. \reft{DetB-eta2} implies
that,
\eq{IDetB-eta}
    \Detgrtet(B_{\even}) \ = \ e^{\xi(\n,g^M,\tet)}\cdot e^{-i\pi\eta(\n,g^M)}.
\end{equation}
This representation of the graded determinant turns out to be very useful, e.g., in computing the metric anomaly of $\Detgrtet(B_{\even})$.

If the connection $\n$ is Hermitian, then \refe{IcalB} coincides with the well known expression for the logarithm of the Ray-Singer torsion
$\TRS= \TRS(\n)$. Hence, for a Hermitian connection $\n$ we have
\[
    \xi(\n,g^M,\tet) \ = \ \log\,\TRS(\n).
\]
If $\n$ is not Hermitian but is sufficiently close (in $C^0$-topology) to an acyclic Hermitian connection, then \reft{DetB-TRS} states that
\eq{IDetB-TRS}
    \log \TRS(\n) \ = \ \RE\,\xi(\n,g^M,\tet).
\end{equation}

Combining \refe{IDetB-TRS} and \refe{IDetB-eta}, we get
\eq{IDetB-TRSherm}
    \big|\,  \Detgrtet(B_{\even}) \,\big| \ = \ \TRS(\n)\cdot e^{\pi\IM\eta(\n,g^M)}.
\end{equation}
If $\n$ is Hermitian, then the operator $B_{\even}$ is self-adjoint (cf. \refss{openset}) and $\eta= \eta(\n,g^M)$ is real. Hence, cf.
\refc{DetB-TRS}, for the case of an acyclic Hermitian connection we obtain from \refe{IDetB-TRSherm}
\eq{IDetB-TRSherm2}\notag
    \big|\,  \Detgrtet(B_{\even}) \,\big| \ = \ \TRS(\n).
\end{equation}

%-------------------------------------------
\subsection{Metric anomaly of the graded determinant}\Label{SS:Ideponmet}
The graded determinant of the odd signature operator is not a differential invariant of the connection $\n$ since, in general, it depends on the
choice of the Riemannian metric $g^M$. Hence, we first investigate the {\em metric anomaly} of the graded determinant and then use it to
``correct" the graded determinant and construct a  differential invariant -- the refined analytic torsion.

Suppose an acyclic connection $\n$ is given. We call a Riemannian metric $g^M$ on $M$ {\em admissible for $\n$} if the operator $B_{\even}= B_{\even}(\n,g^M)$
satisfies Assumption~II of \refss{Ioddsign}. We denote the set of admissible metrics by $\calM(\n)$. The set $\calM(\n)$ might be empty. However, \refp{neigofh}
implies that admissible metrics exist for all flat connections in an open neighborhood (in $C^0$-topology) of the set of acyclic Hermitian connections.

For each admissible metric $g^M\in \calM(\n)$ choose an Agmon angle $\tet$ satisfying (AG1)-(AG2). Then the reduction of $\xi(\n,g^M,\tet)$ modulo $\pi\ZZ$ depends
neither of the choice of $\tet$ nor on the choice of $g^M\in \calM(\n)$, cf. \refp{detB2metric}.

The dependence of $\eta= \eta(\n,g^M)$ on $g^M$ has been analyzed in \cite{APS2} and \cite{Gilkey84}. In particular, it follows from the results in these papers that
(cf. \refp{deponmetr}) {\em
\begin{itemize}
\item
If $\dim M\equiv 1\ (\MOD 4)$ then the reduction of $\eta(\n,g^M)$ modulo $\ZZ$ is independent of the choice of the admissible metric $g^M$;
\item
Suppose $\dim M\equiv 3\ (\MOD 4)$ and let $N$ be a compact oriented manifold whose oriented boundary is diffeomorphic to two disjoint copies of $M$ (since $\dim{}M$
is odd, such a manifold always exists, cf. \cite{Wall60}, \cite[Th.~IV.6.5]{Rudyak_book}). Denote by $L(p)= L_N(p)$ the Hirzebruch $L$-polynomial in the Pontrjagin
forms of a Riemannian metric on $N$ which is a product near $\partial{N}= M\sqcup{M}$. Then, modulo $\ZZ$,
\[
    \eta\,-\, \frac{\rank E}2\,\int_{N}\,L(p)\qquad \text{and, hence, }\qquad \IM\eta
\]
are independent of the choice of the metric on $N$. Note that for different choices of $N$ satisfying $\partial{N}=M\sqcup{M}$ the integral $\int_{N}L(p)$ differs by
an integer.
\end{itemize}}

%-------------------------------------------
\subsection{Definition of the refined analytic torsion}\Label{SS:Irefinedtorsion}
The {\em refined analytic torsion} $T(\n)$ corresponding to an acyclic connection $\n$, satisfying $\calM(\n)\not=\emptyset$, is defined as follows: fix an admissible
Riemannian metric $g^M\in \calM(\n)$ and let  $\tet\in (-\pi,0)$ be an Agmon angle for $B_{\even}(\n,g^M)$.
\begin{enumerate}
{\em
\item
If $\dim M\equiv 1\ (\MOD 4)$ then
\[
   T(\n)\ =\ T(M,E,\n)\ :=\ \Detgrtet\big(B_{\even}(\n,g^M)\big) \ \in \ \CC\backslash{0}.
\]
\item
If $\dim M\equiv 3 (\MOD 4)$ choose a smooth compact oriented manifold $N$ whose oriented boundary is diffeomorphic to two disjoint copies of
$M$. Then
\begin{equation}\Label{E:Irefantor3}\notag
    T(\n) \ = \ T(M,E,\n,N) \ := \
        \Detgrtet(B_{\even})\cdot \exp\Big(\,i\pi\,\frac{\rank E}2\,  \int_{N}\, L(p)\,\Big) \ \in \ \CC\backslash{0}.
\end{equation}}
\end{enumerate}

If $\n$ is close enough to an acyclic Hermitian connection, then $\calM(\n)\not=\emptyset$ and, it follows from the discussion of the metric anomaly of the graded
determinant of $B_{\even}$ in \refss{Ideponmet},  {\em $T(\n)$ is independent of the choice of the admissible metric $g^M$}. Moreover, as $\Detgrtet(B_{\even})$ is
independent of the choice of the Agmon angle $\tet\in (-\pi,0)$, so is $T(\n)$. However, if $\dim M\equiv 3 (\MOD 4)$, then the refined analytic torsion {\em  does
depend on the choice of the manifold $N$}. The quotient of the refined torsions corresponding to different choices of $N$ is a complex number of the form
$i^{k\cdot\rank{E}}$ $(k\in \ZZ)$. Hence, if $\rank{}E$ is even then $T(\n)$ is well defined up to a sign, and if $\rank{}E$ is divisible by 4, then  $T(\n)$ is a
well defined complex number. (Here a quantity being well defined means that it depends only on $M$, $E$ and $\n$.)

A simple example in \refss{example} shows that {\em even when the connection $\n$ is Hermitian, the refined analytic torsion can have an
arbitrary phase}.

%-----------------------------------------
\subsection{Comparison with the Ray-Singer torsion}\Label{SS:IlocalRS}
The equality \refe{IDetB-TRSherm} implies that, if $\n$ is $C^0$-close to an acyclic Hermitian connection, then
\eq{ITRS-T}
    \log\,\frac{|T(\n)|}{\TRS(\n)} \ = \ \pi\,\IM\,\eta(\n,g^M)
\end{equation}
In particular, if $\n$ is an acyclic Hermitian connection, then
\eq{|T|}\notag
    \big|\,T(\n)\,\big| \ = \ \TRS(\n).
\end{equation}

\reft{RaySinger} provides a local expression for the right hand side of \refe{ITRS-T}. Following Farber, \cite{Farber00AT}, we denote by
$\Arg_\n$ the unique cohomology class $\Arg_\n\in H^1(M,\CC/\ZZ)$ such that for every closed curve $\gam\in M$ we have
\eq{IArg}\notag
    \det\big(\,\Mon_\n(\gam)\,\big) \ = \ \exp\big(\, 2\pi i\<\Arg_\n,[\gam]\>\,\big),
\end{equation}
where $\Mon_\n(\gam)$ denotes the monodromy of the flat connection $\n$ along the curve $\gam$ and $\<\cdot,\cdot\>$ denotes the natural pairing
 \(
    H^1(M,\CC/\ZZ)\,\times\, H_1(M,\ZZ) \to  \CC/\ZZ.
 \)
\reft{RaySinger} states that, if $\n$ is $C^0$-close to an acyclic Hermitian connection, then
\eq{IRaySinger3}
    \log\,\frac{|T(\n)|}{\TRS(\n)} \ \ = \ \ \pi\,\big\<\, [L(p)]\cup \IM\Arg_{\n},[M]\,\big\>.
\end{equation}

If \/ $\dim M\equiv 3\ (\MOD\ 4)$, then $L(p)$ has no component of degree $\dim{M}-1$ and, hence,
 \(
    |T(\n)|  =  \TRS(\n).
  \)

%-----------------------------------------------
\subsection{The refined analytic torsion as a holomorphic function on the space of representations}\Label{SS:Iholomorphic}
One of the main properties of the refined analytic torsion $T(\n)$ is that, in an appropriate sense, it depends holomorphically on the
connection. Note, however, that the space of connections is infinite dimensional and one needs to choose an appropriate notion of a holomorphic
function on such a space. A possible choice is explained in \refss{analofconnections}. As an alternative one can view the refined analytic
torsion as a holomorphic function on a finite dimensional space, which we shall now explain.

The set $\Rep$ of all $n$-dimensional complex representations of $\p$ has a natural structure of a complex algebraic variety, cf. \refss{Rep}.
Each representation $\alp\in \Rep$ gives rise to a vector bundle $E_\alp$ with a flat connection $\na$, whose monodromy is isomorphic to $\alp$,
cf. \refss{Rep}. Let $\Repo\subset \Rep$ denote the set of all representations $\alp\in \Rep$ such that the connection $\na$ is acyclic. We also
denote by $\Reph\subset \Rep$ the set of all unitary representations and set $\Repho= \Reph\cap\Repo$.

Denote by $V\subset \Repo$ the set of representations $\alp$ for which there exists a metric $g^M$ so that the odd signature operator
$B_{\even}(\n,g^M)$ is bijective (i.e., Assumption~II of \refss{Ioddsign} is satisfied). It is not difficult to show, cf. \refss{DetRep}, that
$V$ is an open neighborhood of the set $\Repho$ of acyclic unitary representations.

For every $\alp\in V$ one defines the refined analytic torsion $T_\alp:= T(\na)$. \refc{Tanalytic} states that the function $\alp\mapsto T_\alp$
is holomorphic on the open set of all non-singular points of $V$.

%----------------------------------------------------
\subsection{Comparison with Turaev's torsion}\Label{SS:ITuraev}
In \cite{Turaev86,Turaev90}, Turaev introduced a refinement $\TTur_\alp(\eps,\gro)$ of the combinatorial torsion associated to a representation
$\alp$ of $\pi_1(M)$. This refinement depends on an additional combinatorial data, denoted by $\eps$ and called the {\em Euler structure} as
well as on the {\em cohomological orientation} of $M$, i.e., on the orientation $\gro$ of the determinant line of the cohomology $H^\b(M,\RR)$
of $M$. There are two versions of the Turaev torsion -- the homological and the cohomological one. In this paper it turns out to be more
convenient to use the cohomological Turaev torsion as it is defined in Section~9.2 of \cite{FarberTuraev00}. For $\alp\in \Repo$ the
cohomological Turaev torsion $\TTur_\alp(\eps,\gro)$ is a non-vanishing complex number.

Theorem~10.2 of \cite{FarberTuraev00} computes the quotient of the Turaev and the Ray-Singer torsions. Combined with \refe{IRaySinger3} this
result leads to the following  equality (cf. \refss{T-TTur}):

Let $c(\eps)\in H_1(M,\ZZ)$ be the characteristic class of the Euler structure $\eps$, cf. \cite{Turaev90} or Section~5.2 of \cite{FarberTuraev00}.
Denote by $\hatL(p)\in H_\b(M,\ZZ)$ the Poincar\'e dual of the cohomology class $[L(p)]$ and let $\hatL_1\in H_1(M,\ZZ)$ be the component of
$\hatL(p)$ in $H_1(M,\ZZ)$. Then there exists an open neighborhood $V'\subset V$ of $\Repho$ such that for every $\alp\in V'$
\eq{IRelogTTTur}
   \left|\,\frac{T_\alp}{\TTur_\alp(\eps,\gro)}\,\right|^2 \ = \ \left|\, e^{-2\pi i\<\Arg_\alp,c(\eps)+\hatL_1\>}\,\right|,
\end{equation}
where $\Arg_\alp:= \Arg_{\na}\in H^1(M,\CC/\ZZ)$ is as in \refss{IlocalRS} and $\<\cdot,\cdot\>$ denotes the natural pairing
$H^1(M,\CC/\ZZ)\times{}H_1(M,\ZZ) \to \CC/\ZZ$.

Let $\Sig$ denote the set of singular points of the complex analytic set $\Rep$. By \refc{Tanalytic}, the refined analytic torsion $T_\alp$ is a
non-vanishing holomorphic function of $\alp\in V\backslash\Sig$. By the very construction \cite{Turaev86,Turaev90,FarberTuraev00} the Turaev
torsion is a non-vanishing holomorphic function of $\alp\in \Repo$. Hence, $(T_\alp/{\TTur_\alp})^2$ is a holomorphic function on
$V'\backslash\Sig$. By construction of the cohomology class $\Arg_\alp$, for every homology class $z\in H_1(M,\ZZ)$, the expression
 \(
   e^{2\pi\,i\,\<\Arg_\alp,z\>}
 \)
is a holomorphic function on $\Rep$.

If the absolute values of two non-vanishing holomorphic functions are equal on a connected open set then the functions must be equal up to a
factor $\phi\in \CC$ with $|\phi|=1$. This observation and \refe{IRelogTTTur} lead to the following generalization of the Cheeger-M\"uller
theorem, cf. \reft{TTTur}: \
\th{ITTTur}
For each connected component $\ccomp$ of \/ $V'$, there exists a constant $\phi_\ccomp= \phi_\ccomp(\eps,\gro)\in \RR$, depending on $\eps$ and
$\gro$, such that
\eq{IlogTTTur}
    \frac{T_\alp}{\TTur_\alp(\eps,\gro)} \ = \ e^{i\phi_\ccomp}\,{}_{{}_\ccomp}\sqrt{e^{-2\pi i\<\Arg_\alp,c(\eps)+\hatL_1\>}},
\end{equation}
where ${{}_{{}_\ccomp}}\sqrt{\ }$ is the analytic square root defined in \refe{sqrtArg}.
\eth

%------------------------------------------------
\subsection{Application: Phase of the Turaev torsion of a unitary representation}\Label{SS:Farber}
We denote the phase of a complex number $z$ by $\Ph(z)\in [0,2\pi)$ so that $z= |z|e^{i\Ph(z)}$. Set $\eta_\alp:= \eta(\na,g^M)$.

Suppose  $\alp_1,\, \alp_2\in \Repho$ are unitary representations which lie in the same connected component of\, $V'$, where $V'\subset V$ is
the open neighborhood of $\Repho$ defined in \refss{ITuraev}. As an application of \refe{IlogTTTur} one obtains, cf. \reft{extofFarber}, that,
modulo $\pi\,\ZZ$,
\meq{IextofFarber2}
    \Ph(\TTur_{\alp_1}(\eps,\gro))  + \pi\, \eta_{\alp_1}  -
    \pi\,\big\<\Arg_{\alp_1},c(\eps)+\hatL_1\big\>
    \\ \equiv \ \Ph(\TTur_{\alp_2}(\eps,\gro))  + \pi\, \eta_{\alp_2} -
    \pi\, \big\<\Arg_{\alp_2},c(\eps)+\hatL_1\big\>.
\end{multline}

%Moreover, if an additional condition \refe{4z} is satisfied then \refe{IextofFarber2} holds modulo $2\pi\,\ZZ$.

%-------------------------
\subsection{Sign of the absolute torsion and a theorem of Farber}\Label{SS:IsignofAT}
Suppose that the Stiefel-Whitney class
 \(
    w_{d-1}(M)\ \in \ H^{d-1}(M,\ZZ_2)
 \)
vanishes (this is always the case when $\dim{M}\equiv 3\ (\MOD\,4)$, cf. \cite{Massey60}). Then one can choose an Euler structure $\eps$ such
that $c(\eps)=0$, cf. \cite[\S3.2]{FarberTuraev99}. Assume, in addition, that the first Stiefel-Whitney class $w_1(E_\alp)$, viewed as a
homomorphism $H_1(M,\ZZ)\to \ZZ_2$, vanishes on the 2-torsion subgroup of $H_1(M,\ZZ)$. In this case there is also a canonical choice of the
cohomological orientation $\gro$, cf.  \cite[\S3.3]{FarberTuraev99}. Then the Turaev torsions $\TTur_\alp(\eps,\gro)$ corresponding to different
choices of $\eps$ with $c(\eps)=0$ and the canonically chosen $\gro$ will be the same.

If the above assumptions on $w_{d-1}(M)$ and $w_1(E_\alp)$ are satisfied, then the number
\[
    \Tabs_\alp \ := \ \TTur_\alp(\eps,\gro)\ \in\ \CC, \qquad (\,c(\eps)\,=\,0\,),
\]
is canonically defined, i.e., is independent of any choices. It was introduced by Farber and Turaev, \cite{FarberTuraev99}, who called it the
{\em absolute torsion}. If $\alp\in \Repho$, then $\Tabs_\alp\in \RR$, cf. Theorem~3.8 of \cite{FarberTuraev99}. In \refss{signofAT} we show
that, under  the above assumptions, \refe{IextofFarber2} implies that {\em if $\alp_1,\, \alp_2\in \Repho$ are unitary representations which lie
in the same connected component of\, $V'$ then the following statements hold:
\begin{enumerate}
\item
in the case\/ $\dim{M}\equiv 3$ ($\MOD\,4$)
\[
    \sign\big(\,\Tabs_{\alp_1})\cdot e^{i\pi\eta_{\alp_1}} \ = \ \sign\big(\,\Tabs_{\alp_2})\cdot e^{i\pi\eta_{\alp_2}}.
\]
\item
in the case\/ $\dim{M}\equiv 1$ ($\MOD\,4$)
\[
    \sign\big(\,\Tabs_{\alp_1})\cdot e^{i\pi\big(\, \eta_{\alp_1}- \< [L(p)]\cup \Arg_{\alp_1},[M]\>\,\big)}
    \ = \ \sign\big(\,\Tabs_{\alp_2})\cdot e^{i\pi\big(\, \eta_{\alp_2}- \< [L(p)]\cup \Arg_{\alp_2},[M]\>\,\big)}.
\]
\end{enumerate}}

For the special case when there is a real analytic path $\alp_t$ of unitary representations connecting $\alp_1$ and $\alp_2$ such that the
twisted deRham complex \refe{deRham} is acyclic for all but finitely many values of $t$, \reft{Farber} was established by Farber, using a completely different method,%
\footnote{Note that Farber's definition of the $\eta$-invariant differs from ours by a factor of 2 and also that the sign in front of $\<
[L(p)]\cup \Arg_{\alp_1},[M]\>$ is wrong in \cite{Farber00AT}.} see \cite{Farber00AT}, Theorems~2.1 and 3.1.

%------------------------------------------------------
%------------------------------------------------------
\section{Preliminaries on Determinants of Elliptic Operators}\Label{S:determinants}

In this section we briefly review the main facts about the $\zet$-regularized determinants of elliptic operators. At the end of the section (cf.
\refss{grdeterminant}) we define a sign-refined version of the graded determinant -- a notion, which plays a central role in this paper.

%---------------------------
\subsection{Setting}\Label{SS:setting0}
Throughout this paper let $E$ be a complex vector bundle over a smooth compact manifold $M$ and let $D:C^\infty(M,E)\to C^\infty(M,E)$ be an
elliptic differential operator of order $m\ge 1$. Denote by $\symb(D)$ the leading symbol of $D$.

%-------------------------
\subsection{Choice of an angle}\Label{SS:agmonangle}
Our aim is to define the $\zeta$-function and the determinant of $D$. For this we will need to define the complex
powers of $D$. As usual, to define complex powers we need to choose a {\em spectral cut} in the complex plane. We will
restrict ourselves to the simplest spectral cuts given by a ray
\eq{Rtet}
    R_\tet \ = \ \big\{\, \rho e^{i\tet}:\, 0\le \rho<\infty\, \big\},
    \qquad 0\le \tet < 2\pi.
\end{equation}
Consequently, we have to choose an angle $\tet\in [0,2\pi)$.

%--------------
\defe{prangle}
The angle $\tet$ is a {\em principal angle} for an elliptic operator $D$ if
\eq{prangle}\notag
    \spec\big(\, \symb(D)(x,\xi)\, \big)\ \cap\ R_\tet \ = \ \emptyset,
    \qquad \text{for all}\quad x\in M,\ \xi\in T^*_xM\backslash\{0\}.
\end{equation}
\edefe

If $\calI\subset \RR$ we denote by $L_{\calI}$ the solid angle
\eq{Lab}\notag
    L_{\calI} \ = \  \big\{\, \rho e^{i\tet}:\, 0 < \rho<\infty,\, \tet\in \calI\,  \big\}.
\end{equation}

\defe{Agmon}
The angle $\tet$ is an {\em Agmon angle} for an elliptic operator $D$ if it is a principal angle for $D$ and there exists $\eps>0$ such that
\[
    \spec(D) \, \cap L_{[\tet-\eps,\tet+\eps]} \ = \ \emptyset.
\]
\edefe

%-------------------------
\subsection{$\zet$-function and determinant}\Label{SS:zet-det}
Let $\tet$ be an Agmon angle%
\footnote{The existence of an Agmon angle is an additional assumption on $D$, though a very mild one. In particular, if $D$ possesses a principal angle
 it also possesses an Agmon angle, cf. the discussion in \refss{det-tet}.}
for $D$. Assume, in addition, that $D$ is injective. In this case, the $\zeta$-function $\zeta_\tet(s,D)$ of $D$ is defined as follows.

Since $D$ is invertible, there exists a small number $\rho_0>0$ such that
\eq{rho0}\notag
    \spec(D) \, \cap \, \big\{\, z\in \CC;\, |z|<2\rho_0\, \big\} \ = \ \emptyset.
\end{equation}
Define the contour $\Gam= \Gam_{\tet,\rho_0}\subset \CC$ consisting of three curves $\Gam= \Gam_1\cup \Gam_2\cup
\Gam_3$, where
\begin{gather}\Label{E:Gamtetrho}\notag
    \Gam_1 \ = \ \big\{\, \rho e^{i\tet}:\, \infty >\rho\ge \rho_0\, \big\},
    \quad
    \Gam_2 \ = \ \big\{\, \rho_0 e^{i\alp}:\, \tet< \alp<\tet+2\pi\, \big\},\\
    \quad
    \Gam_3 \ = \ \big\{\, \rho e^{i(\tet+2\pi)}:\, \rho_0\le \rho<\infty\, \big\}.
\end{gather}

For $\RE s> \frac{\dim M}m$, the operator
\eq{Ds}
    D_\tet^{-s} \ = \ \frac{i}{2\pi}\, \int_{\Gam_{\tet,\rho_0}}\, \lam^{-s}(D-\lam)^{-1}\, d\lam
\end{equation}
is a pseudo-differential operator with continuous kernel $D_\tet^{-s}(x,y)$, cf. \cite{Seeley67, ShubinPDObook}. In particular, the operator
$D_\tet^{-s}$ is of trace class. When the angle $\tet$ is fixed we will often write $D^{-s}$ for $D^{-s}_\tet$.

We define
\eq{zeta}
    \zeta_{\tet}(s,D) \ = \ \Tr D_\tet^{-s} \ = \ \int_M\, \tr D_\tet^{-s}(x,x)\, dx,
    \qquad \RE s> \frac{\dim M}m.
\end{equation}
It was shown by Seeley \cite{Seeley67} (see also \cite{ShubinPDObook}) that $\zeta_{\tet}(s,D)$ has a meromorphic
extension to the whole complex plane and that 0 is a regular value of $\zeta_{\tet}(s,D)$.

More generally, let $Q$ be a pseudo-differential operator of order $q$. We set
\eq{zetaQ}
    \zeta_{\tet}(s,Q,D) \ = \  \Tr\, Q\,D_\tet^{-s},
    \qquad \RE s> (q+\dim M)/m.
\end{equation}
This function also has a meromorphic extension to the whole complex plane, cf. \cite[\S3.22]{Wodzicki87}, \cite[Th.~2.7]{Grubbseeley95}, and \cite{Guillemin85}.
Moreover, if $Q$ is a $0$-th order pseudo-differential projection, i.e. a 0-th order pseudo-differential operator satisfying $Q^2=Q$, then by \cite[\S7]{Wodzicki84},
\cite{Wodzicki87} (see also \cite{BruningLesch99,Ponge-asymetry}), $\zeta_{\tet}(s,Q,D)$ is regular at 0.

%Finally, we are able to define the $\zeta$-regularized determinant of $D$.
\defe{determinant}
The  $\zeta$-regularized determinant of $D$ is defined by the formula
\eq{logdet}
    \Det_\tet (D) \ := \ \exp\,\left(\,-\frac{d}{ds}\big|_{s=0}\zeta_{\tet}(s,D)\,\right).
\end{equation}
\edefe
Roughly speaking, \refe{logdet} says that the logarithm $\log\Det_\tet(D)$ of the determinant of $D$ is equal to $-\zeta_{\tet}'(0,D)$. However, the
logarithm is a multivalued function. Hence,  $\log\Det_\tet(D)$ is defined only up to a multiple of $2\pi{i}$, while $-\zeta_{\tet}'(0,D)$ is a well
defined complex number. We denote by $\LD_\tet(D)$ the particular value of the logarithm of the determinant such that
\eq{logdet2}
    \LD_\tet(D) \ = \ -\zeta_{\tet}'(0,D).
\end{equation}
Let us emphasis that the equality \refe{logdet2} is the definition of the number $\LD_\tet(D)$.

We will need the following generalization of \refd{determinant}.
\defe{logdetV}
Suppose $Q$ is a $0$-th order pseudo-differential projection commuting with $D$. Then $V:= \IM{}Q$ is a $D$ invariant subspace of $C^\infty(M,E)$.
The $\zet$-regularized determinant of the restriction $D{|_V}$ of $D$ to $V$ is defined by the formula
\eq{logdetV}
    \Det_\tet (D{|_V}) \ := \ e^{\LD_\tet(D|_V)},
\end{equation}
where
\eq{logdetV2}
    \LD_\tet(D|_V) \ = \ -\frac{d}{ds}\big|_{s=0}\zeta_{\tet}(s,Q,D).
\end{equation}
\edefe

%%%
\rem{DetDV}
From the representation of $\zet_\tet(s,Q,D)$ for $\RE{s}>\frac{\dim{}M}m$ by the eigenvalues of $D|_V$, cf. \refe{zet-zet} below, it follows that
the right hand side of \refe{logdetV2} is independent of $Q$ except through $\IM(Q)$. This justifies the notation  $\LD_\tet(D|_V)$. However, we need
to know that $V$ is the image of a 0-th order pseudo-differential projection $Q$ to ensure that $\zeta_{\tet}(s,D)$ has a meromorphic extension to
the whole $s$-plane with  $s=0$ being a regular point.
\erem

%-------------------------------------------
\subsection{Case of a self-adjoint leading symbol}\Label{SS:det-sa}
Let us assume now that
\eq{close2sa}
    \symb(D)^*(x,\xi)  \ = \ \symb(D)(x,\xi), \qquad (x,\xi)\in T^*M,
\end{equation}
where $\symb(D)^*(x,\xi)$ denotes the adjoint of the linear operator $\symb(D)(x,\xi)$ with respect to some fixed scalar product on the fibers on $E$. This assumption
implies that $D$ can be written as a sum $D=D'+A$ where $D'$ is a self-adjoint differential operator of order $m$ and $A$ is a differential operator of order smaller
than $m$.

Though the operator $D$ is not self-adjoint in general, the assumption \refe{close2sa} guarantees that it has nice spectral properties. More
precisely, cf. \cite[\S{}I.6]{Markus88}, the space $L^2(M,E)$ of square integrable sections of $E$ is the closure of the algebraic direct sum of
finite dimensional $D$-invariant subspaces
\eq{L=sumLam}
    L^2(M,E) \ = \ \overline{\bigoplus_{k\ge 1} \Lam_k}
\end{equation}
such that the restriction of $D$ to $\Lam_k$ has a unique eigenvalue $\lam_k$ and $\lim_{k\to\infty}|\lam_k|=\infty$. In general, the sum
\refe{L=sumLam} is not a sum of mutually orthogonal subspaces.

The space $\Lam_k$ are called the {\em space of root vectors of $D$ with eigenvalue $\lam_k$}. We call the dimension of the space $\Lam_k$ the {\em
(algebraic) multiplicity} of the eigenvalue $\lam_k$ and we denote it by $m_k$.

Assume now that $\tet$ is an Agmon angle for $D$. As, for $\RE{s}>\frac{\dim{M}}m$ the operator $D_\tet^{-s}$ is of trace class, we conclude by
Lidskii's theorem, \cite{Lidskii59}, \cite[Ch.~XI]{Retherford93}, that the $\zet$-function \refe{zeta} is equal to the sum (including the algebraic
multiplicities)  of the eigenvalues of $D_\tet^{-s}$. Hence,
\eq{zeta-lam}
    \zet_{\tet}(s,D) \ = \ \sum_{k=1}^\infty\, m_k\ (\lam_k)_\tet^{-s}
    \ = \ \sum_{k=1}^\infty\, m_k\ e^{-s \log_\tet \lam_k},
\end{equation}
where $\log_\tet(\lam)$ denotes the branch of the logarithm in $\CC\backslash{R_\tet}$ with $\tet<\IM\log_\tet(\lam)<\tet+2\pi$.

%-------------------------------------
\subsection{Dependence of the determinant on the angle}\Label{SS:det-tet}
Assume now that  $\tet$ is only a principal angle for $D$. Then, cf. \cite{Seeley67, ShubinPDObook}, there exists
$\eps>0$ such that $\spec(D)\cap L_{[\tet-\eps,\tet+\eps]}$ is finite and $\spec(\symb(D))\cap
L_{[\tet-\eps,\tet+\eps]}=\emptyset$. Thus we can choose an Agmon angle $\tet'\in (\tet-\eps,\tet+\eps)$ for $D$. In
this subsection we show that $\Det_{\tet'}(D)$ is independent of the choice of this angle $\tet'$. For simplicity, we
will restrict ourselves to the case when $D$ has a self-adjoint leading  symbol.

Let $\tet''> \tet'$  be another Agmon angle for $D$ in $(\tet-\eps,\tet+\eps)$. Then there are only finitely many eigenvalues $\lam_{r_1}\nek
\lam_{r_l}$ of $D$ in the solid angle $L_{[\tet',\tet'']}$. We have
\eq{tet-tet}
    \log_{\tet''}\lam_k \ = \
    \begin{cases}
        \log_{\tet'}\lam_k, \quad&\text{if}\ \ k\not\in \{r_1\nek r_l\};\\
        \log_{\tet'}\lam_k+2\pi i, \quad&\text{if}\ \ k\in \{r_1\nek r_l\}.
    \end{cases}
\end{equation}
Hence
\eq{zet-zet}
    \zet_{\tet'}'(0,D) \ - \ \zet_{\tet''}'(0,D) \ = \
    \frac{d}{ds}\big|_{s=0}\,
     \sum_{i=1}^l\, m_{r_i}\,  e^{-s\log_{\tet'}(\lam_{r_i})}(1-e^{-2\pi i s})
    \ = \ 2\pi i\, \sum_{i=1}^l\, m_{r_i}
\end{equation}
and, by \refd{determinant},
\eq{det-det}
        \Det_{\tet''} (D) \ = \ \Det_{\tet'}(D).
\end{equation}

Note that the equality \refe{det-det} holds only because both angles $\tet'$ and $\tet''$ are close to a given principal angle $\tet$ so that the
intersection \/ $\spec(D)\cap L_{[\tet',\tet'']}$\, is finite. If there are infinitely many eigenvalues of $D$ in the solid angle
$L_{[\tet',\tet'']}$ then $\Det_{\tet''}(D)$ and $\Det_{\tet'}(D)$ might be different.

%-------------------------------------------
\subsection{Graded determinant}\Label{SS:grdeterminant}
Let $D:C^\infty(M,E)\to C^\infty(M,E)$ be a differential operator with a self-adjoint leading symbol. Suppose that $Q_j:C^\infty(M,E)\to C^\infty(M,E)$ ($j=0\nek d$)
are $0$-th order pseudo-differential projections commuting with $D$. Set $V_j:= \IM{}Q_j$ and assume that
\[
    C^\infty(M,E)= \bigoplus_{j=0}^dV_j.
\]
%---------
\defe{grdeterminantV}
Assume that $D$ is injective and that $\tet\in [0,2\pi)$ is an Agmon angle for the operator $(-1)^jD|_{V_j}$, for every $j=0\nek d$. The {\em graded determinant}
$\Detgrtet(D)$ of $D$ (with respect to the grading defined by the pseudo-differential projections $Q_j$) is defined by the formula
\eq{grdeterminantV}
    \Detgrtet(D) \ := \ e^{\LDgrtet(D)}.
\end{equation}
where
\eq{grldeterminantV}
    \LDgrtet(D) \ := \ \sum_{j=0}^d\, (-1)^j\,\LD_\tet\,\big(\, (-1)^jD|_{V_j}\,\big).
\end{equation}
\edefe

The following is an important example of the above situation: Let $E=\bigoplus_{j=0}^d E_j$ be a {\em graded} vector bundle over $M$. Suppose that for each $j=0\nek
d$, there is an injective elliptic differential operator
\[
    D_j:\,C^\infty(M,E_j) \ \longrightarrow \  C^\infty(M,E_j),
\]
such that $\tet\in [0,2\pi)$ is an Agmon angle for $(-1)^jD_j$ for all $j=0\nek d$. We denote by
\eq{D=bigoplus}
    D \ = \ \bigoplus_{j=0}^d\,D_j:\, C^\infty(M,E)\ \longrightarrow \  C^\infty(M,E)
\end{equation}
the direct sum of the operators $D_j$. Then \refe{grldeterminantV} reduces to
\eq{grldeterminant}
    \LDgrtet(D) \ := \ \sum_{j=0}^d\, (-1)^j\,\LD_\tet\,\big(\, (-1)^jD_j\,\big).
\end{equation}

%------------------------------------------------------
%------------------------------------------------------
\section{The $\eta$-invariant of a non Self-Adjoint Operator and the Determinant}\Label{S:etainv}

It is well known, cf. \cite{Singer85Dirac,Wojciechowski99}, that the phase of the determinant of a self-adjoint elliptic differential operator $D$
can be expressed in terms of the $\eta$-invariant of $D$ and the $\zet$-function of $D^2$. In this section we extend this result to non self-adjoint
operators.

Throughout this section we use the notation introduced in \refs{determinants} and assume that $D:C^\infty(M,E)\to
C^\infty(M,E)$ is an elliptic differential operator of order $m$ with self-adjoint leading symbol, cf. \refss{det-sa}.
We also assume that $0$ is not in the spectrum of $D$.

%--------------------------
\subsection{$\eta$-invariant}\Label{SS:etainv}
First, we recall the definition of the $\eta$-function of $D$ for a non-self-adjoint operator, cf. \cite{Gilkey84}.
\defe{eta}
Let $\tet$ be an Agmon angle for $D$, cf. \refd{Agmon}. Using the spectral decomposition of $D$ defined in \refss{det-sa}, we define the
$\eta$-function of $D$ by the formula
\eq{eta}
    \eta_{\tet}(s,D) \ = \
    \sum_{\RE\lam_k>0}\, m_k\ (\lam_k)_\tet^{-s}  \ - \
    \sum_{\RE\lam_k<0}\, m_k\ (-\lam_k)_\tet^{-s}
\end{equation}
\edefe
Note that, by definition, the purely imaginary eigenvalues of $D$ do not contribute to $\eta_\tet(s,D)$.

It was shown by Gilkey, \cite{Gilkey84}, that $\eta_\tet(s,D)$ has a meromorphic extension to the whole complex plane $\CC$ with isolated simple
poles, and that it is regular at $0$. Moreover, the number $\eta_\tet(0,D)$ is independent of the Agmon angle $\tet$.

Since the leading symbol of $D$ is self-adjoint, the angles $\pm\pi/2$ are principal angles for $D$, cf.
\refd{prangle}. In particular, there are at most finitely many eigenvalues of $D$ on the imaginary axis.

Let $m_+$ (respectively, $m_-$) denote the number of eigenvalues (counted with their algebraic multiplicities, cf.
\refss{det-sa}) of $D$ on the positive (respectively, negative) part of the imaginary axis.

%------
\defe{etainv}
The $\eta$-invariant $\eta(D)$ of $D$ is defined by the formula
\eq{etainv}
    \eta(D) \ = \ \frac{\eta_\tet(0,D)+m_+-m_-}2.
\end{equation}
In view of \refe{tet-tet}, $\eta(D)$ is independent of the angle $\tet$.
\edefe

Let $D(t)$ be a smooth 1-parameter family of operators. Then $\eta(D(t))$ is not smooth but may have integer jumps when eigenvalues cross the
imaginary axis. Because of this, the $\eta$-invariant is usually considered modulo integers. However, in this paper we will be interested in the
number $e^{i\pi\eta(D)}$, which changes its sign when $\eta(D)$ is changed by an odd integer. Hence, we will consider the $\eta$-invariant as a
complex number.

%------------
\rem{etainv}
Note that our definition of $\eta(D)$ is slightly different from the one suggested by Gilkey in \cite{Gilkey84}. In fact, in our notation, Gilkey's definition is
$\eta(D)+m_-$. Hence, reduced modulo integers the two definitions coincide. However, the number $e^{i\pi\eta(D)}$ will be multiplied by $(-1)^{m_-}$ if we replace one
definition by the other. In this sense, \refd{etainv} can be viewed as a {\em sign refinement} of the definition given in \cite{Gilkey84}.
\erem

%----------------------------
\subsection{Relationship between the $\eta$-invariant and the determinant}\Label{SS:det-eta}
Since the leading symbol of $D$ is self-adjoint, the angles $\pm\pi/2$ are principal for $D$. Hence, cf. \refss{det-tet}, there exists an Agmon
angle $\tet\in (-\pi/2,0)$ such that there are no eigenvalues of $D$ in the solid angles $L_{(-\pi/2,\tet]}$ and $L_{(\pi/2,\tet+\pi]}$. Then
$2\tet$ is an Agmon angle for the operator $D^2$.

%-------
\th{det-eta}
Let $D:C^\infty(M,E)\to C^\infty(M,E)$ be a bijective elliptic differential operator of order $m$ with self-adjoint leading symbol. Let $\tet\in (-\pi/2,0)$ be an
Agmon angle for $D$ such that there are no eigenvalues of $D$ in the solid angles $L_{(-\pi/2,\tet]}$ and $L_{(\pi/2,\tet+\pi]}$ (hence, there are no eigenvalues of
$D^2$ in the solid angle $L_{(-\pi,2\tet]}$). Then \footnote{Recall that we denote by $\LD_{\tet}(D)$ the particular branch of the logarithm of the determinant of $D$
defined by formula \refe{logdet2}.}
\eq{det-eta}
   \LD_{\tet}(D) \ = \ \frac12\,\LD_{2\tet}(D^2)  \ - \
   i\pi\,\Big(\,\eta(D)- \frac12\zet_{2\tet}(0,D^2)\,\Big).
\end{equation}
In particular,
\eq{det-eta2}
 \Det_{\tet}(D) \ = \ e^{-\frac12\,\zet_{2\tet}'(0,D^2)}\cdot e^{-i\pi\,(\,\eta(D)- \frac12\zet_{2\tet}(0,D^2)\,)}.
\end{equation}
\eth

%---------------------
\rem{det-eta}
a. \ Let $\tet$ be as in \reft{det-eta} and suppose that $\tet'\in (-\pi,0)$ is another angle such that both $\tet'$ and $\tet'+\pi$ are Agmon angles for $D$. Then,
by \refe{zet-zet} and \refe{det-det},
\eq{tet-tet'}
  \begin{aligned}
    \Det_{\tet'}(D) \ &= \ \Det_{\tet}(D), \\ \zet_{2\tet}'(0,D^2)\ &\equiv \ \zet_{2\tet'}'(0,D^2) \qquad \MOD\, 2\pi i.
  \end{aligned}
\end{equation}
In particular,
\eq{tet-tet'2}
    e^{-\frac12\,\zet_{2\tet'}'(0,D^2)} \ = \ \pm\, e^{-\frac12\,\zet_{2\tet}'(0,D^2)}.
\end{equation}

Clearly, $\zet_{\tet_1}(0,D^2)= \zet_{\tet_2}(0,D^2)$ if there are finitely many eigenvalues of $D^2$ in the solid angle $L_{[\tet_1,\tet_2]}$. Hence,
$\zet_{2\tet}(0,D^2)= \zet_{2\tet'}(0,D^2)$. We then conclude from \refe{det-eta2}, \refe{tet-tet'}, and \refe{tet-tet'2} that
\eq{set-etatet'}
 \Det_{\tet'}(D) \ = \ \pm\,e^{-\frac12\,\zet_{2\tet'}'(0,D^2)}\cdot e^{-i\pi\,(\,\eta(D)- \frac{\zet_{2\tet'}(0,D^2)}2\,)}.
\end{equation}
In other words, for \refe{det-eta2} to hold we need the precise assumption on $\tet$ which are specified in \reft{det-eta}. But ``up to a sign" it holds for every
spectral cut in the lower half plane.

b. \ If instead of the spectral cut $R_{\tet}$ in the lower half-plane we use the spectral cut $R_{\tet+\pi}$ in the upper half-plane we will
get a similar formula
\eq{det-eta1}
   \LD_{\tet+\pi}(D) \ = \ \frac12\,\LD_{2\tet}(D^2)  \ + \
   i\pi\,\Big(\,\eta(D)- \frac12\zet_{2\tet}(0,D^2)\,\Big),
\end{equation}
whose proof is a verbatim repetition of the proof of \refe{det-eta}, cf. below.

c. \ {\em If the dimension of $M$ is odd}, then  the $\zet$-function of an elliptic differential operator of even order vanishes at 0, cf.
\cite{Seeley67}. In particular, $\zet_{2\tet}(0,D^2)=0$. Hence, \refe{det-eta} simplifies to
\eq{det-eta-odd}
   \LD_{\tet}(D) \ = \ \frac12\,\LD_{2\tet}(D^2)  \ - \ i\pi\,\eta(D).
\end{equation}
\erem

%-----------------------------------
\prf
Let $\Pi_+$ and $\Pi_-$ denote the spectral projections of $D$ corresponding to the solid angles $L_{(-\pi/2,\pi/2)}$ and $L_{(\pi/2,3\pi/2)}$ respectively. Let $P_+$
and $P_-$ denote the spectral projections of $D$ corresponding to the rays $R_{\pi/2}$ and $R_{-\pi/2}$ respectively (here we use the notation introduced in
\refe{Rtet}). Set $\tilPi_\pm= \Pi_\pm+P_\pm$. Since $D$ is injective
 \(
    \Id \ = \ \tilPi_+ \ + \  \tilPi_-.
 \)
Clearly
\[
 \begin{aligned}
    \zet_{\tet}(s,D) \ = \ \Tr\,&\big[\, \tilPi_+\,D_{\tet}^{-s}\,\big] \ + \ e^{-i\pi s}\, \Tr\, \big[\tilPi_-\,(-D)_{\tet}^{-s}\,\big];\\
    \zet_{2\tet}(s/2,D^2) \  =& \ \Tr\,\big[\, \tilPi_+\,D_{\tet}^{-s}\,\big] \ + \  \Tr\, \big[\tilPi_-\,(-D)_{\tet}^{-s}\,\big].
 \end{aligned}
\]
Hence,  using the notation introduced in \refe{zetaQ}, we obtain
\eq{zet-Pizet}
 \begin{aligned}
    \zet_{\tet}(s,D) \ = \
    \zet_{\tet}&(s,\tilPi_+,D) \ + \ e^{-i\pi s}\, \zet_{\tet}(s,\tilPi_-,-D);\\
    \zet_{2\tet}(s/2,D^2) \  =& \
    \zet_{\tet}(s,\tilPi_+,D) \ + \  \zet_{\tet}(s,\tilPi_-,-D).
 \end{aligned}
\end{equation}
As, by assumption, the solid angles $L_{(-\pi/2,\tet]}$ and $L_{(\pi/2,\tet+\pi]}$ do not contain eigenvalues of $D$, it follows that
\eq{eta-Pizet}
    \eta(s,D) \ = \ \Tr\,\big[\, \Pi_+ D_{\tet}^{-s}\,\big] \ - \ \Tr \,\big[\, \Pi_- (-D)_{\tet}^{-s}\,\big]
    \ = \
    \zet_{\tet}(s,\Pi_+,D) \ - \ \zet_{\tet}(s,\Pi_-,-D).
\end{equation}

Recall that  the projectors $P_\pm$ have finite rank, which we denoted by $m_\pm$, cf. \refss{etainv}. Hence,
\[
    \zet_{\tet}(0,P_\pm,\pm D) \ = \ \RANK P_\pm  \ = \ m_\pm.
\]
Combining this equality with \refe{eta-Pizet}, and using \refe{etainv}, we obtain
\eq{etaD-Pizet}
    \eta(D) \ = \ \ \frac{\zet_{\tet}(0,\tilPi_+,D) \ - \ \zet_{\tet}(0,\tilPi_-,-D)}2.
\end{equation}

From \refe{zet-Pizet} and \refe{etaD-Pizet}, we get
\begin{multline}\Label{E:zet'=2}
    \zet_{\tet}'(0,D) \ = \
    \zet_{\tet}'(0,\tilPi_+,D) \ + \ \zet_{\tet}'(0,\tilPi_-,-D) \ - \ i\pi\,\zet_{\tet}(0,\tilPi_-,-D)
    \\ = \
    \frac12\,\zet_{2\tet}'(0,D^2) \ - \ i\pi\,\Big(\,
     \frac{\zet_{\tet}(0,\tilPi_+,D)+\zet_{\tet}(0,\tilPi_-,-D)}2 - \frac{\zet_{\tet}(0,\tilPi_+,D)-\zet_{\tet}(0,\tilPi_-,-D)}2\,\Big)
    \\ = \
    \frac12\,\zet_{2\tet}'(0,D^2) \ - \ i\pi\,\Big(\,\frac12\zet_{2\tet}(0,D^2)-\eta(D)\,\Big).
\end{multline}
\eprf

%---------------------------------
\subsection{Determinant of a self-adjoint operator and the $\eta$-invariant}\Label{SS:symmetric}
If the operator $D$ is, in addition, self-adjoint, then $\eta(D)$ and $\zet_{2\tet}(0,D^2)$ are real and the number $\Det_{2\tet}(D^2)$ is
positive, cf. \refc{detsa} in \refs{det-eta-sa}. Hence, formula \refe{det-eta2} leads to
\begin{align}
    |\Det_{\tet}(D)| \ &= \ \sqrt{\Det_{2\tet}(D^2)}, \Label{E:det-eta-sa}\\
    \Ph\big( \Det_{\tet}(D)\big) \ &\equiv \ -\pi\,\Big(\,\eta(D)- \frac12\zet_{2\tet}(0,D^2)\,\Big),
    \qquad \mod\, 2\pi, \Label{E:det-eta-sa2}
\end{align}
where $\Ph\big( \Det_{\tet}(D)\big)$ denotes the phase of the complex number $\Det_{\tet}(D)$.

If $D$ is not self-adjoint, \refe{det-eta-sa} is not true in general, because the numbers $\LD_{2\tet}(D^2)$, $\eta(D)$, and
$\zet_{2\tet}(0,D^2)$ need not be real. However, they are real and a version of \refe{det-eta-sa} and \refe{det-eta-sa2} holds for a class of
injective elliptic differential operators {\em whose spectrum is symmetric with respect to the real axis}. Though we will not use this result we
present it in the Appendix~\ref{S:det-eta-sa} for the sake of completeness.

%----------------------------------------------------------
\subsection{$\eta$-invariant and graded determinant}\Label{SS:eta-grdet}
Suppose now that $D= \bigoplus_{j=0}^dD_j$ as in \refe{D=bigoplus}. Choose  $\tet\in (-\pi/2,0)$ such that there are no eigenvalues of $D_j$ in the solid angles
$L_{(-\pi/2,\tet]}$ and $L_{(\pi/2,\tet+\pi]}$ for every $0\le j\le  d$. From the definition of the $\eta$-invariant it follows that
\[
    \eta\big(\, \pm D_j\,\big) \ = \ \pm\, \eta(D_j).
\]
Combining this equality with \refe{grldeterminantV} and \refe{det-eta} we obtain
\eq{grdet-eta}
    \LDgrtet(D) \ = \ \frac12\,\sum_{j=0}^d\, (-1)^j\,\LD_{2\tet}(D_j^2)  \ - \
   i\pi\,\Big(\,\eta(D)- \frac12\sum_{j=0}^d\,(-1)^j\zet_{2\tet}(0,D^2_j)\,\Big),
\end{equation}
where
\[
    \eta(D) \ = \ \sum_{j=0}^d\, \eta(D_j)
\]
is the $\eta$-invariant of the operator $D= \bigoplus_{j=0}^dD_j$.

Finally, note that, as in \refr{det-eta}.c, if the dimension of $M$ is odd, then $\zet_{2\tet}(0,D^2_j)= 0$, and \refe{grdet-eta} takes the form
\eq{grdet-eta-odd}
    \LDgrtet(D) \ = \ \frac12\,\sum_{j=0}^d\, (-1)^j\,\LD_{2\tet}(D^2)  \ - \
   i\pi\,\eta(D).
\end{equation}

%---------------------------------------
\subsection{Generalization}\Label{SS:general}
All the constructions and theorems of this section easily generalize to operators acting on a subspace of the space $C^\infty(M,E)$ of sections of $E$.

Let $D:C^\infty(M,E)\to C^\infty(M,E)$ be an injective elliptic differential operator with a self-adjoint leading symbol. Let $Q:C^\infty(M,E)\to
C^\infty(M,E)$ be a 0-th order pseudo-differential projection commuting with $D$. Then $V:= \IM{Q}\subset C^\infty(M,E)$ is a $D$-invariant subspace.
Hence, the decomposition \refe{L=sumLam} implies that
\[
    V \ = \ \overline{\bigoplus_{k\ge 1}\, \big(\, \Lam_k\cap V\,\big)}.
\]
and the restriction $D{|_V}$ of $D$ to $V$ has the same eigenvalues $\lam_1,\lam_2,\ldots$ as $D$ but with new multiplicities $m^V_1, m^V_2,\ldots$.
Note, that now $m_i^V\ge 0$ might vanish for certain $i$'s. Let $m^V_+$ (respectively, $m^V_-$) denote the number of eigenvalues (counted with their
algebraic multiplicities) of $D{|_V}$ on the positive (respectively, negative) part of the imaginary axis. Set
\eq{etaV}
 \begin{aligned}
    \eta_\tet(s,D|_V) \ &= \ \sum_{\RE\lam_k>0}\, m^V_k\ (\lam_k)_\tet^{-s}  \ - \
    \sum_{\RE\lam_k<0}\, m^V_k\ (-\lam_k)_\tet^{-s}, \\
    \eta(D|_V) \ &= \ \frac{\eta_\tet(0,D|_V)+m^V_+-m^V_-}2.
 \end{aligned}
\end{equation}
A verbatim repetition of the proof of \reft{det-eta} implies
\eq{det-etaV}
   \LD_{\tet}(D|_V) \ = \ \frac12\,\LD_{2\tet}(D^2|_V)  \ - \
   i\pi\,\Big(\,\eta(D|_V)- \frac12\zet_{2\tet}(0,D^2|_V)\,\Big),
\end{equation}
where we used the notation
\eq{zetDV}
    \zet_{2\tet}(s,D^2|_V) \ = \ \zet_{2\tet}(s,Q,D^2),
\end{equation}
cf. \refe{zetaQ}.

Finally, suppose that $V=\bigoplus_{j=0}^dV_j$ is given as in \refd{grdeterminantV}. Then
\eq{grdet-etaV}
    \LDgrtet(D) \ = \ \frac12\,\sum_{j=0}^d\, (-1)^j\,\LD_{2\tet}(D^2|_{V_j})  \ - \
   i\pi\,\Big(\,\eta(D)- \frac12\sum_{j=0}^d\,(-1)^j\zet_{2\tet}(0,D^2|_{V_j})\,\Big).
\end{equation}
Note, however, that an analogue of \refe{grdet-eta-odd} does not necessarily hold in this case even if $\dim{M}$ is odd, because
$\zet_{2\tet}(s,D^2|_{V_j})$ defined by \refe{zetDV}, is not a $\zet$-function of a {\em differential} operator and does not necessarily vanish
at $0$.

%----------------------------------------------------------------------------------------------
%----------------------------------------------------------------------------------------------
\section{Determinant as a Holomorphic Function}\Label{S:detanalytic}

In this section we explain that the determinant can be viewed as a holomorphic function on the space of elliptic differential operators. We also
discuss some applications of this result, which will be used in \refs{analytic} to show that the refined analytic torsion is a holomorphic
function and in \refs{deponmet} for studying the dependence of the graded determinant on the Riemannian metric.

%-----------------------------------
\subsection{Holomorphic curves in a Fr\'echet space}\Label{SS:hol}
Let $\calE$ be a complex Fr\'echet space and let $\calO\subset \CC$ be an open set. Recall (cf., e.g., \cite[Def.~3.30]{Rudin_FA}) that a map
$\gam:\calO\to \calE$ is called {\em holomorphic} if for every $\lam\in \calO$ the following limit exists,
\eq{holcurve}\notag
    \lim_{\mu\to \lam}\, \frac{\gam(\mu)-\gam(\lam)}{\mu-\lam}.
\end{equation}
We will refer to a holomorphic map $\gam:\calO\to \calE$ as a {\em holomorphic curve} in $\calE$.

Let $\calZ\subset \calE$ be a subset of a complex Fr\'echet space. By a holomorphic curve in $\calZ$ we understand a holomorphic map
$\gam:\calO\to \calE$ such that $\gam(\lam)\in \calZ$ for all $\lam\in \calO$.

Suppose now that $V\subset \CC^n$ is an open set. We call a map $f:V\to \calZ$ holomorphic if for each holomorphic curve $\gam:\calO\to V$ the
composition $f\circ\gam:\calO\to \calZ$ is a holomorphic curve in $\calZ$. Note that if $\calZ= \CC$ then, by Hartogs' theorem (cf., e.g.,
\cite[Th.~2.2.8]{HormanderSCV}), the above definition is equivalent to the standard definition of a holomorphic map.

%-------------------------------------
\subsection{The space of smooth functions as a Fr\'echet space}\Label{SS:smoothf}
The space $C^\infty_b(\RR^d)$ of bounded smooth complex-valued functions on $\RR^d$ with bounded derivatives has a natural structure of a Fr\'echet space (cf., e.g.,
\cite[Ch.~I]{ZimmerFunctAn}) with topology defined by the semi-norms
\eq{CinftyFr}
   \|f\|_{\alp} \ := \
     \sup_{x\in \RR^d}\, \left|\,\pa_x^\alp\,f(x)\,\right|,
\end{equation}
where $\alp=(\alp_1\nek \alp_d)\in (\ZZ_{\ge0})^d$ and $\pa_x^\alp:= \frac{\pa^{\alp_1}}{\pa x_1^{\alp_1}}\dots\frac{\pa^{\alp_d}}{\pa x_d^{\alp_d}}$.

%-------------------------------------
\subsection{A Fr\'echet space structure on the space of differential operators}\Label{SS:difopFr}
Let $M$ be a closed $d$-dimensional manifold and let $E$ be a complex vector bundle over $M$. Denote by $\Diff_m(M,E)$ the set of differential operators
$D:C^\infty(M,E)\to C^\infty(M,E)$ of order $\le{m}$ with smooth coefficients. It has a natural structure of a Fr\'echet space defined as follows. Consider a pair
$(\phi,\Phi)$ where $\phi:U\to \RR^d$ is a diffeomorphism (with $U\subset M$ an open set), and $\Phi:E|_U\to \CC^l\times{U}$ is a bundle map which identifies the
restriction $E|_U$ of $E$ to $U$ with the trivial bundle $\CC^l\times{U}\to U$. We refer to $(\phi,\Phi)$ as a {\em coordinate pair}.

Using the maps $\phi$ and $\Phi$ we can identify the restriction of an operator $D\in \Diff_m(M,E)$ to $U$ with the operator
\eq{DUphipsi}
  D_{(\phi,\Phi)} \ := \
   \sum_{|\bet|\le m}\,
      a_{(\phi,\Phi)}^{\bet}(x)\,\,\pa_x^\bet \ \in \ \Diff_m(\RR^d,\CC^l\times \RR^d),
\end{equation}
where $|\bet|= \sum_{j=1}^d\bet_j$ and
 \(
        a_{(\phi,\Phi)}^{\bet}(x)  =  \big\{\, a_{(\phi,\Phi);i,j}^{\bet}(x)\,\big\}_{i,j=1}^l
 \)
are smooth bounded matrix-valued functions on $\RR^d$, called the {\em coefficients} of $D$ with respect to the coordinate pair $(\phi,\Phi)$.

We now define a structure of a Fr\'echet space on $\Diff_m(M,E)$ using the semi-norms
\eq{DiffFr}
  \|D\|_{(\phi,\Phi);\bet;i,j}^{\alp} \ := \
    \big\|\,a_{(\phi,\Phi);i,j}^{\bet}\,\big\|_{\alp},
\end{equation}
where $(\phi,\Phi)$ runs over all coordinate pairs, $\alp,\bet\in (\ZZ_{\ge0})^d$ with $|\bet|\le{}m$,\ $1\le i,j\le l$, and the norm on the right hand side of
\refe{DiffFr} is defined by \refe{CinftyFr}.

%-------------------------------------
\subsection{Holomorphic curves in the space of differential operators}\Label{SS:holcurveDiff}
Suppose $\calO\subset \CC$ is an open set and consider a map $\gam:\calO\to \Diff_m(M,E)$. For every coordinate pair $(\phi,\Phi)$ we denote by
$a_{(\phi,\Phi)}^{\bet}(x;\lam)$ the coefficients of $\gam(\lam)$ with respect to the coordinate pair $(\phi,\Phi)$.

Clearly, $\gam$ is a holomorphic curve in $\Diff_m(M,E)$ with respect to the Fr\'echet space structure introduced in \refss{difopFr} if and only if for every
coordinate pair $(\phi,\Phi)$, every $\bet\in (\ZZ_{\ge0})^d$ with $|\bet|\le m$,\ and every $1\le i,j\le l$, the map $\lam\mapsto a_{(\phi,\Phi);i,j}^{\bet}(x;\lam)$
is a holomorphic curve in $C^\infty_b(\RR^d)$.

The following lemma follows immediately from the definitions.
\lem{comphol}
Let $\calO\subset \CC$ be an open set and, for $i=1,2$, let $\gam_i:\calO\to \Diff_{m_i}(M,E)$ be a holomorphic curve. Then $\lam\mapsto
\gam(\lam):= \gam_1(\lam)\circ\gam_2(\lam)$ is a holomorphic curve in $\Diff_{m_1+m_2}(M,E)$. Here $\gam_1(\lam)\circ\gam_2(\lam)$ is the
composition of the differential operators $\gam_1(\lam)$ and $\gam_2(\lam)$.
\elem

%------------------------------
\subsection{Determinant of a holomorphic curve of operators}\Label{SS:detoffamily}
Let $\Ell_{m}(M,E)\subset \Diff_m(M,E)$ denote the open set of elliptic differential operators of order $m$ and let  $\Ell_{m,\tet}(M,E)\subset
\Ell_m(M,E)$ be the open subset of operators which have $\tet$ as an Agmon angle. We denote by $\Ell'_{m,\tet}(M,E)$ the open subset of
invertible operators in $\Ell_{m,\tet}(M,E)$.  According to \refss{zet-det}, the function
\[
    \LD_\tet: \Ell'_{m,\tet}(M,E) \ \longrightarrow \ \CC
\]
is well defined. For $D\in \Ell_{m,\tet}(M,E)$ we set
\eq{detnoninv}
    \Det_\tet(D) \ = \ \begin{cases}
                          \exp\,\left(\,\LD_\tet(D)\,\right)\ \in \ \CC\backslash\{0\}, \qquad&\text{if $D$ is invertible};\\
                          0, \qquad&\text{otherwise}.
                       \end{cases}
\end{equation}

Further, we denote by $\Ell_{m,(\tet_1,\tet_2)}(M,E)\subset \Ell_m(M,E)$ the open subset of operators for which all the angles $\tet\in
(\tet_1,\tet_2)$ are principal, cf. \refss{agmonangle}. Any operator $D\in \Ell_{m,(\tet_1,\tet_2)}(M,E)$ has an Agmon angle $\tet\in
(\tet_1,\tet_2)$ and, by \refe{det-det},  the determinant $\Det_\tet(D)$ is independent of the choice of $\tet$ in the interval
$(\tet_1,\tet_2)$. The following theorem is well known, cf., for example, \cite[Corollary~4.2]{KontsevichVishik_short},
\th{dethol}
Let $E$ be a complex vector bundle over a closed manifold $M$ and let $\calO\subset  \CC$ be an open set.

a. \ Suppose $\gam:\calO\to \Ell'_{m,\tet}(M,E)$ is a holomorphic curve in $\Ell'_{m,\tet}(M,E)\subset \Diff_m(M,E)$. Then the function
$\calO\to \CC$, \ $\lam\mapsto \LD_\tet\big(\gam(\lam)\big)\in \CC$ is holomorphic.

b. \ Given angles $\tet_1<\tet_2$ and an operator $D\in \Ell_{m,(\tet_1,\tet_2)}(M,E)$, denote by $\Det_{(\tet_1,\tet_2)}(D)$ the determinant
$\Det_\tet(D)$ defined using any Agmon angle $\tet\in (\tet_1,\tet_2)$. Let $\gam:\calO\to \Ell_{m,(\tet_1,\tet_2)}(M,E)$ be a holomorphic curve
in $\Ell_{m,(\tet_1,\tet_2)}\subset \Diff_m(M,E)$. Then
\eq{dethol}
  \calO\ \longrightarrow \ \CC, \qquad \lam \ \mapsto \ \Det_{(\tet_1,\tet_2)}\big(\gam(\lam)\big)
\end{equation}
is a holomorphic function.
\eth
%The above theorem seems to be well known to the specialists. However, since we were unable to find it in the literature, we present its proof in
%Appendix~\ref{S:prlogdethol}.

\rem{holomorphicFr}
The theorem implies that the function $\Det_{(\tet_1,\tet_2)}(D)$\/ is {\em G\^ateaux holomorphic}\, on\linebreak
$\Ell_{m,(\tet_1,\tet_2)}(M,E)$, cf. \cite[Def.~3.1]{Dineen99}. Moreover, since $\Det_{(\tet_1,\tet_2)}$ is continuous on
$\Ell_{m,(\tet_1,\tet_2)}(M,E)$ it follows that this function is holomorphic in the sense of Definition~3.6 of \cite{Dineen99}. However, since
there seems to be no standard notion of a holomorphic function on a Fr\'echet space, we prefer to avoid this terminology.
\erem
\cor{etahol}
Suppose $E\to M$ is a complex Hermitian vector bundle over a closed manifold $M$. Let $\Ell'_{m,\operatorname{sa}}(M,E)$ denote the set of
invertible elliptic operators of order $m$ with self-adjoint leading symbol and let $\gam:\calO\to \Ell'_{m,\operatorname{sa}}(M,E)$ be a
holomorphic curve in $\Ell'_{m,\operatorname{sa}}(M,E)$. Then the function
\eq{etahol}
  \calO \ \longrightarrow \ \CC, \qquad\lam \ \mapsto \ e^{2\pi i\eta(\gam(\lam))}
\end{equation}
is holomorphic.
\ecor
\prf
By formula \refe{det-eta2} of \reft{det-eta}
\eq{det-etagam}
    e^{2\pi i \eta(\gam(\lam))} \ = \
    \frac{\Det_{(-\pi,0)}\big(\gam(\lam)^2\big)}
    {\Big(\,\Det_{(-\pi/2,0)}\big(\gam(\lam)\big)\,\Big)^2}
    \cdot e^{i\pi\zet_{2\tet}\big(0,\gam(\lam)^2\big)}.
\end{equation}
By \refl{comphol}, \, $\lam\mapsto \gam(\lam)^2$ is a holomorphic curve in $\Ell'_{2m,\operatorname{sa}}(M,E)$. Hence, by \reft{dethol}.b the
quotient on the right hand side of \refe{det-etagam} is a holomorphic function in $\lam$. It remains to show that
$\zet_{2\tet}\big(0,\gam(\lam)^2\big)$ depends holomorphically on $\lam$.

First, note that by \refe{tet-tet}, $\zet_{2\tet}\big(0,\gam(\lam)^2\big)$ is independent of $\tet$. By a result of Seeley \cite{Seeley67} (see
also \cite{ShubinPDObook}), the value $\zet_{2\tet}\big(0,\gam(\lam)^2\big)$ of the zeta-function of $\gam(\lam)^2$ is given by a local formula,
i.e., by an integral over $M$ of a $\CC$-valued differential form $\phi$ whose value at a point $x\in M$ is a rational function  of the symbol
of $\gam(\lam)$ and a finite number of its derivatives. It follows that the function $\calO\to \CC$, \, $\lam\mapsto
\zet_{2\tet}(0,\gam(\lam)^2)$ is holomorphic.
\eprf

Another important consequence of \reft{dethol} is the following
\cor{zeroanal}
Let $V\subset \CC^n$ be an open set and let $f:V\to  \Ell_{m,(\tet_1,\tet_2)}(M,E)$ be a holomorphic map in the sense of \refss{hol}. Then the
set
\[
   \Sig \ := \ \big\{\, \lam \in V:\, f(\lam) \ \text{is not invertible}\,\big\}
\]
is a complex analytic subset of $V$. In particular, if\/ $V$ is connected then so is $V\backslash\Sig$.
\ecor
\prf
In view of Hartogs' theorem (\cite[Th.~2.2.8]{HormanderSCV}), \reft{dethol}.b implies that the function $V\to \CC, \ \lam\mapsto
\Det_{(\tet_1,\tet_2)}(f(\lam))$ is holomorphic on $V$. By \refe{detnoninv},\/ $\Sig= \big\{\lam\in V:\, \Det_{(\tet_1,\tet_2)}(f(\lam))=
0\big\}$.
\eprf

%----------------------------------------------------------------------------------------------------
%-----------------------------------------------------------------------------------------------------
\section{Graded Determinant of the Odd Signature Operator}\Label{S:grdet}

In this section we define the graded determinant of the Atiyah-Patodi-Singer odd signature operator, \cite{APS2,Gilkey84}, of a flat vector
bundle $E$ over a closed Riemannian manifold $M$. In \refs{grdet-TRS} we show that, if $E$ admits an invariant Hermitian metric, then the
absolute value of this determinant is equal to the Ray-Singer analytic torsion \cite{RaSi1}. There is a similar, though slightly more
complicated, relationship between the graded determinant and the Ray-Singer torsion in the general case, cf. \reft{DetB-TRS}.  Thus, the graded
determinant of the odd signature operator can be viewed as a {\em refinement of the Ray-Singer torsion}.

%------------------------------
\subsection{Setting}\Label{SS:setgrdet}
Let $M$  be a smooth closed oriented manifold of odd dimension $d=2r-1$ and let $E\to M$ be a complex vector bundle over $M$ endowed with a flat connection $\n$. We
denote by $\n$ also the induced differential
\[
    \n:\, \Ome^\b(M,E) \ \longrightarrow \Ome^{\b+1}(M,E),
\]
where $\Ome^k(M,E)$ denotes the space of smooth differential forms of $M$ with values in $E$ of degree $k$.

%--------------------------------------------------------------
\subsection{Odd signature operator}\Label{SS:oddsign}
Fix a Riemannian metric $g^M$ on $M$ and let $*:\Ome^\b(M,E)\to \Ome^{d-\b}(M,E)$ denote the Hodge $*$-operator. Define the {\em chirality operator}
$\Gam:\Ome^\b(M,E)\to \Ome^\b(M,E)$ by the formula
\eq{Gam}
    \Gam\, \ome \ := \ i^r\,(-1)^{\frac{k(k+1)}2}\,*\,\ome, \qquad \ome\in \Ome^k(M,E),
\end{equation}
with $r$ given as above by \/ $r=\frac{d+1}2$. This operator is equal to the operator defined in \S3.2 of \cite{BeGeVe} as follows from applying Proposition~3.58 of
\cite{BeGeVe} in the case $\dim{M}$ is odd. Note that $\Gam^2=1$ and that $\Gam$ is self-adjoint with respect to the scalar product on $\Ome^\b(M,E)$ induced by the
Riemannian metric $g^M$ and by an arbitrary Hermitian metric on $E$.

%--------
\defe{oddsign}
The {\em odd signature operator} is the operator $B=B(\n,g^M):\Ome^\b(M,E)\to \Ome^\b(M,E)$ is defined by
\begin{equation} \Label{E:oddsign}
    B \ = \ \Gam\,\n\ + \ \n\,\Gam.
\end{equation}
We denote by $B_k$ the restriction of $B$ to the space $\Ome^{k}(M,E)$.
\edefe
Explicitly, for  $\ome\in \Ome^{k}(M,E)$ one has
\begin{equation} \Label{E:oddsignnoGam}
    B_k\,\ome \ := \ i^r(-1)^{\frac{k(k+1)}2+1}\,\big(\,(-1)^k*\n-\n*\,\big)\,\ome
    \ \in \ \Ome^{d-k-1}(M,E)\oplus  \Ome^{d-k+1}(M,E).
\end{equation}

The odd signature operator was introduced by Atiyah, Patodi, and Singer, \cite[p.~44]{APS1}, \cite[p.~405]{APS2}, in the case when $E$ is
endowed with a Hermitian metric {\em invariant with respect to $\n$} (i.e. invariant under parallel transport by $\n$). The general case was
studied by Gilkey, \cite[p.~64--65]{Gilkey84}.

\lem{symbolofB}
Suppose that $E$ is endowed with a Hermitian metric $h^{E}$. Denote by $\<\cdot,\cdot\>$ the scalar product on $\Ome^\b(M,E)$ induced by $h^{E}$ and
the Riemannian metric $g^M$ on $M$.

1. \  $B$ is elliptic and its leading symbol is symmetric with respect to the Hermitian metric $h^{E}$.

2. \ If, in addition, the metric $h^{E}$ is invariant with respect to the connection $\n$, then $B$ is symmetric with respect to the scalar
product $\<\cdot,\cdot\>$, $$B^*= B.$$ If\, the metric $h^{E}$ is not invariant, then, in general, $B$ is not symmetric.
\elem
The proof of the lemma is a simple calculation. The first part is already stated in \cite[p.~405]{APS2}. The second part is proven in the Remark on
page~65 of \cite{Gilkey84}.

%------------------------------------------------------------------
\subsection{Assumptions}\Label{SS:assumptions}
In this paper we study the odd signature operator $B$ and the analytic torsion under the following simplifying assumptions. The general case is addressed in
\cite{BrKappelerRATdetline}.

\subsection*{Assumption~I}
The connection $\n$ is {\em acyclic}, i.e., the {\em twisted deRham complex}
\eq{deRham}
    \begin{CD}
        0 \ \to \Ome^0(M,E) @>\n>> \Ome^1(M,E) @>\n>>\cdots @>\n>> \Ome^d(M,E) \ \to \ 0
    \end{CD}
\end{equation}
is acyclic,
\[
   \IM\big(\n|_{\Ome^{k-1}(M,E)}\big)\ =\ \Ker\big(\n|_{\Ome^k(M,E)}\big) \qquad\text{for every}\quad k=1\nek d.
\]

\subsection*{Assumption~II}
The odd signature operator $B:\Ome^\b(M,E)\to \Ome^\b(M,E)$ is bijective.

%-----------------------------------------------------------------------
\subsection{Hermitian connection}\Label{SS:openset}
Suppose that there exists a Hermitian metric $h^E$ on $E$ invariant with respect to $\n$ (in this case we say that the connection $\n$ is {\em
Hermitian}). Then Assumption~II follows from Assumption~I. Indeed, in this case the operator $B$ is symmetric with respect to the scalar product
$\<\cdot,\cdot\>$, defined by the metrics $g^M$ and $h^{E}$, cf. \refl{symbolofB}. Hence, we only need to show that $\Ker{}B= \{0\}$. Let $\n^*$
denote the formal adjoint of the operator $\n:\Ome^\b(M,E)\to \Ome^{\b+1}(M,E)$ with respect to the scalar product $\<\cdot,\cdot\>$. Since the
metric $h^{E}$ is flat, we obtain, cf. \cite[\S6.1]{Warner},
\eq{n*=}
    \n^* \ = \ \Gam\,\n \, \Gam.
\end{equation}
Using this identity and  the definition \refe{oddsign} of $B$ we see that
\eq{B2}
    B^2 \ = \ \n^*\,\n \ + \ \n\,\n^*
\end{equation}
is the Laplacian. Thus $\Ker B= \Ker B^2$ is isomorphic to the cohomology of the complex \refe{deRham}, and, hence, is trivial by Assumption~I. Conversely, these
arguments show at the same time that, in the case considered, Assumption~II implies Assumption~I.

%-----------------------------------------------------------------------
\subsection{Connections which are close to a Hermitian connection}\Label{SS:openset1}
In this paper we are interested in the study of connections which are {\em close to a Hermitian connection} in the following sense:

Let $\Ome^1(M,\End(E))$ denote the space of differential 1-forms on $M$ with values in the bundle $\End(E)$ of endomorphisms of $E$. A Hermitian
metric on $E$ and a Riemannian metric on $M$ define a natural norm $|\cdot|$ on the bundle $\Lam^1(T^*M)\otimes\End(E)\to M$. Using this norm we
define the sup-norm
\[
  \|\ome\|_{\sup} \ := \ \max_{x\in M}\, |\ome(x)|, \qquad \ome\in \Ome^1(M,\End E)
\]
on $\Ome^1(M,\End E)$. The topology defined by this norm is independent of the metrics and is called the $C^0$-topology on $\Ome^1(M,\End E)$.

Let $\calC(E)$ denote the space of connections on the bundle $E$. By choosing a connection $\n_0$ we can identify this space with $\Ome^1(M,\End
E)$ associating to a connection $\n\in \calC(E)$ the 1-form $\n-\n_0\in \Ome^1(M,\End E)$. By this identification the $C^0$-topology on
$\Ome^1(M,\End E)$ provides a topology on $\calC(E)$ which is independent of the choice of $\n_0$ and is called the {\em $C^0$-topology} on the
space of connections.

Finally, we denote by $\Flat(E)\subset \calC(E)$ the set of flat connections on $E$ and by $\Flat'(E,g^M)\subset \Flat(E)$ the set of flat
connections satisfying Assumption~I and II of \refss{assumptions}. The topology induced on these sets by the $C^0$-topology on $\calC(E)$ is
also called the $C^0$-topology. The discussion of the previous subsection implies that {\em $\Flat'(E,g^M)$ contains  all the acyclic Hermitian
connections}.

\prop{neigofh}
$\Flat'(E,g^M)$ is a $C^0$-open subset of\, $\Flat(E)$, which contains all acyclic Hermitian connections on $E$.
\eprop
\prf
We already know that $\Flat'(E,g^M)$ contains all acyclic Hermitian connections on $E$. Hence it is enough to show that $\Flat'(E,g^M)$ is open
in $C^0$-topology.

Let $\n\in \Flat'(E,g^M)$ and suppose that $\n'\in \Flat(E)$ is sufficiently close to $\n$ in $C^0$-topology. Let $B= B(\n,g^M),\ B'=B(\n',g^M)$
denote the odd signature operators associated to the connections $\n$ and $\n'$, respectively. Then $B-B'$ is a 0'th order differential operator
on $\Ome^\b(M,E)$ and, hence, is bounded. Moreover, if $\n$ is close to $\n'$ in the $C^0$-topology, then $B'-B:\Ome^\b(M,E)\to \Ome^\b(M,E)$ is
small in the operator norm, when $\Ome^\b(M,E)$ is endowed with the $L^2$-norm induced by the Riemannian metric on $M$ and the Hermitian metric
on $E$. We refer to this operator norm as the {\em standard} operator norm and denote it by $\|\cdot\|$.

Since the operator $B$ satisfies Assumption~II, its inverse $B^{-1}$ can be viewed as a bounded operator on the $L^2$-completion
$L^2\Ome^\b(M,E)$ of $\Ome^\b(M,E)$. If $B-B'$ is sufficiently small so that $\|(B'-B)B^{-1}\|<1$, then $B'$, viewed as an unbounded operator on
$L^2\Ome^\b(M,E)$, has a bounded inverse given by the formula
\[
   (B')^{-1} \ = \ B^{-1}\, \big(\,\Id \ + \ (B'-B)B^{-1}\,\big)^{-1}.
\]
By elliptic theory, $(B')^{-1}$ maps the space of smooth forms $\Ome^\b(M,E)$ to itself. Hence, $B'$ satisfies Assumption~II.
\eprf

%----------------------------------------
\subsection{Decomposition of the odd signature operator}\Label{SS:decompos}
Set
\begin{equation}\Label{E:decompos}\notag
 \begin{aligned}
    {}&\Ome^{\even}(M,E) \ := \ \bigoplus_{p=0}^{r-1}\, \Ome^{2p}(M,E), \quad
    \Ome^{\odd}(M,E) \ := \ \bigoplus_{p=1}^r\, \Ome^{2p-1}(M,E),\\
    {}&B_{\even} \ := \ \bigoplus_{p=0}^{r-1}\, B_{2p}:\, \Ome^{\even}(M,E)\ \longrightarrow \ \Ome^{\even}(M,E), \\
    {}&B_{\odd} \ := \ \bigoplus_{p=1}^{r}\, B_{2p-1}:\, \Ome^{\odd}(M,E)\ \longrightarrow \ \Ome^{\odd}(M,E),
 \end{aligned}
\end{equation}

Using that $\Gam^2=1$ we obtain
\eq{B=*B*}
    B_{\odd}\ =\ \Gam\circ B_{\even}\circ \Gam\,\big|_{\Ome^\odd(M,E)}.
\end{equation}
Hence, the whole information about the odd signature operator is encoded in its {\em even part} $B_{\even}$. The operator $B_{\even}$ can be expressed by the
following formula, which is slightly simpler than \refe{oddsignnoGam}:
\eq{oddsigneven}
    B_{\even }\,\ome \ := \ i^r(-1)^{p+1}\,\big(\,*\n-\n*\,\big)\,\ome,\qquad\text{for}\quad
                 \ome\in \Ome^{2p}(M,E).
\end{equation}

Assume now that $\n\in \Flat'(E,g^M)$, i.e., that Assumption~I and II of \refss{assumptions} are satisfied. From Assumption~I we conclude that
the kernel and the image of the operator $\n:\Ome^\b(M,E)\to \Ome^\b(M,E)$ coincide. Hence,
\eq{ker-im}
 \begin{aligned}
    \Ker\,(\Gam\,\n) \ &= \ \Ker\,\n \ = \ \IM\,\n\ = \ \IM\,(\n\,\Gam)\\
    \Ker\,(\n\,\Gam)\ &= \ \Gam\,\big(\,\Ker\,\n\,\big)\ =\ \Gam\,\big(\,\IM\,\n\,\big) \ = \ \IM\,(\Gam\,\n).
 \end{aligned}
\end{equation}
We set
\eq{ome+-}
 \begin{aligned}
  \Ome^k_+(M,E) \ &:= \ \Ker\,(\n\,\Gam)\,\cap\,\Ome^k(M,E) \ = \ \big(\Gam\,\Ker\,\n\big)\,\cap\,\Ome^k(M,E);\\
  \Ome^k_-(M,E) \ &:= \ \Ker\,(\Gam\,\n)\,\cap\,\Ome^k(M,E)
  \ = \ \Ker\,\n\,\cap\,\Ome^k(M,E),
 \end{aligned}
\end{equation}
and refer to $\Ome^k_+(M,E)$ and $\Ome^k_-(M,E)$ as the positive and negative subspaces of $\Ome^k(M,E)$.
\footnote{Note, that our grading is opposite to the one considered in \cite[\S2]{BFK3}.}

Assumption~II of \refss{assumptions} then implies that
\eq{OmecapOme}
    \Ker\Big(\, \n\,\Gam{\big|_{\Ome^k(M,E)}}\,\Big) \,\cap\,
    \Ker\Big(\, \Gam\,\n{\big|_{\Ome^k(M,E)}}\,\Big)  \ = \ \{0\}, \qquad k=0\nek d
\end{equation}
(as $B$ is one-to-one) and
\eq{Ome+Ome}
    \IM\Big(\, \n\,\Gam{\big|_{\Ome^{d-k+1}(M,E)}}\,\Big) \ + \
    \IM\Big(\, \Gam\,\n{\big|_{\Ome^{d-k-1}(M,E)}}\,\Big)  \ = \ \Ome^k(M,E), \qquad k=0\nek d
\end{equation}
(as $B$ is onto). Combining \refe{OmecapOme} and \refe{Ome+Ome} with \refe{ker-im} and \refe{ome+-} we conclude that
\eq{Ome=2}
    \Ome^k(M,E)\ =\ \Ome^k_+(M,E)\,\oplus\, \Ome^k_-(M,E).
\end{equation}
Clearly,
\eq{*Ome}
    \Gam:\,\Ome^k_\pm(M,E) \ \longrightarrow \ \Ome^{d-k}_\mp(M,E), \qquad\qquad k=0\nek d.
\end{equation}
%%%%%%
\rem{Ome=}
Suppose that $h^E$ is a flat Hermitian metric on $E$. Then, using \refe{n*=} and \refe{ker-im}, we obtain
\eq{Ome-=}
  \Ome^k_+(M,E) \ = \ \Ker\,\n^*\,\cap\,\Ome^k(M,E),
\end{equation}
where $\n^*$ is the adjoint of $\n$ with respect to the scalar product induced by the metrics $g^M$ and $h^{E}$.
\erem
Let $B_{k}^\pm$ denote the restriction of $B$ to $\Ome^k_\pm(M,E)$. By \refe{oddsign} and \refe{ker-im},
\eq{Bpm}
 \begin{aligned}
    B_{k}^+ \ = \ \Gam\,\n &:\, \Ome^k_+(M,E) \ \longrightarrow \ \Ome^{d-k-1}_+(M,E),\\
    B_{k}^- \ = \ \n\,\Gam &:\, \Ome^k_-(M,E) \ \longrightarrow \ \Ome^{d-k+1}_-(M,E).
 \end{aligned}
\end{equation}
It follows from Assumption~II of \refss{assumptions} that both, $B_{k}^+$ and $B_{k}^-$, are invertible.

%--------------------
\subsection{Graded determinant of the odd signature operator}\Label{SS:drdetoddsign}
Set
\[
    \Ome^{\even}_\pm(M,E)\ =\ \bigoplus_{p=0}^{r-1}\, \Ome^{2p}_\pm(M,E)
\]
and let $B_{{\even}}^\pm$ denote the restriction of $B_{{\even}}$ to the space $\Ome^{\even}_\pm(M,E)$. Then
\[
    B_{\even}^\pm:\, \Ome^{\even}_\pm(M,E)\ \longrightarrow \  \Ome^{\even}_\pm(M,E).
\]

We consider $\Ome^{\even}(M,E)$ as a graded vector space
\[
    \Ome^{\even}(M,E) \ = \ \Ome^{\even}_+(M,E)\,\oplus\, \Ome^{\even}_-(M,E).
\]
By \refd{grdeterminantV}, the graded determinant of the odd signature operator is
\eq{grdetB}
    \Detgrtet(B_{\even}) \ := \ e^{\LDgrtet(B_{\even})},
\end{equation}
where $\tet\in (-\pi,0)$ is an Agmon angle for the operator $B=B_{\even}\oplus B_{\odd}$, cf. \refd{Agmon}, and
\eq{grldetB}
    \LDgrtet(B_{\even}) \ := \ \LD_\tet\big(B_{\even}^+\big) \ - \ \LD_\tet\big(-B_{\even}^-\big)
    \ \in \ \CC.
\end{equation}
%As the leading symbol of $B$ is symmetric (cf. \refr{oddsign}.1), we can choose an Agmon angle in the interval $(-\pi/2,0)$.
%%

According to \refe{det-det}, $\Detgrtet(B_{\even})$ is independent of the choice of the Agmon angle $\tet\in (-\pi,0)$.

%------------------------------------------------------------------------
%------------------------------------------------------------------------
\section{Relationship with the $\eta$-invariant}\Label{S:grdet-eta}

In this section we use the notations of the previous section and assume that, for a given pair $(\n,g^M)$, Assumptions~I and II of \refss{assumptions} are satisfied.
In particular the operator $B= B(\n,g^M)$ is bijective. It follows that the operators $B_{\even}$ and $B_k^+:\Ome^k_+(M,E)\to \Ome^{d-k-1}_+(M,E)$ \ ($0\le k\le d-1$)
are also invertible.
%--------------------
\subsection{Graded determinant and $\eta$-invariant}\Label{SS:grdetB-eta}
To simplify the notation set
\eq{eta=eta}
    \eta \ = \ \eta(\n,g^M) \ := \  \eta(B_{\even}),
\end{equation}
and
\eq{calB}
\begin{aligned}
    \xi \ = \ \xi(\n,g^M,\tet) \ &:= \ \frac12\,\sum_{k=0}^{d-1}\,(-1)^k\LD_{2\tet}\big(\,B_{d-k-1}^+\circ B_{k}^+\,\big)
    \\ &= \ \frac12\,\sum_{k=0}^{d-1}\,(-1)^k\LD_{2\tet} \Big(\,
    {(\Gam\,\n)^2}{\big|_{\Ome^k_+(M,E)}}\,\Big).
\end{aligned}
\end{equation}

%-----------
\th{DetB-eta2}
Let\, $\n\in \Flat'(E,g^M)$ be a flat connection on a vector bundle $E$ over a closed oriented Riemannian manifold $(M,g^M)$ of odd dimension $d=2r-1$. Let $\tet\in
(-\pi/2,0)$ be an Agmon angle for $B$ such that there are no eigenvalues of the operator $B$ in the solid angles $L_{(-\pi/2,\tet]}$ and $L_{(\pi/2,\tet+\pi]}$. Then
\eq{DetB-eta2}
    \LDgrtet(B_{\even}) \ = \ \xi \ - \  i\pi\,\eta.
\end{equation}
\eth
The rest of this section is devoted to the proof of \reft{DetB-eta2}. By \refe{grdet-etaV}, it is enough to show that
\begin{gather}
    2\,\xi \ = \ \LD_{2\tet}\,(B_{\even}^+)^2 \ - \ \LD_{2\tet}\,(B_{\even}^-)^2; \Label{E:calBalp=}\\
    \zet_{2\tet}\big(0,(B_{\even}^+)^2\big) \ - \ \zet_{2\tet}\big(0,(B_{\even}^-)^2\big) \ = \ 0. \Label{E:zetB2=}
\end{gather}

%----------------------------------------
\subsection{Calculation of $\zet_{2\tet}\big(s,(B_{\even}^+)^2\big)-\zet_{2\tet}\big(s,(B_{\even}^-)^2\big)$}\Label{SS:zetBpm}
Set
\eq{A+}\notag
     A_k^+ \ := \ B_{k}^+\ + \ B_{d-k-1}^+:\, \Ome^k_+(M,E)\oplus \Ome^{d-k-1}_+(M,E)
     \ \longrightarrow \ \Ome^k_+(M,E)\oplus \Ome^{d-k-1}_+(M,E),
\end{equation}
for $k=0\nek r-2$, and
\eq{A+r-1}\notag
     A_{r-1}^+ \ := \ B_{r-1}^+:\,\Ome^{r-1}_+(M,E) \ \longrightarrow \ \Ome^{r-1}_+(M,E).
\end{equation}
Similarly, set
\eq{A-}\notag
     A_k^- \ := \ B_{k}^-\ + \ B_{d-k+1}^-:\, \Ome^k_-(M,E)\oplus \Ome^{d-k+1}_-(M,E)
     \ \longrightarrow \ \Ome^k_-(M,E)\oplus \Ome^{d-k+1}_-(M,E)
\end{equation}
for $k=1\nek r-1$, and
\eq{A-r}
     A_{r}^- \ := \ B_{r}^-:\,\Ome^{r}_-(M,E) \ \longrightarrow \ \Ome^{r}_-(M,E).
\end{equation}
Then
\eq{A2=}
% \begin{aligned}
    (A^+_k)^2 \ = \ {(\Gam\,\n)^2}{\big|_{\Ome^k_+(M,E)\oplus\Ome^{d-k-1}_+(M,E)}},\qquad
    (A^-_k)^2 \ = \ {(\n\,\Gam)^2}{\big|_{\Ome^k_-(M,E)\oplus\Ome^{d-k+1}_-(M,E)}}.
% \end{aligned}
\end{equation}
Hence,
\eq{A2=*n+n*}
    \zet_{2\tet}\big(s,(A^+_k)^2\big) \ = \ \zet_{2\tet}\big(s,{(\Gam\,\n)^2}{\big|_{\Ome^k_+(M,E)}}\big)
    \ + \
    \zet_{2\tet}\big(s,{(\Gam\,\n)^2}{\big|_{\Ome^{d-k-1}_+(M,E)}}\big).
\end{equation}

From \refe{B=*B*} and \refe{*Ome} we get
\eq{A+=*A-*}
    A_k^+ \ = \ \Gam\circ A_{k-1}^-\circ \Gam.
\end{equation}
Hence,
\eq{zetA+=zetA-}
    \zet_{2\tet}\big(s,(A_k^+)^2\big) \ = \ \zet_{2\tet}\big(s,(A_{k-1}^-)^2\,\big).
\end{equation}
Since $B_{\even}^\pm$ is a direct sum of the operators $A^\pm_{2p}$, we obtain from \refe{zetA+=zetA-} that
\meq{zetB2=2}\notag
    \zet_{2\tet}\big(s,(B_{\even}^+)^2\big) \ - \ \zet_{2\tet}\big(s,(B_{\even}^-)^2\big) \ = \
    \sum_{p=0}^{[(r-1)/2]}\,\zet_{2\tet}\big(s,(A_{2p}^+)^2\big) \ - \
    \sum_{p=1}^{[r/2]}\,\zet_{2\tet}\big(s,(A_{2p}^-)^2\big)
    \\ = \
    \sum_{p=0}^{[(r-1)/2]}\,\zet_{2\tet}\big(s,(A_{2p}^+)^2\big) \ - \
    \sum_{p=1}^{[r/2]}\,\zet_{2\tet}\big(s,(A_{2p-1}^+)^2\big)
    \ = \
    \sum_{k=0}^{r-1}\,(-1)^k\,\zet_{2\tet}\big(s,(A_{k}^+)^2\big).
\end{multline}
Combining this equality with \refe{A2=*n+n*} we get
\eq{zetB2=3}
    \zet_{2\tet}\big(s,(B_{\even}^+)^2\big) \ - \ \zet_{2\tet}\big(s,(B_{\even}^-)^2\big) \ = \
    \sum_{k=0}^{d-1}\,(-1)^k\,\zet_{2\tet}\big(s,{(\Gam\,\n)^2}{\big|_{\Ome^k_+(M,E)}}\big).
\end{equation}

%--------------------------
\lem{zetPQ=zetQP}
Let $F_1, F_2$ be vector bundles over $M$ and let $P:C^\infty(M,F_1)\to C^\infty(M,F_2)$ and $Q:C^\infty(M,F_2)\to C^\infty(M,F_2)$ be invertible elliptic
pseudo-differential operators such that $\phi$ is an Agmon angle for $PQ$ and $QP$. Then, every regular point $s\in \CC$ of the function $s\mapsto \zet_{\phi}(s,PQ)$
is also a regular point of $s\mapsto \zet_{\phi}(s,QP)$ and
\eq{zetPQ=zetQPs}
    \zet_{\phi}(s,PQ) \ = \ \zet_{\phi}(s,QP).
\end{equation}
In particular,
\eq{zetPQ=zetQP}
    \zet_{\phi}(0,PQ) \ = \ \zet_{\phi}(0,QP).
\end{equation}
\elem
\prf
For every elliptic operator $D$ with Agmon angle $\phi$
\[
    QD^{-s-1}_{\phi}Q^{-1} \ = \ \big(\,QDQ^{-1}\,\big)_{\phi}^{-s-1}.
\]
Hence,
\eq{PQs=QPs}
    Q(PQ)_{\phi}^{-s-1} \ = \ \big[\,Q(PQ)_{\phi}^{-s-1}Q^{-1}\,\big]\, Q \ = \ (QP)_{\phi}^{-s-1}Q.
\end{equation}

Recall that if\/ $T$ and $S$ are operators such that the composition $TS$ is of trace class, then $ST$ is also of trace class and $\Tr(TS)= \Tr(ST)$, cf.
\cite[Ch.~III, Th.~8.2]{GohbergKrein_book69}. Using this equality and \refe{PQs=QPs} we obtain
\eq{QP-PQ}
 \begin{aligned}
    \zet_{\phi}(s,PQ) \ &= \ \Tr\, (PQ)^{-s}_{\phi} \ = \  \Tr\, \big[\,(PQ)^{-s-1}_{\phi} PQ\,\big]
    \\ &= \ \Tr\, \big[\,Q(PQ)^{-s-1}_{\phi} P\,\big] \ = \ \Tr\, (QP)^{-s}_{\phi} \ = \ \zet_{\phi}(s,QP).
 \end{aligned}
\end{equation}

Since both the left and the right hand sides of \refe{QP-PQ} are analytic in $s$, this equality holds for all regular points of the function
$s\mapsto \zet_{\phi}(s,PQ)$.
\eprf

From \refe{zetPQ=zetQPs}, we conclude that for all regular points of the function $s\mapsto \zet_{2\tet}\big(s,{(\Gam\n)^2}{\big|_{\Ome^k_+(M,E)}}\big)$ the following
equality holds
\meq{zet*n=zetn*}
    \zet_{2\tet}\big(s,{(\Gam\,\n)^2}{\big|_{\Ome^k_+(M,E)}}\big) \ = \
    \zet_{2\tet}\big(s,(\Gam\,\n\,\Gam){\big|_{\Ome^{k+1}_-(M,E)}}{\n}{\big|_{\Ome^k_+(M,E)}}\big)
    \\ = \
    \zet_{2\tet}\big(s,{\n}{\big|_{\Ome^k_+(M,E)}}(\Gam\,\n\,\Gam){\big|_{\Ome^{k+1}_-(M,E)}}\big)
    \ = \
    \zet_{2\tet}\big(s,{(\n\,\Gam)^2}{\big|_{\Ome^{k+1}_-(M,E)}}\big).
\end{multline}

From \refe{A2=}, \refe{OmecapOme}, and \refe{Ome+Ome}, we get
\begin{multline}\notag
    \zet_{2\tet}\big(s,{B^2}{\big|_{\Ome^k(M,E)}}\big) \ = \
    \zet_{2\tet}\Big(s,\big(\,{(\Gam\,\n)^2}+{(\n\,\Gam)^2}\,\big){\big|_{\Ome^k(M,E)}}\Big)
    \\ = \
    \zet_{2\tet}\big(s,{(\Gam\,\n)^2}{\big|_{\Ome^k_+(M,E)}}\big) \ + \
    \zet_{2\tet}\big(s,{(\n\,\Gam)^2}{\big|_{\Ome^k_-(M,E))}}\big).
\end{multline}
Using this equality, \refe{zetB2=3}, and \refe{zet*n=zetn*}, we obtain
\eq{zetB2=5}
    \zet_{2\tet}\big(s,(B_{\even}^+)^2\big) \ - \ \zet_{2\tet}\big(s,(B_{\even}^-)^2\big) \\ = \
    \sum_{k=0}^d\, (-1)^{k+1}\,k\,\,
    \zet_{2\tet}\big(s,{B^2}{\big|_{\Ome^k(M,E)}}\big).
\end{equation}

%-------------------------------------------
\subsection{Proof of \reft{DetB-eta2}}\Label{SS:prDetB-eta2}
From \refe{calB} and \refe{zetB2=3} we conclude that
\eq{calBalp=2}
    2\,\xi \ = \ \frac{d}{ds}{\big|_{s=0}}\Big[\,
      \zet_{2\tet}\big(s,(B_{\even}^+)^2\big) \ - \ \zet_{2\tet}\big(s,(B_{\even}^-)^2\big) \,\Big].
\end{equation}
Hence \refe{calBalp=} is established.

Since the dimension of $M$ is odd, the $\zet$-function of every elliptic differential operator of even order vanishes at 0, cf. \cite{Seeley67}.
Hence, by \refl{symbolofB}.1, the equality \refe{zetB2=5} implies \refe{zetB2=}. \hfill$\square$

%------------------------------------------------------------------------
%------------------------------------------------------------------------
\section{Comparison with the Ray-Singer Torsion}\Label{S:grdet-TRS}

%-------------------
\subsection{Ray-Singer torsion}\Label{SS:RaySinger0}
Let $E\to M$ be a complex vector bundle over a closed oriented manifold $M$ of odd dimension $d=2r-1$ and let $\n$ be an acyclic flat connection on $E$. Fix a
Riemannian metric $g^M$ on $M$ and a Hermitian metric $h^{E}$ on $E$. Let $\n^*$ denote the adjoint of $\n$ with respect to the scalar product $\<\cdot,\cdot\>$ on
$\Ome^\b(M,E)$ defined by $h^{E}$ and the Riemannian metric $g^M$. If $\n$ is acyclic (i.e., Assumption~I of \refss{assumptions} is satisfied) the {\em Ray-Singer
torsion} $\TRS$ of $E$, \cite{RaSi1,BisZh92,BFK3}, is defined by
\eq{RaySingertor}
    \TRS \ = \TRS(\n) := \ \exp\,\Big(\,
     \frac12\,\sum_{k=0}^d\,(-1)^{k+1}\,k\, \LD_{-\pi}\big(\, (\n^*\n+\n\n^*){\big|_{\Ome^k(M,E)}}\,\big)\,\Big),
\end{equation}
where $\n^*$ denotes the adjoint of $\n$ with respect to the scalar product $\<\cdot,\cdot\>$ induced by $g^M$ and $h^{E}$. Note that, as $\n$ is assumed to be
acyclic, $(\n^*\n+\n\n^*)_{|_{\Ome^k(M,E)}}$ is a strictly positive operator and, therefore,  $\LD_{-\pi}\big(\, (\n^*\n+\n\n^*){\big|_{\Ome^k(M,E)}}\big)$ is well
defined.

The Ray-Singer torsion is a positive number, which, {\em  in the case considered}, is independent of the the Hermitian metric $h^{E}$ and the
Riemannian metric $g^M$, cf. \cite{RaSi1,BisZh92}. We denote by $\log\TRS$ the value at $\TRS$ of the principal branch of the logarithm.

The determinants in \refe{RaySingertor} are defined using the spectral cut $R_{-\pi}$ along the negative real axis. Since the spectrum of the
operator $\n^*\n+\n\n^*$ lies on the positive real axis, we can replace it with a spectral cut $R_\phi$ for any $\phi\not=0$ without changing
the formula. In particular, we can take the spectral cut along $R_{2\tet}$, where $\tet\in (-\pi/2,0)$ is an Agmon angle for the odd signature
operator $B$.

Using the decomposition \refe{Ome=2}, we have
\eq{RaySingertor2}
  \begin{aligned}
    \log\,\TRS \ &= \ \frac12\,\sum_{k=0}^d\,(-1)^k\, \LD_{2\tet}\big(\, {\n^*\n}{\big|_{\Ome^k_+(M,E)}}\,\big)
    \\ &= \
    \frac12\,\sum_{k=0}^d\,(-1)^{k+1}\, \LD_{2\tet}\big(\, {\n\n^*}{\big|_{\Ome^k_-(M,E)}}\,\big).
  \end{aligned}
\end{equation}
This formula is proven, for example, on page~340 of \cite{BFK3}.

Suppose now that the Hermitian metric $h^{E}$ on $E$ is invariant with respect to $\n$. From \refe{n*=}, we obtain
\[
    \n^*\,{\n}{\big|_{\Ome^k_+(M,E)}} \ = \ {(\Gam\,\n)^2}{\big|_{\Ome^k_+(M,E)}}.
\]
Hence, we can rewrite \refe{RaySingertor2} as
\eq{RaySingertor3}
    \log\,\TRS \ = \
    \frac12\,\sum_{k=0}^d\,(-1)^k\LD_{2\tet} \Big(\, {(\Gam\,\n)^2}{\big|_{\Ome^k_+(M,E)}}\,\Big)
    \ = \ \xi(\n,g^M,\tet),
\end{equation}
where the number $\xi=\xi(\n,g^M,\tet)$ is defined in \refe{calB}.

The next theorem generalizes this result to the case where $h^E$ is not necessarily invariant. Recall that the set $\Flat'(M,g^M)$ is defined in
\refss{openset1}.
%%%
\th{DetB-TRS}
Assume that $M$ is a closed oriented manifold and $g^M$ is a Riemannian metric on $M$. Then there exists an $C^0$-open (cf. \refss{openset1}) neighborhood $U\subset
\Flat(E)$ of the set of acyclic Hermitian connections on $E$, such that for every connection $\n\in U$ we have $\n\in \Flat'(M,g^M)$ and
\eq{DetB-TRS}
    \log\,\TRS(\n) \ = \
    \frac12\,\RE\,\sum_{k=0}^d\,(-1)^k\LD_{2\tet} \Big(\, {(\Gam\,\n)^2}{\big|_{\Ome^k_+(M,E)}}\,\Big)
    \  = \ \RE\,\xi(\n,g^M,\tet).
\end{equation}
Hence, in view of \refe{DetB-eta2}, for every $\n\in U$, we obtain
\eq{DetB-TRS2}
    \Big|\,\Detgrtet(B_{\even})\,\Big| \ = \ \TRS(\n)\cdot e^{\pi\IM\eta(\n,g^M)}.
\end{equation}
\eth
If $\n$ is an acyclic Hermitian connection, then the operator $B_{\even}$ is self-adjoint with respect to the inner product given by $g^M$ and the invariant Hermitian
metric $h^E$ on $E$. Thus, the $\eta$-invariant $\eta= \eta(\n,g^M)$ is real. Hence, \reft{DetB-TRS} implies the following
\cor{DetB-TRS}
If \/ $\n$ is an acyclic Hermitian connection then
\eq{DetB-TRSherm}
    \big|\,  \Detgrtet(B_{\even}) \,\big| \ = \ \TRS(\n).
\end{equation}
\ecor
In particular, \refc{DetB-TRS} implies that, under the given assumptions, $\Detgrtet(B_{\even})$ contains all the information about the
Ray-Singer torsion, and, hence, ``refines" it by having a phase.

The rest of this section is occupied with the proof of \reft{DetB-TRS}.

%--------------------------------
\subsection{Alternative formula for $\xi$}\Label{SS:altcalB}
From \refe{calBalp=2} and \refe{zetB2=5} we obtain the following analogue of \refe{RaySingertor}
\eq{detB23}
     \xi(\n,g^M,\tet) \ = \  \frac12\,\sum_{k=0}^d\,(-1)^{k+1}\,k\,\LD_{2\tet}\Big[\,
    \Big(\, {(\Gam\,\n)^2} \ + \ {(\n\,\Gam)^2}\,\Big){\Big|_{\Ome^k(M,E)}}\,\Big],
\end{equation}
where $\tet\in (-\pi/2,0)$ is an Agmon angle for $B$ so that there are no eigenvalues of $B$ in the solid angles $L_{(-\pi/2,\tet]}$ and $L_{(\pi/2,\tet+\pi]}$. Note
that this condition implies that for all $k=0\nek d$, $2\tet$ is an Agmon angle for
 \(
    B^2{\big|_{\Ome^k(M,E)}}.
 \)
%-----------------------
\subsection{Choice of the spectral cut}\Label{SS:choicetet}
It follows from \refe{det-det} that $\Detgrtet(B_{\even})$ does not depend on the choice of the Agmon angle $\tet\in (-\pi,0)$. It is convenient
for us to work with an angle $\tet\in (-\pi/2,0)$ such that there are no eigenvalues of the operator $B$ in the solid angles
$L_{(-\pi/2,\tet]}$, $L_{[-\tet,\pi/2)}$, $L_{(\pi/2,\tet+\pi]}$ and $L_{[-\tet-\pi,-\pi/2)}$. We will fix such an angle till the end of this
section.

%-------------------------------------
\subsection{The dual connection}\Label{SS:n'}
Fix a Hermitian metric $h^{E}$ on $E$. Denote by $\n'$ the {\em connection on $E$ dual to the connection $\n$}. It is defined by the formula
\[
    dh^{E}(u,v) \ = \ h^{E}(\n u,v) \ + \ h^{E}(u,\n' v),
    \qquad u,v\in C^\infty(M,E).
\]
From the definition of the scalar product $\<\cdot,\cdot\>$ on $\Ome^\b(M,E)$ it then follows that
\eq{n*=2}
    \n^* \ = \ \Gam\,\n'\,\Gam,\qquad    (\n')^* \ = \ \Gam\,\n\,\Gam.
\end{equation}
Since  $\Gam^2=\Id$, \refe{n*=2} implies
\eq{(*na)*}
    \Big(\,{(\Gam\,\n)^2}\,\Big)^* \ = \ (\Gam\,\n')^2, \quad
     \Big(\,{(\n\,\Gam)^2}\,\Big)^* \ = \ (\n'\,\Gam)^2.
\end{equation}
Let $B'$ denote the odd signature operator associated to the connection $\n'$. Using \refe{Bpm} and \refe{n*=2} one readily sees that
\eq{Balp*-Balp'}\notag
    B^* \ = \ B'.
\end{equation}
Therefore, if the connection $\n$ satisfies Assumption~I and II of \refss{assumptions}, then so does the connection $\n'$. Our choice of the
angle $\tet$ in \refss{choicetet} guarantees that $\pm2\tet$ are Agmon angles for the operator
\[
    (\Gam\,\n')^2 \ = \ \big((\Gam\,\n)^2\big)^*.
\]
In particular, the number $\xi(\n',g^M,\tet)$ can be defined by  formula \refe{calB}, using the same angle $\tet$ and replacing everywhere $\n$
by $\n'$.

%%%
\lem{Balp-Balp'}
Using the notation introduced above, we have
\eq{Balp-Balp'}
    \xi(\n',g^M,\tet) \ \equiv \ \overline{\xi(\n,g^M,\tet)}, \qquad \ \mod\, \pi i,
\end{equation}
where $\oz$ denotes the complex conjugate of the number $z\in \CC$.
\elem
\prf
Set
\[
    D_k \ := \ \Big(\,{(\Gam\,\n)^2}+{(\n\,\Gam)^2}\,\Big){\Big|_{\Ome^k(M,E)}}:\,
    \Ome^k(M,E) \ \longrightarrow \ \Ome^k(M,E).
\]
Then, by \refe{(*na)*},
\[
    D_k^* \ = \ \Big(\,{(\Gam\,\n')^2}+{(\n'\,\Gam)^2}\,\Big){\Big|_{\Ome^k(M,E)}}:\,
    \Ome^k(M,E) \ \longrightarrow \ \Ome^k(M,E).
\]
With $\tet$ given as in \refss{choicetet}, we have
\eq{log=olog}
    \LD_{2\tet}\, D_k^* \ = \ \overline{\LD_{-2\tet}\,D_k}.
\end{equation}
Note that $D_k$ has a positive-definite leading symbol. Hence, at most  finitely many eigenvalues of $D_k$ lie in the solid angle
$L_{[-2\tet,2\tet+2\pi]}$ (which contains the negative real axis).%
\footnote{By our assumptions on $\tet$, all these eigenvalues must lie on the real axis. But we don't use this fact here.}
Hence, by \refe{zet-zet},
\eq{LDequivLD}
    \LD_{-2\tet} D_k \ \equiv \ \LD_{2\tet}D_k \qquad \MOD\ 2\pi i.
\end{equation}

Using \refe{log=olog} and \refe{LDequivLD}, we obtain now from \refe{detB23}, that
\[
    \xi(\n',g^M,\tet) \ = \
     \frac12\,\sum_{k=0}^d\,(-1)^{k+1}k\,\LD_{2\tet}\,D_k^*
    \ \ \equiv \  \
    \frac12\,\sum_{k=0}^d\,(-1)^{k+1}k\,\,{\overline{\LD_{2\tet}D_k}}
    \ = \ \overline{\xi(\n,g^M,\tet)}, \qquad  \MOD\ \pi i.
\]
\eprf

%%%
\lem{TRSalp-alp'}
For every\/ $\n\in \Flat'(M,g^M)$ we have\/ $\TRS(\n') = \TRS(\n)$.
\elem
\prf
From \refe{n*=2}, we obtain
\[
  \begin{aligned}
    \n^*\n \ &= \ \Gam\,\n'\,\Gam\,\n \ = \ \Gam\,(\,\n'\,\Gam\,\n\,\Gam\,)\,\Gam \ = \ \Gam\,\n'\,(\n')^*\,\Gam;\\
    \n\n^* \ &= \ \n\,\Gam\,\n'\,\Gam \ = \ \Gam\,(\,\Gam\,\n\,\Gam\,\n'\,)\,\Gam \ = \ \Gam\,(\n')^*\,\n'\,\Gam.
  \end{aligned}
\]
From  \refe{RaySingertor} we obtain
\eq{logTRS=sum}
 \begin{aligned}
    \log\,\TRS(\n) \ &= \
    \frac12\,\sum_{k=0}^d\,(-1)^{k+1}\,k\, \LD_{-\pi}\Big(\, (\n^*\n+\n\n^*){\big|_{\Ome^k(M,E)}}\,\Big)
    \\ &= \
    \frac12\,\sum_{k=0}^d\,
      (-1)^{k+1}\,k\, \LD_{-\pi}\Big(\, \Gam\,\big(\,\n'(\n')^*+(\n')^*\n'\,\big)\,\Gam{\big|_{\Ome^k(M,E)}}\,\Big)
    \\ &= \
    \frac12\,\sum_{k=0}^d\,(-1)^{k+1}\,k\,
      \LD_{-\pi}\Big(\, \big(\,\n'(\n')^*+(\n')^*\n'\,\big){\big|_{\Ome^{d-k}(M,E)}}\,\Big)
    \\ &= \
    \frac12\,\sum_{k=0}^d\,(-1)^{k}\,(d-k)\,
      \LD_{-\pi}\Big(\, \big(\,\n'(\n')^*+(\n')^*\n'\,\big){\big|_{\Ome^{k}(M,E)}}\,\Big).
 \end{aligned}
\end{equation}
By \refe{RaySingertor2},
\meq{sum=0}\notag
    \sum_{k=0}^d\,(-1)^{k}\, \LD_{-\pi}\Big(\, \big(\,\n'(\n')^*+(\n')^*\n'\,\big){\big|_{\Ome^{k}(M,E)}}\,\Big)
    \\ = \
    \sum_{k=0}^d\,(-1)^{k}\, \LD_{-\pi}\Big(\, {(\n')^*\n'}{\big|_{\Ome^k_+(M,E)}}\,\Big) \ + \
    \sum_{k=0}^d\,(-1)^{k}\, \LD_{-\pi}\Big(\, {\n'(\n')^*}{\big|_{\Ome^k_-(M,E)}}\,\big)
    \ = \ 0.
\end{multline}
Hence, from \refe{logTRS=sum}, we obtain
\[
    \log\,\TRS(\n) \ = \
    \frac12\,\sum_{k=0}^d\,(-1)^{k+1}\,k\,
       \LD_{-\pi}\Big(\, \big(\,\n'(\n')^*+(\n')^*\n'\,\big){\big|_{\Ome^{k}(M,E)}}\,\Big)
    \ = \ \log\,\TRS(\n').
\]
\eprf

%-------------------------------------
\subsection{Proof of \reft{DetB-TRS}}\Label{SS:prDetB-TRS}
In the case $\n$ is an acyclic Hermitian connection the statement has been already proved, cf. \refe{RaySingertor3}.

In the general case, let
\[
   \tiln \ = \ \begin{pmatrix}
                        \n&0\\
                        0&\n'
                       \end{pmatrix},
\]
denote the flat connection on $E\oplus{E}$ obtained as a direct sum of the connections $\n$ and $\n'$.  From Lemmas~\ref{L:Balp-Balp'} and
\ref{L:TRSalp-alp'} we obtain
\eq{alp+alp'}\notag
  \begin{aligned}
    \TRS(\tiln) \ &= \ \TRS(\n)\cdot\TRS(\n') \ = \ \big(\,\TRS(\n)\,\big)^2,\\
    \xi(\tiln,g^M,\tet) \ &= \ \xi(\n,g^M,\tet)\ + \ \xi(\n',g^M,\tet) \ \equiv \
    2\,\RE\,\xi(\n,g^M,\tet) \qquad \mod\, \pi i.
  \end{aligned}
\end{equation}
Hence, to prove \reft{DetB-TRS}, it is enough to show that
\eq{TRS=calB}
    \xi(\tiln,g^M,\tet) \ \equiv \ \log\,\TRS(\tiln) \qquad \MOD\, \pi i.
\end{equation}

We will prove \refe{alp+alp'} by a deformation argument. For $t\in [-\pi/2,\pi/2]$ introduce the rotation $U_t$ on
\[
    \Ome^\b \ := \ \Ome^\b(M,E)\,\oplus\,\Ome^\b(M,E)
\]
given by
\[
    U_t \ = \ \begin{pmatrix}
    \ \cos t&-\,\sin t \ \\
    \ \sin t& \ \ \ \cos t \
    \end{pmatrix}.
\]
Note that $U_t^{-1}= U_{-t}$.

Consider two one-parameter families of operators $\tilB(t), \hatB(t):\Ome^\b\to \Ome^\b$  ($t\in [-\pi/2,\pi/2]$):
\[
    \tilB(t) \ := \ \Gam\,U_t\,\tiln\,U_t^{-1} \ + \ \tiln\,\Gam;\qquad
    \hatB(t) \ := \ \Gam\,\tiln \ + \  U_t\, \tiln  U_t^{-1}\,\Gam.
\]
Note that $\tilB(0)= \hatB(0)= B(\tiln,g^M)$. If the Hermitian metric $h^E$ is invariant with respect to $\n$ then $\n'= \n$ and
\eq{tilB=hatB}
  \tilB(t) \ = \ \hatB(t) \ = \ B(\tiln,g^M) \ = \ B(\n,g^M)\oplus B(\n',g^M)
\end{equation}
for all $t\in  [-\pi/2,\pi/2]$. It follows then from the Assumption~II of \refss{assumptions} that the operator \refe{tilB=hatB} is invertible.

Suppose now that $\n$ is sufficiently close to an acyclic Hermitian connection $\n_0$ in the $C^0$-topology, cf. \refss{openset1}, and that the
metric $h^E$ is chosen to be invariant with respect to the connection $\n_0$. Then $\n'$ is also close to $\n_0$. Since both
$\tilB(t)-B(\tiln,g^M)$ and $\hatB(t)-B(\tiln,g^M)$ are 0'th order differential operators, it follows that they are small in the standard
operator norm (cf. proof of \refp{neigofh}) for all $t\in [-\pi/2,\pi/2]$. Therefore the operators $\tilB(t), \hatB(t)$ are invertible for all
$t\in [-\pi/2,\pi/2]$.

We denote by $V\subset \Flat(E)$ the set of connections, which satisfy the following property: there exists a Hermitian metric $h^E$ such
that the operators $\tilB(t)$ and $\hatB(t)$ are invertible for all $t\in [-\pi/2,\pi/2]$.%
\footnote{Recall that $\tilB(t)$ and $\hatB(t)$ depend on $h^E$ since the dual connection $\n'$ does.} Then $V$ is open in $\Flat(E)$. Moreover,
since for every $\n\in V$ the operator $\tilB(0)= B(\n,g^M)\oplus{}B(\n',g^M)$ is invertible, it follows that $V\subset \Flat'(E,g^M)$.

The above discussion shows that $V$ contains the set of acyclic Hermitian connections. In the rest of the proof we assume that $\n\in V$ and that $h^E$ is chosen so
that the operators $\tilB(t)$ and $\hatB(t)$ are invertible.

Set
\eq{Ome+-t}\notag
 \begin{aligned}
    \Ome^\b_+(t) \ &:= \ \Ker\, U_t \tiln  U_t^{-1}\,\Gam \ = \ \left\{\, \Gam\, U_t \,
    \left(\begin{smallmatrix}
        \ome\\ \ome'
    \end{smallmatrix}\right):\ \ome\in \Ker\n, \ \ome'\in \Ker\n'\,\right\};\\
        \Ome^\b_- \ &:= \ \Ker\tiln \ = \ \Ker \n\oplus \Ker\n'.
 \end{aligned}
\end{equation}
Note that $\Ome^\b_-$ is independent of $t$.

Since the range of $\Gam U_t \tiln  U_t^{-1}$ is contained in $\Ome^\b_+(t)$ whereas the range of $\tiln\,\Gam$ is contained in $\Ome^\b_-$, it follows from the
surjectivity of $\tilB(t)$ that
\eq{Omet+-1}
    \Ome^\b_+(t)\, +\, \Ome^\b_- \ = \ \Ome^\b,\qquad t\in [-\pi/2,\pi/2].
\end{equation}
Similarly, since, by definition, the kernel of\, $U_t \tiln U_t^{-1}\Gam$ is equal to  $\Ome^\b_+(t)$ whereas the kernel of $\Gam\,\tiln$ is equal to $\Ome^\b_-$, it
follows from injectivity of\, $\hatB(t)$ that
\eq{Omet+-2}
    \Ome^\b_+(t)\cap \Ome^\b_- \ = \ \{0\},\qquad t\in [-\pi/2,\pi/2].
\end{equation}

Combining \refe{Omet+-1} and \refe{Omet+-2} we obtain
\eq{Omet+-}
    \Ome^\b \ = \ \Ome^\b_+(t)\oplus \Ome^\b_-,\qquad t\in [-\pi/2,\pi/2].
\end{equation}

For each $t\in [-\pi/2,\pi/2]$ define $\xi(t)\in  \CC/\pi{i}\ZZ$ by the formula
\eq{calBt}
    \xi(t) \ \equiv \ \frac12\,\sum_{k=0}^d\,(-1)^k\LD_{\tet'}\Big(\,\Gam\, U_t\, \tiln\,  U_t^{-1}\,\Gam\,\tiln{\big|_{\Ome^k_+(t)}}\,\Big),
    \qquad \mod\, \pi i,
\end{equation}
where $\tet'\in L_{(-2\tet,2\pi+2\tet)}$ is any Agmon angle for the operators%
\footnote{Recall from \refss{det-tet} that a different choice of $\tet'\in  L_{(-2\tet,2\pi+2\tet)}$ changes the number
 \(
    \LD_{\tet'}(\,\Gam U_t \tiln  U_t^{-1}\Gam\tiln{\big|_{\Ome^k_+(t)}}\,)
 \)
by a multiple of $2\pi{i}$.}\ \
 \(
    \Gam U_t \tiln  U_t^{-1}\Gam\tiln{\big|_{\Ome^k_+(t)}}
 \)\
($k=0\nek N$).

Since
\[
    \Gam\, U_0\, \tiln\,  U_0^{-1}\,\Gam\,\tiln{\big|_{\Ome^k_+(0)}} \ = \
    \begin{pmatrix}
    \Gam\n\Gam\n\big|_{\Ome^k_+(M,E)}&0\\ 0 & \Gam\n'\Gam\n'\big|_{\Ome_+^k(M,E)}
    \end{pmatrix},
\]
for $t=0$, \refe{calBt} coincides with \refe{calB} with $\n$ replaced by $\tiln$.

Similarly, since
\[
    \Gam\,U_{\pi/2}\, \tiln\,  U_{\pi/2}^{-1}\,\Gam\,\tiln{\big|_{\Ome^k_+({\pi/2})}} \ = \
    \begin{pmatrix}
    \n^*\n\big|_{\Ome^k_+(M,E)}&0\\ 0 & \n'^*\n'\big|_{\Ome_+^k(M,E)}
    \end{pmatrix}
\]
for $t=\pi/2$ the right hand side of \refe{calBt} coincides with \refe{RaySingertor2}. Summarizing, we conclude that
\eq{calB0calBpi}
    \xi(0) \equiv \ \xi(\tiln,g^M,\tet), \quad
    \xi(\pi/2) \ \equiv \ \log\,\TRS(\tiln), \qquad \mod\ \pi i.
\end{equation}
We will finish the proof of \refe{TRS=calB} (and, hence, of \reft{DetB-TRS}) by showing that
\eq{ddtcalB}
    \frac{d}{dt}\,\,\xi(t) \ = \ 0.
\end{equation}

This is done by applying the arguments of the standard proof of the independence of the Ray-Singer torsion on the Hermitian metric. First we need the following
notation (cf., for example, Section~2 of \cite{BFK3}): Suppose $f(s)$ is a function of a complex parameter $s$ which is meromorphic near $s=0$. We call the zero order
term in the Laurent expansion of $f$ near $s=0$ the {\em finite part of $f$ at 0} and denote it by $\Fp f(s)$. Then, cf. Lemma~3.7 of \cite{BFK92} or formula (1.13)
of \cite{KontsevichVishik_long},
\meq{ddtdettilB}
    \frac{d}{dt}\,\LD_{\tet'} \Big(\, \Gam\, U_t\, \tiln\,  U_t^{-1}\,\Gam\,\tiln{\big|_{\Ome^k_+(t)}}\,\Big)
    \\ = \
    \Fp \Tr\,\left[\, { \left(\frac{d}{dt}\,\big(\Gam\, U_t\, \tiln\,  U_t^{-1}\,\Gam\,\tiln\big)\,\right)\,
    \left(\,\Gam\, U_t\, \tiln\,  U_t^{-1}\,\Gam\,\tiln\,\right)_{\tet'}^{-s-1}}{\big|_{\Ome^k_+(t)}}\,\right].
\end{multline}

One has
\[
    \frac{d}{dt}\,\big(\,\Gam\, U_t\, \tiln\,  U_t^{-1}\,\Gam\,\tiln\,\big) \ = \ \dot{U}_tU_t^{-1}\big(\Gam\, U_t \tiln  U_t^{-1}\Gam\,\tiln\big){\big|_{\Ome^k_+(t)}}
    \ - \ \big(\Gam\,U_t\tiln U_t^{-1}\Gam\,\big){\big|_{\Ome^{k+1}_-}}\dot{U}_tU_t^{-1}\tiln{\big|_{\Ome^{k}_+(t)}}
\]
By \refl{zetPQ=zetQP}, for $\RE{}s>d/2$,
\meq{onOme-}
    \Tr\, \left[\,\big(\Gam\,U_t\tiln U_t^{-1}\Gam\,\big){\big|_{\Ome^{k+1}_-}}\,\dot{U}_tU_t^{-1}\,\tiln{\big|_{\Ome^{k}_+(t)}}\,
       \left(\,\Gam\, U_t \tiln  U_t^{-1}\Gam\,\tiln\,\right)_{\tet'}^{-s-1}{\big|_{\Ome^k_+(t)}}\,\right]
    \\ = \
    \Tr\,\left[\, \dot{U}_tU_t^{-1}\tiln{\big|_{\Ome^{k}_+(t)}}\,
       \left(\,\Gam\, U_t \tiln  U_t^{-1}\Gam\,\tiln\,\right)_{\tet'}^{-s-1}{\big|_{\Ome^k_+(t)}}\,
            \big(\Gam\,U_t\tiln U_t^{-1}\Gam\,\big){\big|_{\Ome^{k+1}_-}}\,\right]
    \\ = \
    \Tr\,\left[\, \dot{U}_tU_t^{-1}\,\big(\, \tiln\,\Gam\,U_t\tiln\,U_t^{-1}\Gam\,\big)\big|_{\Ome^{k+1}_-}\,\big)_{\tet'}^{-s}\,\right]
\end{multline}
Hence, \refe{ddtdettilB} implies that
\meq{ddtdettilB2}
    \frac{d}{dt}\,\LD_{\tet'} \Big(\, \Gam\, U_t \tiln  U_t^{-1}\Gam\,\tiln{\big|_{\Ome^k_+(t)}}\,\Big)
    \\ = \
    \Fp \Tr\,\Big[\,\dot{U}_tU_t^{-1}\,\Big(\,\Gam\, U_t \tiln  U_t^{-1}\Gam\,\tiln{\big|_{\Ome^k_+(t)}}\,\Big)_{\tet'}^{-s}
    \, - \,
    \dot{U}_tU_t^{-1}\,\Big(\,\tiln\Gam\, U_t \tiln  U_t^{-1}\Gam{\big|_{\Ome^{k+1}_-}}\,\Big)_{\tet'}^{-s}\,\Big]
\end{multline}

Consider the operator
\eq{tilLap}
    \tilDel_k(t) \ := \ \Gam\, U_t\, \tiln\,  U_t^{-1}\,\Gam\,\tiln\,{\big|_{\Ome^k_+(t)}} \ + \
    \tiln\,\Gam\, U_t\,\tiln \, U_t^{-1}\,\Gam\,{\big|_{\Ome^{k}_-}}.
\end{equation}
It is a second order elliptic differential operator on\/ $\Ome^k(M,E\oplus{}E)$, whose leading symbol is equal to the leading symbol of the Laplacian
$\tiln^*\,\tiln+\tiln\tiln^*\,$. In other words, $\tilDel_k(t)$ is a {\em generalized Laplacian} in the sense of \cite{BeGeVe}.

The decomposition \/ $\Ome^k(M,E\oplus{}E) \ = \ \Ome^k_+(t)\,\oplus\,\Ome^k_-$\, implies that
\eq{tilDels2}\notag
    \tilDel_k(t)_{\tet'}^{-s} \ = \ \Big(\,\Gam\, U_t\, \tiln\, U_t^{-1}\,\Gam\,\tiln\,{\big|_{\Ome^k_+(t)}}\,\Big)_{\tet'}^{-s} \ + \
    \Big(\,\tiln\,\Gam\, U_t\, \tiln\,  U_t^{-1}\,\Gam\,{\big|_{\Ome^{k}_-}}\,\Big)_{\tet'}^{-s}.
\end{equation}

Hence, from \refe{ddtdettilB2}, we obtain
\eq{ddtdettilB3}
    \frac{d}{dt}\,\xi(t) \ = \ \frac12\,\sum_{k=0}^d\,(-1)^k\, \Fp \Tr\,\Big[\,\dot{U}_tU_t^{-1}\,\tilDel_k(t)_{\tet'}^{-s}\,\Big].
\end{equation}
By a slight generalization of a result of Seeley \cite{Seeley67}, which is discussed in \cite{Wodzicki87}, the right hand side of
\refe{ddtdettilB3} is given by a {\em local formula}, i.e., by an integral
\eq{intphi}
   \int_M\, \phi
\end{equation}
of a differential form $\phi$, whose value at a point $x\in M$ depends only on the full symbol of $\tilDel$ and a finite number of its
derivatives at the point $x$. Moreover, since the dimension of the manifold $M$ is odd, the differential form $\phi$ vanishes identically.
\hfill $\square$

%-------------------------------------------------------------------
%-------------------------------------------------------------------
\section{Dependence of the Graded Determinant on the Riemannian Metric}\Label{S:deponmet}

As already mentioned, one can consider the graded determinant $\Detgrtet(B_{\even})$, defined in \refe{grdetB}, as a refinement of the Ray-Singer torsion. However, in
general, $\Detgrtet(B_{\even})$ depends on the choice of the Riemannian metric $g^M$ on $M$. In this section we investigate this dependence. In particular, we show
that, if $\dim M= 2r-1\equiv 1 \ (\MOD\ 4)$, then $\Detgrtet(B_{\even})$ is independent of $g^M$. Later we will use the results of this section to construct a
refinement of the Ray-Singer torsion which is a diffeomorphism invariant of the pair $(E,\n)$ (i.e. is independent of the metric).

\defe{calBn}
A Riemannian metric $g^M$ on $M$ is called {\em admissible for a given acyclic connection $\n$} if the odd signature operator $B= B(\n,g^M)$ satisfies Assumption~II
of \refss{assumptions}. We denote the set of admissible metrics by $\calM(\n)$.
\edefe

Recall from \reft{DetB-eta2} that, for $g^M\in \calM(\n)$,
\eq{DetB-eta22}
    \LDgrtet\big(\,B_{\even}(\n,g^M)\,\big) \ = \ \xi(\n,g^M,\tet) \ - \  i\pi\,\eta(\n,g^M),
\end{equation}
where $\tet\in (-\pi/2,0)$ is an Agmon angle for $B$ such that there are no eigenvalues of the operator $B$ in the solid angles $L_{(-\pi/2,\tet]}$ and
$L_{(\pi/2,\tet+\pi]}$.

%----------------------------------------
\subsection{Dependence of the $\eta$-invariant on the metric}\Label{SS:etamet}
First, we study the dependence of the $\eta$-invariant $\eta= \eta(B_{\even}(\n,g^M))$ on the metric $g^M$. Fortunately, this was essentially done in \cite{APS2} and
\cite{Gilkey84}. Below we present a brief review of the relevant results.

Let $B_{\text{trivial}}:\Ome^{\even}(M)\to \Ome^{\even}(M)$ denote the even part of the odd signature operator corresponding to the trivial line bundle over $M$
endowed with the trivial connection. It is shown on page~52 of \cite{Gilkey84} (see also Theorem~2.4 of \cite{APS2} where the case of unitary connection is
established) that modulo $\ZZ$ the difference
\eq{tileta}
    \eta\big(B_{\even}(\n,g^M)\big) \ - \ (\rank E)\ \eta\big(B_{\text{trivial}}(g^M)\big)
\end{equation}
is independent of the Riemannian metric.

Let us describe the dependence of $\eta\big(B_{\text{trivial}}(g^M)\big)$ on the metric.

%---------------------------------------------
\subsubsection{Case when $M$ bounds an oriented manifold $N'$}\Label{SS:Mbounds}
Suppose, first, that $M$ is the oriented boundary of a smooth compact oriented manifold $N'$. Let $\sign(N')$ denote the signature of $N'$, cf. \cite{APS1}. This is
an integer defined in purely cohomological terms. In particular, it is independent of the metric. The signature theorem for manifolds with boundary (cf. Theorem~4.14
of \cite{APS1} and Theorem~2.2 of \cite{APS2}) states that
\eq{signth}
    \sign(N') \ = \ \int_{N'}\, L(p)  \ - \ \eta(B_{\text{trivial}}),
\end{equation}
where $L(p):= L_{N'}(p)$ is the Hirzebruch $L$-polynomial in the Pontrjagin forms of a Riemannian metric on $N'$ which is a product near $M$. It follows from
\refe{signth} that $\int_{N'}L(p)$ is independent of the choice of the Riemannian metric on $N'$ among those that are product-like near $\partial{N'}$. Note that if
$\dim{M}\equiv 1 \ (\MOD\, 4)$ then $L(p)$ does not have a term of degree $\dim{N'}$ and, hence, $\int_{N'}\,L(p)=0$.

Combining \refe{signth} with the metric independence of\, $\sign(N')$ and \refe{tileta}, we conclude that, modulo $\ZZ$,
\eq{eta-L}
    \eta \ - \ (\rank E)\, \int_{N'}\, L(p)
\end{equation}
is independent of the metric $g^M$.  Since for different choices of $N'$, the integral $\int_{N'}L(p)$ differs by an integer, the expression \refe{eta-L}, modulo
$\ZZ$, is also independent of $N'$.

%---------------------------------------------
\subsubsection{General case ($M$ does not necessarily bound an oriented manifold)}\Label{SS:Mnotbounds}
In general, there might be no smooth oriented manifold whose oriented boundary is diffeomorphic to $M$. However, since $\dim{}M$ is odd, there exists an oriented
manifold $N$ whose oriented boundary is the disjoint union of two copies of $M$ (with the same orientation), cf. \cite{Wall60}, \cite[Th.~IV.6.5]{Rudyak_book}. Then
the same arguments as above show that, modulo $\ZZ$,
\eq{eta-L2}
    \eta \ - \ \frac{\rank E}2\, \int_{N}\, L(p)
\end{equation}
is metric independent. In particular, if $\dim{M}\equiv 1 \ (\MOD\, 4)$, then the reduction of $\eta$ modulo $\ZZ$ is metric independent.

%-------------
\rem{eta-L} Note again that replacing $\eta$ by \refe{eta-L2} removes the dependence on the metric but creates a new dependence on the
choice of the manifold $N$. For different choices of $N$ the integrals $\int_{N}\, L(p)$ might differ by an integer.
\erem

%------------------------
The following proposition summarizes our discussion of the dependence of the $\eta$-invariant on the Riemannian  metric.
\prop{deponmetr}
Let $N$ be a smooth compact oriented manifold whose oriented boundary is equal to two disjoint copies of $M$. Denote by $L(p)= L_N(p)$ the Hirzebruch $L$-polynomial
in the Pontrjagin forms of a Riemannian metric on $N$ which is a product near $\partial{N}= M\sqcup{M}$. Then, modulo $\ZZ$,
\[
    \eta\,-\, \frac{\rank E}2\,\int_{N}\,L(p)
\]
is independent of the choice of the metric on $M$. For different choices of $N$ satisfying $\partial{N}=M\sqcup{M}$ the integral $\int_{N}L(p)$ differs by an integer.
If\, $\dim{M}\equiv1 (\MOD\,4)$, then $\int_{N}L(p)=0$.

The imaginary part\/ $\IM\eta$ of\/ $\eta$ is independent of\/ $g^M$.
\eprop

%--------------------------------
We are now ready to formulate the main result of this section.

%-------
\th{grdetmet}
Let $E$ be a flat vector bundle over a closed oriented odd-dimensional manifold $M$ and let $\n$ be the flat connection on $E$. Let $N$ be an oriented manifold whose
oriented boundary is the disjoint union of two copies of $M$. For each admissible Riemannian metric $g^M\in \calM(\n)$ consider the number
\eq{metricindep}
    \Detgrtet\big(\,B_\even(\n,g^M)\,\big)\cdot e^{i\pi\,\frac{\rank E}2\,\int_N\,L(p,g^M)}\in \CC\backslash \{0\},
\end{equation}
where $\tet\in (-\pi/2,0)$ is an Agmon angle for $B_\even(\n,g^M)$ and $L(p,g^M)$ is the Hirzebruch $L$-polynomial in the Pontrjagin forms of any Riemannian metric on
$N$ which near $M$ is the product of $g^M$ and the standard metric on the half-line. Then the number \refe{metricindep} is independent of\/ $g^M\in \calM(\n)$ and
$\tet\in (-\pi/2,0)$.

In particular, if\/ $\dim{}M\equiv 1 (\MOD 4)$, then $\int_NL(p,g^M)=0$ and, hence, $\Detgrtet\big(B_\even(\n,g^M)\big)$ is independent of $g^M$.
\eth

To prove \reft{grdetmet} we need to study the dependence of\/ $\xi=\xi(\n,g^M,\tet)$ on $g^M$.

%---------------------------------
\subsection{Dependence of $\xi$ on the Riemannian metric}\label{SS:xi(t)}
For \refe{DetB-eta22} to hold we need to assume that there are no eigenvalues of the operator $B$ in the solid angles $L_{(-\pi/2,\tet]}$ and $L_{(\pi/2,\tet+\pi]}$.
However, for the study of the dependence of $\xi$ on $g^M$ it will be convenient for us to work with $\xi(\n,g^M,\tet)$ with the only assumption that both, $\tet\in
(-\pi,0)$ and $\tet+\pi$, are Agmon angles for $B_{\even}$. If $\tet_1$ and $\tet_2$ are two such angles then, by \refe{zet-zet} and \refe{calB},
\eq{xitet1tet2}
    \xi(\n,g^M,\tet_1) \ \equiv \ \xi(\n,g^M,\tet_2) \qquad \MOD\ \pi i.
\end{equation}

%---------------------
\prop{detB2metric}
Suppose $g^M_0,\ g^M_1\in \calM(\n)$ are admissible Riemannian metrics on\/ $M$ and let $\tet_0, \tet_1\in (-\pi/2,0)$  be such that, for $j=0,1$, both, $\tet_j$ and
$\tet_j+\pi$ are Agmon angles for $B(\n,g^M_j)$. Then
\eq{xitet1=xitet2}
    \xi(\n,g^M_1,\tet_1) \ \equiv \ \xi(\n,g^M_0,\tet_0) \qquad \MOD\ \pi i.
\end{equation}
\eprop
\rem{detB2-RS}
If, in addition, $\n$ is Hermitian and $h^{E}$ is a $\n$-invariant Hermitian metric on $E$, then, from \refe{n*=} we obtain
 \(
    (\Gam\,\n)^2 \ = \ \n^*\,{\n}.
 \)
By \refe{RaySingertor2}, $\xi(t,\tet_0)$ coincides, in this case, with the Ray-Singer torsion. Hence, in the case of a Hermitian connection, the statement of the
proposition reduces to the classical result about the independence of the Ray-Singer torsion on the Riemannian metric.
\erem

For the proof of the proposition we, first, consider the case when $g^M_0$ and $g^M_1$ belong to the same path-connected component of the set $\calM(\n)$ of
admissible metrics.

Suppose that $g^M_t\in \calM(\n)$, $t\in\RR$, is a smooth family of admissible Riemannian metrics on $M$ and let $B_t= B(\n,g_t^M)$  be the corresponding odd
signature operator.  To simplify the notation set
\eq{xi(t)eta(t)}\notag
    \xi(t,\tet) \ := \ \xi(\n,g^M_t,\tet).
\end{equation}

Fix $t_0\in \RR$ and let $\tet_0\in (-\pi/2,0)$ be an Agmon angle for $B_{t_0}$ such that there are no eigenvalues of $B_{t_0}$ in the solid angles
$L_{(-\pi/2,\tet_0]}$ and $L_{(\pi/2,\tet_0+\pi)}$. Choose $\del>0$ so that for every $t\in (t_0-\del,t_0+\del)$  both, $\tet_0$ and $\tet_0+\pi$, are Agmon angles of
$B_t$. For $t\not= t_0$ it might happen that there are eigenvalues of $B_t$ in $L_{(-\pi/2,\tet_0)}$ and/or $L_{(\pi/2,\tet_0+\pi)}$. Hence, \refe{DetB-eta22} is not
necessarily true, in general, for $t\not=t_0$. However, from \refe{xitet1tet2}, we conclude that for every $t\in (t_0-\del,t_0+\del)$ and  $\tet\in (-\pi/2,0)$, such
that $\tet$ and $\tet+\pi$ are Agmon angles for $B_t$,
\eq{xitettet'}
    \xi(t,\tet) \ \equiv \ \xi(t,\tet_0) \qquad \MOD\ \pi i,
\end{equation}

\lem{detB2metric}
Under the above assumptions, $\xi(t,\tet_0)$ is independent of\/ $t\in (t_0-\del,t_0+\del)$.
\elem
\prf
Let $\Gam_t$ denote the chirality operator corresponding to the metric $g_t^M$. Then, cf. Lemma~3.7 of \cite{BFK92} or formula (1.13) of \cite{KontsevichVishik_long},
\eq{ddtdet}
    \frac{d}{dt}\,\LD_{2\tet_0} \Big(\, {(\Gam_t\n)^2}{\big|_{\Ome^k_+(M,E)}}\,\Big)
    \ = \
    \Fp \Tr\,\Big[\, { \big(\frac{d}{dt}\,(\Gam_t\n)^2\,\big)\,\big(\,(\Gam_t\n)^2\,\big)_{2\tet_0}^{-s-1}}{\big|_{\Ome^k_+(M,E)}}\,\Big],
\end{equation}
where we use the notation $\Fp$ introduced in \refss{prDetB-TRS}.

We denote by $\dot\Gam_t$ the derivative of $\Gam_t$ with respect to the parameter $t$. Then
\eq{ddt*n}
    \frac{d}{dt}\,{(\Gam_t\n)^2}{\big|_{\Ome^k_+(M,E)}} \ = \
    \dot\Gam_t\, \Gam_t\,{(\Gam_t\n)^2}{\big|_{\Ome^k_+(M,E)}} \ + \
    (\Gam_t\n){\big|_{\Ome^{d-k-1}_+(M,E)}}\,\dot\Gam_t\,\Gam_t\,(\Gam_t\n){\big|_{\Ome^k_+(M,E)}},
\end{equation}
where we used that $\Gam_t^2=1$. Using \refe{ddt*n} and the equality $\Tr AB=\Tr BA$, we obtain from \refe{ddtdet} that
\meq{ddtdet2}
    \frac{d}{dt}\,\LD_{2\tet_0} \Big(\, {(\Gam_t\,\n)^2}{\big|_{\Ome^k_+(M,E)}}\,\Big)
    \\ = \
    \Fp \Tr\,\Big[\, \dot\Gam_t\Gam_t\,\big(\,{(\Gam_t\,\n)^2}{\big|_{\Ome^k_+(M,E)}}\,\big)_{2\tet_0}^{-s} +
    \dot\Gam_t\,\Gam_t\,\big(\,{(\Gam_t\,\n)^2}{\big|_{\Ome^{d-k-1}_+(M,E)}}\,\big)_{2\tet_0}^{-s}\,\Big].
\end{multline}
Hence,
\meq{ddtsumdet}
    \frac{d}{dt}\,\sum_{k=0}^d\,(-1)^k\LD_{2\tet_0} \Big(\,{(\Gam_t\,\n)^2}{\big|_{\Ome^k_+(M,E)}}\,\Big)
    \\ = \
    2\, \sum_{k=0}^d\,(-1)^k\,\Fp \Tr\,\Big[\, \dot\Gam_t\Gam_t\,\big(\,{(\Gam_t\,\n)^2}{\big|_{\Ome^k_+(M,E)}}\,\big)_{2\tet_0}^{-s}\,\Big].
\end{multline}

Similarly,
\meq{ddtsumdet2}
    \frac{d}{dt}\, \sum_{k=0}^d\,(-1)^{k-1}\LD_{2\tet_0} \Big(\, {(\n\Gam_t)^2}{\big|_{\Ome^{k}_-(M,E)}}\,\Big)
    \\ = \
    2\, \sum_{k=0}^d\,(-1)^{k-1}\,\Fp \Tr\,\Big[\, \Gam_t\dot\Gam_t\,\big(\,{(\n\Gam_t)^2}{\big|_{\Ome^{k}_-(M,E)}}\,\big)_{2\tet_0}^{-s}\,\Big].
\end{multline}
From \refe{calB}, we see that \refe{ddtsumdet} is equal to $2\frac{d}{dt}\xi(t,\tet_0)$. By \refe{zet*n=zetn*}, the left hand sides of \refe{ddtsumdet} and
\refe{ddtsumdet2} are equal.  Hence \refe{ddtsumdet2} is also equal to $2\frac{d}{dt}\xi(t,\tet_0)$. We conclude that
\meq{ddtsum=ddtsum}
   \frac{d}{dt}\,\xi(t,\tet_0) \ = \  \sum_{k=0}^d\,(-1)^k\,\Fp \Tr\,\Big[\, \dot\Gam_t\Gam_t\,\big(\,{(\Gam_t\,\n)^2}{\big|_{\Ome^k_+(M,E)}}\,\big)_{2\tet_0}^{-s}\,\Big]
    \\ = \
    \sum_{k=0}^d\,(-1)^{k-1}\Fp \Tr\,\Big[\, \Gam_t\dot\Gam_t\,\big(\,{(\n\Gam_t)^2}{\big|_{\Ome^{k}_-(M,E)}}\,\big)_{2\tet_0}^{-s}\,\Big].
\end{multline}
Hence,
\eq{dddtxi}
    2\,\frac{d}{dt}\,\xi(t,\tet_0) \ = \
      \sum_{k=0}^d\,(-1)^k\,\Fp \Tr\,\Big[\, \dot\Gam_t\Gam_t\,\big(\,{(\Gam_t\,\n)^2}{\big|_{\Ome^k_+(M,E)}}\,\big)_{2\tet_0}^{-s}
      \,-\,  \Gam_t\dot\Gam_t\,\big(\,{(\n\Gam_t)^2}{\big|_{\Ome^k_-(M,E)}}\,\big)_{2\tet_0}^{-s}\,\Big]
\end{equation}

Since $\Gam_t^2=1$, we obtain $\dot\Gam_t\Gam_t + \Gam_t\dot\Gam_t= \frac{d}{dt}\Gam_t^2 = 0$. Hence, \refe{dddtxi} can be rewritten as
\eq{ddtsumdetfinal}
    2\,\frac{d}{dt}\,\xi(t,\tet_0)
    \ = \
    \sum_{k=0}^d\,(-1)^{k}\,\Fp \Tr\,\Big[\, \dot\Gam_t\Gam_t\,\tilDel_k(t)_{2\tet_0}^{-s}\,\Big],
\end{equation}
where $\tilDel_k(t)= (\Gam_t\n)^2 +  (\n\Gam_t)^2$   (k=0\nek d).
%Note that, if $\n$ is Hermitian we conclude from \refe{n*=} that $\tilDel_k(t)$ coincides with the
%Laplacian $\n^*\n+\n\n^*$. In the general case, where $\n$ is not necessarily Hermitian,  $\tilDel_k(t)$ is a second order elliptic differential operator, whose
%leading symbol is equal to the leading symbol of the Laplacian, i.e., $\tilDel_k(t)$ is a {\em generalized Laplacian} in the sense of \cite{BeGeVe}.
By a slight generalization of a result of Seeley \cite{Seeley67}, which is discussed in \cite{Wodzicki87}, the right hand side of \refe{ddtsumdetfinal} is given by a
{\em local formula}, i.e., by the integral \refe{intphi} of a differential form $\phi_t$, whose value at any point $x\in M$ depends only on the values of the
components of the metric tensor $g_t^M$ and a finite number of their derivatives at $x$. Moreover, since the dimension of the manifold $M$ is odd, the differential
form $\phi_t$ vanishes identically.  Hence, $\frac{d}{dt}\xi(t,\tet_0)=0$ for all $t\in (t_0-\del,t_0+\del)$.
\eprf

%-------------------------------------------------------
\subsection{Proof of \refp{detB2metric}}\Label{SS:prdetB2metric2}
Set
\[
  g_t^{M} \ = \ (1-t)\,g^{M}_0 \ + \ t\,g^{M}_1, \qquad t\in [0,1],
\]
and let $\Gam_t$ denote the chirality operator corresponding to the metric $g^M_t$. The operators $\Gam_t$ depend real analytically on $t$ and we can extend their
definition to all $t$ in some connected open connected neighborhood $U\subset \CC$ of $[0,1]$. Hence, the operator
\[
    B_t \ := \ \Gam_t\,\n\ +\ \n\,\Gam_t
\]
is well defined for all $t\in U$ and is holomorphic on $U$ in the sense of \refss{holcurveDiff}.  If we choose the neighborhood $U$ to be small enough then $B_t^2\in
\Ell_{m,(-3\pi/4,-\pi/4)}(M,E)$ for all $t\in U$. By \refc{zeroanal} the set
\[
  \Sig \ := \ \big\{\,t\in U:\,B_t \ \ \text{is not invertible}\,\big\}
\]
is a complex analytic subset of $U$. Thus $U\backslash\Sig$ is connected. Since the metrics $g_0^M$ and $g_1^M$ are admissible, it follows that $0$ and $1$ are in
$U\backslash\Sig$. For $\tet\in (-\pi/2,0)$  such that both, $\tet$ and $\tet+\pi$, are Agmon angles for $B_t$, set
\eq{xit}
   \xi(t,\tet)\ := \  \ \frac12\,\sum_{k=0}^d\, (-1)^{k+1}\,k\,\LD_{2\tet}\,\Big[\,B_t^2\big|_{\Ome^k(M,E)}\,\Big].
\end{equation}
If $t\in [0,1]$ is such that the metric $g^M_t$ is admissible, then $\xi(t,\tet)= \xi(\n,g^M_t,\tet)$.

By \reft{dethol}.a, the function $t\mapsto \xi(t,\tet)$ is holomorphic on the open set
\[
    U_\tet \ := \ \big\{\,t\in U\backslash\Sig:\ \ \tet,\ \tet+\pi \ \text{are Agmon angles for}\  B_t\,\big\}.
\]
By \refe{xitettet'}, if $t\in U_{\tet_1}\cap U_{\tet_2}$ then $\xi(t,\tet_1)\equiv \xi(t,\tet_2)$ modulo $\pi\ZZ$. Hence, we can define a multivalued analytic
function on $U\backslash\Sig$ by the formula
\[
    t \ \mapsto \ \xi(t,\tet_t)\ + \ \pi \ZZ,
\]
where $\tet_t\in (-\pi/2,0)$ is any angle such that $t\in U_{\tet_t}$.

Since the metric $g_0^M$ is admissible there exists $\eps>0$ such that for all real $t\in [0,\eps]$ the metric $g_t^M$ is admissible and $\tet_0$ and $\tet_0+\pi$ are
Agmon angles for $B_t$. Hence, by \refl{detB2metric}, the holomorphic function $\xi(t,\tet_0)$ is constant on $[0,\eps]$. Thus, since the set $U\backslash\Sig$ is
connected, our multivalued analytic function $t\mapsto \xi(t,\tet)$ is constant on $U\backslash\Sig$. \hfill$\square$

%-------------------------------------------------------
\subsection{Proof of \reft{grdetmet}}\Label{SS:prgrdetmet}
The fact that \refe{metricindep} is independent of $\tet$ follows immediately from \refe{zet-zet}. Let us prove that it is independent of $g^M\in \calM(\n)$. Suppose
$g^M_0$ and $g^M_1$ are admissible metrics. We shall use the notation introduced in \refss{prdetB2metric2}. For $t\in U\backslash\Sig$, fix $\tet_t\in (-\pi/2,0)$
such that there are no eigenvalues of $B_t$ in the solid angles $L_{(-\pi/2,\tet_t]}$ and $L_{(\pi/2,\tet_t+\pi)}$. As $t$ is not necessarily real, in general, $B_t$
is not an odd signature operator associated to a Riemannian metric. Hence, to calculate $\LD_{\operatorname{gr},\tet_t}(B_t)$ we can not use \reft{DetB-eta2}.
However, a verbatim repetition of the proof of this theorem shows that
\[
    \LD_{\operatorname{gr},\tet_t}\big(\,B_t\,\big) \ = \ \xi(t,\tet_t) \ - \  i\pi\,\eta(B_t),
\]
where $\xi(t,\tet_t)$ is defined by \refe{xit}. It follows now from Propositions~\ref{P:deponmetr} and \ref{P:detB2metric} that
\meq{0equiv1}
    \Detgrteto\big(\,B_\even(\n,g^M_0)\,\big)\cdot e^{i\pi\,\frac{\rank E}2\,\int_N\,L(p,g^M_0)} \\ = \
    \pm\,\Detgrtetone\big(\,B_\even(\n,g^M_1)\,\big)\cdot e^{i\pi\,\frac{\rank E}2\,\int_N\,L(p,g^M_1)},
\end{multline}
where $\tet_j$ $(j=0,1)$ is an Agmon angle for $B_\even(\n,g^M_j)$.

Since the function
\[
    t\ \mapsto\ \Det_{\operatorname{gr},\tet_t}\big(\,B_\even(\n,g^M_t)\,\big)\cdot e^{i\pi\,\frac{\rank E}2\,\int_N\,L(p,g^M_t)}
\]
(which is independent of $\tet_t$) is continuous on the connected set $U\backslash\Sig$, the sign in the right hand side of \refe{0equiv1} must be positive. The
theorem is proven. \hfill$\square$

%-------------------------------------------------------------------
%-------------------------------------------------------------------
\section{Refined Analytic Torsion in the Case $\dim M \equiv 1\ (\text{mod}\ 4)$}\Label{S:equiv3}

By \reft{grdetmet}, if $\dim M= 2r-1\equiv 1 \ (\text{mod}\ 4)$, then $\Detgrtet(B_{\even})$ is independent of the choice of the Riemannian metric on $M$. This
justifies the following %%
\defe{refantor}
Let $M$ be  a closed oriented manifold of dimension \/$\dim M= 2r-1\equiv 1 \ (\MOD\ 4)$. Suppose that there exists a Riemannian metric $g^M$
such that $\n\in \Flat'(M,g^M)$. The {\em refined analytic torsion} $T(\n)$ is defined to be the graded determinant of the odd signature
operator,
\eq{refantor}
    T(\n) \ = \ T(M,E,\n) \ := \ \Detgrtet(B_{\even})  \ \in \ \CC\backslash\{0\},
\end{equation}
where $\tet\in (-\pi,0)$ is an Agmon angle for the operator $B_{\even}= B_{\even}(\n,g^M)$.
\edefe

Remark that $T(\n)\not=0$. Note also that $T(\n)$ does not depend on the choice of $\tet\in (-\pi,0)$, cf. \refe{det-det}. Further, by
\refc{DetB-TRS}, if $\n$ is a Hermitian connection then
\eq{T=TRS}
    |T(\n)| \ = \ \TRS(\n).
\end{equation}

%------------------------------
\subsection{Example}\Label{SS:example}
We now calculate the refined analytic torsion in the simplest possible example, when $M=\RR/2\pi\ZZ$ is the circle and $E= M\times\CC$ is the
trivial line bundle over $M$. Fix a number $a\in \CC\backslash\ZZ$ and define the connection $\n_a$ on $E$ by the formula
\[
    \n_a:\, f \ \mapsto \ df+iafdx, \qquad f\in \Ome^0(M,E)=\Ome^0(M),
\]
where $x\in [0,2\pi)$ is the coordinate on $M$. According to the formula \refe{oddsign}, the odd signature operator is
\[
    B_{\even} \ = \ B_{0}:\, f \ \mapsto \ -i*\n_a f  \ = \ -if'+af.
\]
As $a\in \CC\backslash\ZZ$, the operator $B_{0}:\Ome^0(M)\to \Ome^0(M)$ is invertible and its eigenvalues are given by $a+n$, $n\in \ZZ$. Since
$\RE{a}\in \RR\backslash\ZZ$, the angle $\tet=-\pi/2$ is an Agmon angle for $B_{0}$.

The refined analytic torsion $T(\n_a)=\Det_{\tet}(B_{0})$  can be easily calculated using, for example, the general formula for determinants of
elliptic operators on the circle, obtained in \cite{BFK91} (see also \cite{Wojciechowski99} for an alternative way of calculation). As a result
we obtain the following formula for the refined analytic torsion
\[
    T(\n_a) \ = \ 1 \ - \ e^{2a\pi i} \ = \ 2\, (\sin\pi a)\, e^{i\frac{\pi}2(2a-1)}.
\]
Note that if $a\in \RR\backslash\ZZ$ then $\n_a$ is a Hermitian connection.  We conclude that, even for a Hermitian connection, the refined
analytic torsion is a complex number, which, depending on the value of $a$, can have an arbitrary phase, aside from $\pm\pi/2$.

%-------------------------------------------------------------------
%-------------------------------------------------------------------
\section{Refined Analytic Torsion in the Case $\dim M \equiv 3\ (\text{mod}\ 4)$}\Label{S:equiv1}

If $\dim M \equiv 3\ (\text{mod}\ 4)$ then $\Detgrtet(B_{\even})$ depends on the Riemannian metric on $M$. However, \reft{grdetmet} allows to define a refinement of
the Ray-Singer torsion, which is independent of the choice of an (admissible) metric on $M$.

\defe{refantor2}
Let $M$ be a closed oriented manifold of dimension $\dim M\equiv 3 (\MOD 4)$. Assume that there exists a Riemannian metric $g^M$ on $M$ such that
$\n\in \Flat'(M,g^M)$. Let $\tet\in (-\pi,0)$ be an Agmon angle for $B_{\even}= B_{\even}(\n,g^M)$. Choose a smooth compact oriented manifold $N$
whose oriented boundary is diffeomorphic to two disjoint copies of $M$. Then the {\em refined analytic torsion} $T(\n)$ is the non-vanishing complex
number defined by the formula
\begin{equation}\Label{E:refantor3}
    T(\n) \ = \ T(M,E,\n,N) \ := \
        \Detgrtet(B_{\even})\cdot \exp\Big(\,i\pi\,\frac{\rank E}2\,  \int_{N}\, L(p)\,\Big),
\end{equation}
where $L(p)= L_N(p)$ is the Hirzebruch $L$-polynomial in the Pontrjagin forms of a Riemannian metric on $N$ which is a product near
$\partial{N}$.
\edefe

\rem{TdependsonN}
Note that the refined torsion is independent of the angle $\tet\in (-\pi,0)$ and of the metric. But {\em it does depend on the choice of the manifold $N$}. However,
from \refp{deponmetr}, we conclude that $T(\n)$ is independent of the choice of $N$ up to multiplication by $i^{k\cdot\rank{E}}$ $(k\in \ZZ)$. If $\rank{}E$ is even
then $T(\n)$ is well defined up to a sign, and if $\rank{}E$ is divisible by 4, then  $T(\n)$ is a well defined complex number.

(Here a quantity being well defined means that it depends only on $M$, $E$ and $\n$.)
\erem

\rem{Twhenbounds}
If $M$ is the oriented boundary of a smooth compact oriented manifold $N'$, one can define a version of the refined analytic torsion:
\begin{equation}\Label{E:refantor2}
    T'(\n) \ = \ T'(M,E,\n,N') \ := \
        \Detgrtet(B_{\even})\cdot \exp\Big(\,i\pi\cdot\rank{}E\,  \int_{N'}\, L(p)\,\Big).
\end{equation}
Note that the indeterminacy in the definition of $T'(\n)$ is smaller than the indeterminacy in the definition of $T(\n)$, cf.
\refr{TdependsonN}: $T'(\n)$ is well defined up to a sign. If $\rank{}E$ is even, then $T'(\n)$ is a well defined complex number.
\erem

%----------------------------------------------------------------------
%----------------------------------------------------------------------
\section{Comparison Between the Refined Analytic and the Ray-Singer Torsions}\Label{S:deponalp}

Assume that $M$ is a closed oriented odd-dimensional manifold and $g^M$ is a Riemannian metric on $M$. By \reft{DetB-TRS}, there exists a
$C^0$-open neighborhood $U\subset \Flat(E)$ of the set of acyclic Hermitian connections on $E$, such that, for every $\n\in U$,
\eq{DetB-TRS3}
    \Big|\,\Detgrtet(B_{\even})\,\Big| \ = \ \TRS(\n)\cdot e^{\pi\IM\eta(\n,g^M)}.
\end{equation}
Combining this equality with \refe{DetB-eta2} and the definition of the refined analytic torsion we obtain %%
\th{T-TRS}
Assume that $M$ is a closed oriented odd-dimensional manifold and $g^M$ is a Riemannian metric on $M$. Then there exists a $C^0$-open
neighborhood $U\subset \Flat(E)$ of the set of acyclic Hermitian connections on $E$, such that $U\subset \Flat'(E,g^M)$ and for all $\n\in U$
\eq{T-TRS}
    \log\,\frac{|T(\n)|}{\TRS(\n)} \ = \ \pi\,\IM\,\eta(\n,g^M).
\end{equation}
\eth

In this section we present a local expression for the right hand side of \refe{T-TRS}.

%---------------------------
\subsection{Dependence of the $\eta$-invariant on the connection}\Label{SS:eta(alp)}
Suppose that $t\mapsto \n_t, \ t\in [0,1]$, is a smooth path of connections in $\Flat'(E,g^M)$.
We shall need the following result of Gilkey \cite[Th.~3.7]{Gilkey84}%
\footnote{Note that Gilkey considered the $\eta$-invariant of the full odd signature operator $B= B_{\even}\oplus{}B_{\odd}$. Hence,
our $\eta(\n,g^M)$ is equal to one half of the invariant considered in \cite{Gilkey84}.} %
(see also Theorem~7.6 of \cite{FarberLevine96}%
\footnote{Note, however, that Theorem~7.6 of \cite{FarberLevine96} has a wrong sign.}) %%
\th{Gilkey}
Let $\oeta(\n_t,g^M)\in \CC/\ZZ$ denote the reduction of $\eta(\n_t,g^M)$ modulo $\ZZ$. Then $\oeta(\n_t,g^M)$ depends smoothly on $t$, cf. \cite[\S1]{Gilkey84}.

1. \ If \/ $\dim M\equiv 3\ (\MOD\ 4)$ then $\oeta(\n_t,g^M)$ is independent of \/ $t\in [0,1]$.

2. \ Suppose $\dim M\equiv 1\ (\MOD\ 4)$. Set
\eq{psit}
    \psi_t \ := \ \frac{d}{dt}\,\n_t \ \in \ \Ome^1(M,\End E).
\end{equation}
Then
\eq{Gilkey}
    \frac{d}{dt}\,\oeta(\n_t,g^M) \ = \ \frac{i}{2\pi}\,\int_M\, L(p)\wedge\Tr(\psi_t),
\end{equation}
where $L(p)= L_M(p)$ is the Hirzebruch $L$-polynomial in the Pontrjagin forms of $g^M$.
\eth

%--------------------------------------------
\subsection{Cohomology class $\Arg_\n$}\Label{SS:Argalp}
Following Farber, \cite{Farber00AT}, we denote by $\Arg_\n$ the unique cohomology class $\Arg_\n\in H^1(M,\CC/\ZZ)$ such that for every closed curve
$\gam\in M$ we have
\eq{Arg}
    \det\big(\,\Mon_\n(\gam)\,\big) \ = \ \exp\big(\, 2\pi i\<\Arg_\n,[\gam]\>\,\big),
\end{equation}
where $\Mon_\n(\gam)$ denotes the monodromy of the flat connection $\n$ along the curve $\gam$ and $\<\cdot,\cdot\>$ denotes the natural pairing
\[
    H^1(M,\CC/\ZZ)\,\times\, H_1(M,\ZZ) \ \longrightarrow \ \CC/\ZZ.
\]

\rem{Arghermitian}
The notation $\Arg_\n$ is motivated by the case where $\n$ is a Hermitian connection. In this case, $\Mon_\n(\gam)$ is unitary and
$\Arg_\n\in H^1(M,\RR/\ZZ)$. Therefore, the expression $2\pi\<\Arg_\n,[\gam]\>$ is equal to the phase of the complex number
$\det(\Mon_\n(\gam))$.
\erem

\renewcommand{\d}{\partial}
%------------
\lem{Arg=psi}
Assume that $\n_t$ ($t\in [0,1]$) is a smooth path of connections. Then, using the notation introduced in \reft{Gilkey}.2, we have
\eq{Arg-psi}
    2\pi i \,\frac{d}{dt}\,\Arg_{\n_t} \ = \ -\big[\,\Tr\,\psi_t\,\big] \  \in \ H^1(M,\CC),
\end{equation}
where $\big[\,\Tr\,\psi_t\,\big]$ denotes the cohomology class of the closed differential form $\Tr\,\psi_t$.
\elem
\prf
Let $S^1$ be the standard circle and let $x\in [0,2\pi)$ be the coordinate on $S^1$. Let $\gam:S^1\to M$ be a closed curve. Fix a trivialization
of the bundle $\gam^*E\to S^1$. Let $A_t(x)$ denote the (periodically extended to $\RR$) connection form on $\gam^*E$ induced by $\n_t$. Then,
for each $t\in[0,1]$,  the monodromy $\Mon_{\n_t}(\gam)$ along $\gam$ is given by the matrix $\Phi_t(2\pi)$ where $\Phi_t(x)$ is the matrix
function solving the following initial value problem
\eq{ddxPhi}
    \begin{aligned}
    {}&\frac{\d}{\d x}\,\Phi_t(x) \ + \ A_t(x)\Phi_t(x)\ = \ 0, \qquad x\in \RR;\\
    {}&\Phi_t(0) \ = \ \Id.
    \end{aligned}
\end{equation}

Let $\dot{\Phi}_t(x)$ and $\dot{A}_t(x)$ denote the derivative with respect to $t$ of the matrices $\Phi_t(x)$ and $A_t(x)$ respectively.

We are interested in computing
\eq{ddtlogdetPhi}
    2\pi i \,\frac{d}{d t}\,\big\<\,\Arg_{\n_t},[\gam]\,\big\> \ = \
     \frac{d}{d t}\,\log\det\, \Phi_t(2\pi) \ = \ \frac{d}{d t}\,\Tr\log\,\Phi_t(2\pi) \ = \ \Tr\,\dot{\Phi}_t(2\pi)\,\Phi_t(2\pi)^{-1}.
\end{equation}
Note that, though $\log$ is a multivalued function the derivatives $\frac{\d}{\d t}\,\log\det\,\Phi_t(2\pi)$ and $\frac{\d}{\d
t}\,\Tr\log\,\Phi_t(2\pi)$ are unambiguously defined complex numbers.

From \refe{ddxPhi}, we obtain
\eq{ddxdt0}
    \frac{\d}{\d x}\,\dot{\Phi}_t(x)\ + \ \dot{A}_t(x)\,\Phi_t(x)\ + \ A_t(x)\,\dot{\Phi}_t(x) \ = 0.
\end{equation}
Hence,
\eq{ddxdt}
    \left(\,\frac{\d}{\d x}\,\dot{\Phi}_t(x)\,\right)\,\Phi_t(x)^{-1} \ = \ -\dot{A}_t(x) \ - \ A_t(x)\,\dot{\Phi}_t(x)\,\Phi_t(x)^{-1}.
\end{equation}

On the other side,
\meq{Trddx}
    \Tr\,\frac{\d}{\d x}\,\Big(\,\dot{\Phi}_t(x)\,\Phi_t(x)^{-1}\,\Big) \ = \  \Tr\,\Big(\,\frac{\d}{\d x}\,\dot{\Phi}_t(x)\,\Big)\Phi_t(x)^{-1}
    \ - \ \Tr\,\dot{\Phi}_t(x)\,\Phi_t(x)^{-1}\,\Big(\,\frac{\d}{\d x}\,\Phi_t(x)\,\Big)\,\Phi_t(x)^{-1}
    \\ = \
    \Tr\,\Big(\,\frac{\d}{\d x}\,\dot{\Phi}_t(x)\,\Big)\Phi_t(x)^{-1} \ + \ \Tr\,A_t(x)\,\dot{\Phi}_t(x)\,\Phi_t(x)^{-1},
\end{multline}
where in the last equality we used \refe{ddxPhi}.

Combining \refe{Trddx} with \refe{ddxdt} and using \refe{psit} we get
\eq{ddxPhi2}
    \frac{\d}{\d x}\,\Tr\,\dot{\Phi}_t(x)\,\Phi_t(x)^{-1} \ = \ -\,\Tr\,\dot{A}_t(x) \ = \ -\,\Tr\,\gam^*\psi_t(x).
\end{equation}
Here $\gam^*\psi_t$ denotes the pull-back of the differential form $\psi_t$ under the map $\gam:S^1\to M$.

From \refe{ddtlogdetPhi} and \refe{ddxPhi2} we obtain
\[
    \frac{d}{d t}\,2\pi i\<\Arg_{\n_t},[\gam]\>  \ = \
    \frac{d}{dt}\, \Big(\, \log\,\det\,\Phi_t(2\pi) \ - \ \log\,\det\,\Phi_t(0)\,\Big)
    \ = \ -\,\int_0^{2\pi}\,\Tr\,\dot{A}_t(x)\,dx \ = \ -\big\<\,[\Tr\,\psi_t],[\gam]\,\big\>.
\]
\eprf

From \refl{Arg=psi} and \refe{Gilkey} we obtain
\eq{detadt}
    \frac{d}{dt}\,\oeta(\n_t,g^M) \ = \  \big\<\,[L(p)]\cup \frac{i}{2\pi}[\Tr\,\psi_t],[M]\,\big\> \ = \
    \big\<\,[L(p)]\cup \frac{d}{dt}\Arg_{\n_t},[M]\,\big\>,
\end{equation}
where $\cup$ denotes the cup-product in cohomology.

Assume that $\n_t$ $(t\in [0,1])$ is a smooth path of acyclic connections and that the connection $\n_0$ is Hermitian. By \refr{Arghermitian},
$\IM\Arg_{\n_0}=0$ and, thus, \refe{Arg-psi} leads to
\eq{Arg-psi2}
    \Big[\,\int_0^t\, \Tr\,(\RE\psi_t)\,dt \,\Big]\ = \ 2\pi\,\IM\Arg_{\n_t} \ \in \ H^1(M,\RR).
\end{equation}

%-------------------------------
\subsection{Comparison with the Ray-Singer torsion}\Label{SS:RaySinger}
Let $U\subset \Flat'(E,g^M)$ be as in \reft{T-TRS}. Denote by $U'\subset U$ the set of flat connections satisfying the following condition: for
every $\n\in U'$ there exists a smooth path $t\mapsto \n_t\in U$, $t\in [0,1]$, of connections such that  $\n_0$ is Hermitian, and $\n_1=\n$.
Then $U'\subset \Flat'(E,g^M)$ is an open neighborhood of the set of acyclic Hermitian connections.

\th{RaySinger}
Let $M$ be a closed oriented odd-dimensional manifold and let $g^M$ be a Riemannian metric on $M$. Suppose $\n\in U'$. Then, with $L(p)= L_M(p)$
denoting the Hirzebruch $L$-polynomial in the Pontrjagin forms of a Riemannian metric on $M$,
\eq{RaySinger3}
    \log\,\frac{|T(\n)|}{\TRS(\n)} \ \ = \ \ \pi\,\big\<\, [L(p)]\cup \IM\Arg_{\n},[M]\,\big\>.
\end{equation}

In the case \/ $\dim M\equiv 3\ (\MOD\ 4)$
\eq{RaySinger1}
    |T(\n)| \ = \ \TRS(\n).
\end{equation}
\eth
\rem{RaySinger}
1. \ The advantage of \reft{RaySinger} over \reft{T-TRS} is that the right hand side of \refe{RaySinger3} is given by a local formula. Hence, it
might be possible to effectively compute it in some examples.

2. \ When $\dim M\equiv 3\ (\MOD\ 4)$ the right hand side of \refe{RaySinger3} vanishes since\, $L(p)$ has no component of degree $\dim{M}-1$
and, hence,  $[L(p)]\cup \IM\Arg_{\n}$\, does not have a component of degree $\dim{}M$.
\erem
\prf
Let $\n_t\in U'$  $(0\le t\le1$) be a smooth path of connections such that $\n_0$ is a Hermitian connection and $\n_1=\n$. Then, by \reft{T-TRS},
\eq{logTT=Imeta}
    \log\,\frac{|T(\n_t)|}{\TRS(\n_t)} \ = \ \pi \,\IM\,\eta(\n_t,g^M), \qquad \text{for every}\quad t\in [0,1].
\end{equation}

As the number $\oeta(\n_t,g^M)$ is defined modulo integers, its imaginary part is a well defined real number and
\[
    \IM\, \oeta(\n_t,g^M) \ = \ \IM\,\eta(\n_t,g^M).
\]
Since the connection $\n_0$ is Hermitian, $\IM\eta(\n_0,g^M)= 0$. Hence, from \reft{Gilkey}.1  we conclude that $\IM\, \eta(\n,g^M)= 0$ if $ \dim
M\equiv 3\ (\MOD\, 4)$. From \reft{Gilkey}.2 and \refe{Arg-psi2} we see that
\[
    \IM\, \eta(\n,g^M) \ = \
        \frac1{2\pi}\,\int_0^1\,\Big(\,\int_M\,L(p)\wedge \Tr(\RE\psi_t)\,\Big)\,dt \ = \
        \big\<\, [L(p)]\cup \IM\Arg_{\n},[M]\,\big\>,
\]
if $\dim M\equiv 1\ (\MOD\, 4)$. \reft{RaySinger} follows now from \refe{logTT=Imeta}.
\eprf

%---------------------------------------------------------------------
%---------------------------------------------------------------------
\section{Graded Determinant as a Holomorphic Function of a Representation of $\pi_1(M)$}\Label{S:analytic}

In this section we show, first, that the graded determinant $\Detgrtet(B_{\even}(\n,g^M))$ of the odd signature operator is, in an appropriate
sense, a holomorphic function of the connection $\n$. Then we change the point of view and consider the graded determinant as a function of the
representation of the fundamental group $\pi_1(M)$ of $M$. More precisely, each representation $\alp$ of  $\pi_1(M)$ induces a flat vector
bundle $(E_\alp,\na)$ over $M$ and we denote by $B_\alp= B(\na,g^M)$ the corresponding odd signature operator. The space $\Rep$ of all complex
$n$-dimensional representations of $\pi_1(M)$ has a natural structure of a complex algebraic variety. We show that $\Detgrtet(B_{\alp,\even})$
is a well defined holomorphic function on an open subset of this variety. Throughout the section we assume that the dimension of $M$ is odd.

\newcommand{\DiffLam}{\nolinebreak{\Diff_1(M,\Lam^\b{}T^*M\otimes{E})}}

%--------------------------------------
\subsection{Graded determinant as a holomorphic function}\Label{SS:analofconnections}
Let $M$ be a closed oriented manifold of odd dimension $d=2r-1$ and let $E\to M$ be a flat vector bundle over $M$. Every (not necessarily flat) connection on $E$ can
be viewed as a first order differential operator on $\Ome^\b(M,E)$. Thus the space $\calC(E)$ of all connections on $E$ is an affine subspace of the space $\DiffLam$
of first order differential operators on the complex vector bundle $\Lam^\b{}T^*M\otimes{}E\to M$ and, hence, inherits from $\DiffLam$ the structure of a Fr\'echet
space. See \refss{difopFr} for the definition of the Fr\'echet topology on $\DiffLam$.

Fix a Riemannian metric $g^M$ on $M$. Recall that we denote by $\Flat'(E,g^M)$  the set of flat connections $\n$ on $E$ such that the pair
$(\n,g^M)$ satisfies Assumption~I and II of \refss{assumptions}. By \refe{det-det}, the graded determinant
$\Detgrtet\big(\,B_{\even}(\n,g^M)\,\big)$ is independent of the choice of the Agmon angle $\tet\in (-\pi,0)$. Thus one obtains a function
\eq{detgrn}
    \Detgr: \operatorname{Flat}'(E,g^M) \ \longrightarrow \ \CC\backslash\{0\},
    \qquad \Detgr:\,\n\ \mapsto \ \Detgrtet\big(\,B_{\even}(\n,g^M)\,\big),
\end{equation}
where $\tet$ is any Agmon angle of $B(\n,g^M)$ in the interval $(-\pi,0)$. Recall that the notion of a holomorphic curve has been introduced in
\refss{hol}.

\prop{analyticn0}
Suppose $E$ is a vector bundle over  a closed oriented odd-dimensional Riemannian manifold $(M,g^M)$. Let $\calO\subset \CC$ be an open set and
let $\gam:\calO\to \Flat'(E,g^M)$ be a holomorphic curve in $\Flat'(E,g^M)$. Then the function $\lam\mapsto
\Detgr\big(B_{\even}(\gam(\lam),g^M)\big)$, is holomorphic on $\calO$.
\eprop

In fact, we will need a slightly more general statement. Thus we will, first, generalize \refp{analyticn0} and, then, prove this more general
version.

%---------------------
\subsection{Extension of the graded determinant to connections which are not flat}\Label{SS:extgrdet}
Recall from \reft{DetB-eta2} that, for every connection $\n\in \Flat'(E,g^M)$,
\eq{Detgr=xieta}
  \Detgrtet\big(\,B_{\even}(\n,g^{M})\,\big) \ = \ e^{\xi(\n,g^M,\tet)}\cdot e^{-i\pi\eta(\n,g^M)},
\end{equation}
where $\tet\in (-\pi/2,0)$ is any an Agmon angle for $B(\n,g^M)$ such that there are no eigenvalues of the operator $B(\n,g^M)$  in the solid
angles $L_{(-\pi/2,\tet]}$ and $L_{(\pi/2,\tet+\pi]}$.

Let $\n_0\in \Flat'(M,g^M)$. Then $B(\n_0,g^M)$ is invertible. Formula \refe{oddsign} defines the odd signature operator $B(\n,g^M)$  for an
arbitrary, not necessarily flat, connection. We wish to extend the notion of the graded determinant to operators $B(\n,g^M)$ with $\n$ in some
open neighborhood of $\Flat'(E,g^M)$ in $\calC(E)$.

The same arguments as in the proof of \refp{neigofh} show that  there exists a $C^0$-neighborhood $\calU$ of $\n_0$ in the space $\calC(E)$ of
all connections such that $B(\n,g^M)$ is invertible for all $\n\in \calU$. As in \refl{symbolofB}, the leading symbol of $B(\n,g^M)$ is
symmetric and, hence, $B(\n,g^M)$ admits an Agmon angle for $B(\n,g^M)$ such that there are no eigenvalues of the operator $B(\n,g^M)$  in the
solid angles $L_{(-\pi/2,\tet]}$, $L_{(\pi/2,\tet+\pi]}$. Thus we can use formula \refe{etainv} to define $\eta(\n,g^M)=
\eta\big(B_{\even}(\n,g^M)\big)$ for all $\n\in \calU$. Similarly, we can use the expression \refe{detB23} for $\xi$ to define
$\xi(\n,g^M,\tet)$ to all $\n\in \calU$. We now use \refe{Detgr=xieta} as the {\em definition} of $\Detgrtet\big(\,B_{\even}(\n,g^{M})\,\big)$
for $\n\in \calU$.

\prop{analyticn}
Suppose $E$ is a complex vector bundle over a closed oriented odd-dimensional Riemannian manifold $(M,g^M)$ and let\, $\calU\subset \calC(E)$ be
the $C^0$-open set defined above. Let\, $\calO\subset \CC$ be an open set and let\, $\gam:\calO\to \calU$ be a holomorphic curve in $\calC(E)$
such that there exists $\lam_0\in \calO$ with\, $\gam(\lam_0)\in \Flat'(E,g^M)$. Then the function $\lam\mapsto
\Detgrtet\big(B_{\even}(\gam(\lam),g^M)\big)$ is holomorphic in a neighborhood of $\lam_0$. Here $\tet\in (-\pi/2,0)$ is any Agmon angle for
$B_{\even}(\gam(\lam),g^M)$ such that there are no eigenvalues of the operator $B_{\even}(\gam(\lam),g^M)$  in the solid angles
$L_{(-\pi/2,\tet]}$ and $L_{(\pi/2,\tet+\pi]}$.
\eprop

\prf
Fix an Agmon angle $\tet\in (-\pi/2,0)$  for $B_{\even}(\gam(\lam_0),g^M)$ such that there are no eigenvalues of the operator $B_{\even}(\gam(\lam_0),g^M)$  in the
solid angles $L_{(-\pi/2,\tet]}$ and $L_{(\pi/2,\tet+\pi]}$. Then, for all $\lam$ in a small neighborhood of $\lam_0$, $\tet$ is also an Agmon angle for
$B_{\even}(\gam(\lam),g^M)$ and there are no eigenvalues of $B_{\even}(\gam(\lam),g^M)$  in the solid angles $L_{(-\pi/2,\tet]}$ and $L_{(\pi/2,\tet+\pi]}$.

By \refc{etahol}, the function
\[
  \calO \ \longrightarrow \ \CC, \qquad \lam \ \mapsto \ e^{-2i\pi\eta(\gam(\lam),g^M)}
\]
is holomorphic on $\calO$. Similarly, \reft{dethol}.a and the expression \refe{detB23} for $\xi$ imply that the function $\lam\mapsto
e^{2\xi(\gam(\lam),g^M,\tet)}$ also is holomorphic on $\calO$. Hence,
\[
  F(\lam) \ := \ \Detgrtet\big(\,B_{\even}(\gam(\lam),g^{M})\,\big)^2 \ = \ e^{2\xi(\gam(\lam),g^M,\tet)}\cdot e^{-2i\pi\eta(\gam(\lam),g^M)}
\]
is a non-vanishing holomorphic function on $\calO$.

Since $F(\lam)$ is a continuous function of $\lam$ and $F(\lam_0)\not=0$, we can find a neighborhood $\calO'\subset \calO$ of $\lam_0$ such that
for all $\lam\in \calO'$ we have
\[
    \big|\,F(\lam)-F(\lam_0)\,\big| \ \le \ \frac12\,\big|\,F(\lam_0)\,\big|.
\]
Then $\Detgrtet\big(\,B_{\even}(\lam,g^{M})\,\big)$ coincides on $\calO'$ with one of the two analytic square roots of $F(\lam)$.
\eprf
\rem{xi-eta}
The above arguments show a very close relationship between $e^{\xi(\n,g^M,\tet)}$ and $e^{-i\pi\eta(\n,g^M)}$. Each of these numbers by itself depends on the choice
of the Agmon angle $\tet$. But their product is a well defined holomorphic function. This relationship plays a very important role in the whole paper since it
explains many features of the refined analytic torsion.
\erem

%----------------------------------------
\subsection{Space of representations of the fundamental group}\Label{SS:Rep}
Let $M$  be a closed oriented manifold of odd dimension $d=2r-1$, where $r\ge1$. Denote by $\tilM$ the universal cover of $M$ and by $\pi_1(M)$ the fundamental group
of $M$, viewed as the group of deck transformations of $\tilM\to M$. The set $\Rep$ of all $n$-dimensional complex representations of $\p$ has a natural structure of
a complex algebraic variety. Indeed, $\pi_1(M)$ is a finitely presented group, i.e., it is generated by a finite number of elements $\gam_1\nek \gam_L$, which satisfy
finitely many relations. Hence, a representation $\alp\in \Rep$ is given by $2L$ invertible $n\times{n}$-matrices $\alp(\gam_1)\nek \alp(\gam_L),\
\alp(\gam_1^{-1})\nek \alp(\gam_L^{-1})$ with complex coefficients satisfying finitely many polynomial equations. In other words, a representation $\alp$ is given by
a point of the direct product ${\operatorname{Mat}}_{n\times{n}}(\CC)^{2L}$ of $2L$ copies of the space ${\operatorname{Mat}}_{n\times{n}}$ of $n\times{n}$-matrices
with complex coefficients.

In the sequel, we fix generators $\gam_1\nek \gam_L$ of $\pi_1(M)$ and view $\Rep$ as an algebraic subset of
${\operatorname{Mat}}_{n\times{n}}(\CC)^{2L}$ with the induced topology. For $\alp\in \Rep$, we denote by
\[
  E_\alp \ := \ \tilM\times_\alp \CC^n \ \longrightarrow M
\]
the flat vector bundle induced by $\alp$. Let $\n_\alp$ be the flat connection on $E_\alp$ induced from the trivial connection on
$\tilM\times\CC^n$. We will also denote by $\n_\alp$ the induced differential
\[
    \n_\alp:\, \Ome^\b(M,E_\alp) \ \longrightarrow \Ome^{\b+1}(M,E_\alp),
\]
where $\Ome^\b(M,E_\alp)$ denotes the space of smooth differential forms of $M$ with values in $E_\alp$.

For each connected component%
\footnote{In this paper we always consider the classical (not the Zariski) topology on the complex analytic space $\Rep$.} %
$\calC\subset \Rep$ of the space of representations all the bundles $E_\alp$ are isomorphic, see e.g.  \cite{GoldmanMillson88}.

Let $\Repo\subset \Rep$ denote the (possibly empty) set of all representations $\alp\in \Rep$ such that the connection $\na$ is acyclic. A
representation $\alp\in \Rep$ is called {\em unitary} if there exists a Hermitian scalar product $(\cdot,\cdot)$ on $\CC^n$ which is preserved
by the matrices $\alp(\gam)$ for all $\gam\in \p$. The scalar product $(\cdot,\cdot)$ induces a flat Hermitian metric $h^{E_\alp}$ on the bundle
$E_{\alp}$. We denote the set of unitary representations by $\Reph$. One might think of
\[
    \Reph\ \subset\ \Rep
\]
as {\em the real locus} of the complex algebraic variety $\Rep$. Set
\[
    \Repho \ :=\  \Reph\cap\Repo.
\]

%---------------------------------------------------------
\subsection{Graded determinant of the odd signature operator as a function on the space of representations}\Label{SS:DetRep}
Fix a Riemannian metric $g^M$ on $M$.  Let
\[
    B_\alp \ := \ B(\na,g^M):\, \Ome^\b(M,E_\alp)\ \longrightarrow \ \Ome^\b(M,E_\alp)
\]
and let $B_{\alp,\even}$ denote the restriction of $B_\alp$ to $\Ome^{\even}(M,E_\alp)$.

Suppose that for some representation $\alp_0\in \Repo$ the operator $B_{\alp_0}$ is invertible (in other words, we assume that
$(\n_{\alp_0},g^M)$ satisfies Assumption~I and II of \refss{assumptions}). Then there exists an open neighborhood (in classical topology)
$V\subset \Rep$ of the set of acyclic unitary representations such that, for all $\alp\in V$ the pair $(\na,g^M)$ satisfies Assumption~I and II
of \refss{assumptions}. Thus, for all $\alp\in V$, the graded determinant $\Detgrtet(B_{\alp,\even})$ is defined, where $\tet\in (-\pi,0)$ is an
Agmon angle for $B_\alp$.

%--------
%------------
\th{analytic}
Let $M$ be a closed oriented odd-dimensional manifold and let $g^M$ be a Riemannian metric on $M$. Let $\calO\subset \CC$ be a connected open
set and let\, $\gam:\calO\to \Repo$ be a holomorphic curve. Assume that for  $\lam_0\in \calO$ the connection $\n_{\gam(\lam_0)}$ on\,
$E_{\gam(\lam_0)}\to M$\, satisfies Assumption~II of \refss{assumptions} (with respect to the given metric $g^M$). Then the function
\eq{weakanalytic}
   \lam \ \mapsto \Detgrtet\big(\,B_{\gam(\lam),\even}\,\big)
\end{equation}
is holomorphic in a neighborhood of $\lam_0$.
\eth
\newcommand{\Mat}{\operatorname{Mat}}

\prf
First, we need to introduce some additional notations. Let $E$ be a vector bundle over $M$ and let $\n$ be a (not necessarily flat) connection
on $E$. Fix a base point $x_*\in M$ and let $E_{x_*}$ denote the fiber of $E$ over $x_*$. We will identify $E_{x_*}$ with $\CC^n$ and
$\pi_1(M,x_*)$ with $\pi_1(M)$.

For a closed  path $\phi:[0,1]\to M$ with $\phi(0)=\phi(1)= x_*$, we denote by $\Mon_\n(\phi)\in \End{}E_{x_*}\simeq \Mat_{n\times n}(\CC)$ the
monodromy of $\n$ along $\phi$, cf. \refe{ddxPhi}. Note that, if $\n$ is flat then $\Mon_\n(\phi)$ depends only on the class $[\phi]$ of $\phi$
in $\pi_1(M)$. Hence, if $\n$ is flat, then the map $\phi\mapsto \Mon_\n(\phi)$ defines an element of $\Rep$, called  the {\em monodromy
representation} of $\n$.

Suppose now that $\calO\subset \CC$ is a connected open set. Let $\gam:\calO\to \Repo$ be a holomorphic curve. By Proposition~4.5 of
\cite{GoldmanMillson88}, all the bundles $E_{\gam(\lam)}, \ \lam\in \calO$, are isomorphic to each other. In other words, there exists a vector
bundle $E\to M$ and a family of flat connections $\n_\lam,\ \lam\in \calO$, on $E$, such that the monodromy representation of $\n_\lam$ is
isomorphic to $\gam(\lam)$ for all $\lam\in \calO$. Moreover, the family $\n_\lam$ can be chosen to be real differentiable, i.e., such that for
every $\lam\in \calO$ there exist $\ome_1,\ome_2\in \Ome^1(M,\End E)$ with
\eq{nmu-nlam}
    \n_\mu \ = \ \n_\lam \ + \ \RE(\mu-\lam)\cdot\ome_1 \ + \ \IM(\mu-\lam)\cdot\ome_2 \ + \ o(\mu-\lam),
\end{equation}
where $o(\mu-\lam)$ is understood in the sense of the Fr\'echet topology on $\calC(E)$ introduced in \refss{analofconnections}.

By \refl{realanal}  there exist a smooth form $\ome\in \Ome^1(M,\End E)$  with $\n_{\lam} \ome = 0$ and a family  $G(\mu)\in \End{}E$ ($\mu\in
\calO$) of gauge transformations such that $G(\lam)=\Id$ and
\eq{guage}\notag
   \n_\lam \ + \ (\mu-\lam)\,\ome \ = \ G(\mu)\cdot\n_\mu\cdot G(\mu)^{-1} \ + \ o(\mu-\lam).
\end{equation}
Note that the connection $\n_\lam + (\mu-\lam)\,\ome$ is not necessarily flat.

From the definition of the the odd signature operator it then follows that
\eq{Bguage}
    B\,\big(\n_\lam+(\mu-\lam)\,\ome,g^M\big) \ = \ G(\mu)\cdot B(\n_\mu,g^M)\cdot G(\mu)^{-1} \ + \ o(\mu-\lam),
\end{equation}
where $o(\mu-\lam)$ is understood in the sense of the Fr\'echet topology introduced in \refss{difopFr}.

Suppose now that $\lam$ is close enough to $\lam_0$ so that the connection $\n_{\gam(\lam)}$ satisfies Assumption~II of \refss{assumptions}
(with respect to the metric $g^M$). Recall that in \refss{extgrdet} we extended the definition of the graded determinant of $B_{\even}(\n,g^M)$
to the case when the connection $\n$ is not necessarily flat. Thus $$\Detgrtet\big(\,B_{\even}(\n_\lam+(\mu-\lam)\,\ome,g^M)\,\big)$$ is defined
for all $\mu\in \CC$ close enough to $\lam$.

By \refp{analyticn}, the map
\[
     \mu\ \mapsto \ \Detgrtet\big(B_{\even}(\n_\lam+(\mu-\lam)\ome,g^M)\big)
\]
is holomorphic near $\lam$. Hence, there exists a number $a\in \CC$ such that
\eq{complder}
    \Detgrtet\big(\,B_{\even}(\n_\lam+(\mu-\lam)\,\ome,g^M)\,\big) \ = \ \Detgrtet\big(\,B_{\even}(\n_\lam,g^M)\,\big) \ +\ a\cdot(\mu-\lam)\ + \
    o(\mu-\lam).
\end{equation}
On the other side, \refe{Bguage} implies that
\eq{guagenmu}
\begin{aligned}
    \Detgrtet\big(\,B_{\even}(\n_\lam+(\mu-\lam)\,\ome,g^M)\,\big)
    \ &= \ \Detgrtet\big(\,G(\mu)\cdot B_{\even}(\n_\mu,g^M)\cdot G(\mu)^{-1}\,\big) \ + \ o(\mu-\lam)
    \\ &= \  \Detgrtet\big(\, B_{\even}(\n_\mu,g^M)\,\big) \ + \ o(\mu-\lam).
\end{aligned}
\end{equation}
Combining \refe{complder} with \refe{guagenmu} we obtain
\[
    \Detgrtet\big(\, B_{\even}(\n_\mu,g^M)\,\big) \ = \ \Detgrtet\big(\,B_{\even}(\n_\lam,g^M)\,\big) \ +\ a\cdot(\mu-\lam)\ + \
    o(\mu-\lam).
\]
Since the above equality holds for all $\lam$ close enough to $\lam_0$ the theorem is proven.
\eprf

%------------
\cor{analytic}
Let $M$ be a closed oriented odd-dimensional manifold. Let $V\subset \Rep$ be an open set such that for every $\alp\in V$ there exists a
Riemannian metric $g^M$ such that the connection $\na\in \Flat'(E_\alp,g^M)$ (cf. \refss{openset1}). Assume, further, that all the points of $V$
are regular points of the complex algebraic set $\Rep$. Then the map
\eq{D(alp)}
    \Det:\, V \ \longrightarrow \ \CC, \qquad \Det:\,\alp \  \mapsto\ \Det(\alp) \ := \ \Detgrtet(B_{\alp,\even}).
\end{equation}
is holomorphic. Here $\tet\in (-\pi,0)$ is an Agmon angle for $B_{\alp,\even}$.
\ecor

\prf
By Hartogs' theorem (cf., for example, \cite[Th.~2.2.8]{HormanderSCV}), a function on a smooth algebraic variety is holomorphic if its
restriction to each holomorphic curve is holomorphic. Hence, the corollary follows immediately from \reft{analytic}.
\eprf

\rem{analytic2}
In \refs{comb} below we mostly view the graded determinant of the odd signature operator as a function on the space of representations rather
than as a function on the space of flat connections. As $\Rep$ is a {\em finite dimensional}  algebraic variety, we can use the methods of
complex analysis of holomorphic functions on finite dimensional varieties.
\erem

By the definition of the refined analytic torsion, cf. Definitions~\ref{D:refantor} and \ref{D:refantor2}, \reft{analytic} and \refc{analytic}
imply now the following
\cor{Tanalytic}
Let $M$ be a closed oriented odd-dimensional manifold.

1. \ Let $\calO\subset \CC$ be an open set and let $\gam:\calO\to \Repo$ be a holomorphic curve. Assume that for $\lam_0\in \calO$ there exists
a Riemannian metric $g^M$ so that the connection $\n_{\gam(\lam_0)}$ satisfies Assumption~II of \refss{assumptions}. Then the function
\eq{weakanTalytic}
   \lam \ \mapsto T(\n_{\gam(\lam)}\,\big)
\end{equation}
is holomorphic in a neighborhood of $\lam_0$.

2. \ Let $V\subset \Repo$ denote the open set of all representations $\alp$ such that for some Riemannian metric $g^M$ on $M$ the connection $\na\in
\Flat'(E_\alp,g^M)$. Let $\Sig\subset \Rep$ denote the set of singular points of the complex algebraic variety $\Rep$. Then $\alp\mapsto T(\na)$ is a
holomorphic function on $V\backslash\Sig$.
\ecor

%----------------------------------------------------------------------
%----------------------------------------------------------------------
\section{Comparison with Turaev's Refinement of the Combinatorial Torsion}\Label{S:comb}

In \cite{Turaev86,Turaev90}, Turaev introduced a refinement $\TTur_\alp(\eps,\gro)$ of the combinatorial torsion associated to the representation
$\alp$ of $\pi_1(M)$. This refinement depends on an additional combinatorial data, denoted by $\eps$ and called the {\em Euler structure} as well as
on the {\em cohomological orientation} of $M$, i.e. on the orientation $\gro$ of the determinant line of the cohomology $H^\b(M,\RR)$ of $M$. There
are two versions of the Turaev torsion -- the homological and the cohomological one. In this paper it is more convenient for us to use the
cohomological Turaev torsion as it is defined in Section~9.2 of \cite{FarberTuraev00}. For $\alp\in \Repo$ the cohomological Turaev torsion
$\TTur_\alp(\eps,\gro)$
is a non-vanishing complex number. If, moreover, $\alp\in \Repho$ the absolute value of the Turaev torsion is equal to the Reidemeister torsion.%
\footnote{Here the Reidemeister torsion is understood as the positive real number defined, for example, in Definition~1.1 of \cite{RaSi1}.} One
can view  \reft{DetB-TRS} as an analytic analogue of this result, where the role of the Reidemeister torsion is played by the Ray-Singer
torsion. Another property of the Turaev torsion is that it is a holomorphic function of $\alp\in \Repo$. In \refc{Tanalytic} we established the
same property for the refined analytic torsion.

Though, in general, the refined analytic torsion $T_\alp= T(\na)$ and the Turaev torsion $\TTur_\alp(\eps,\gro)$ are not equal they are very closely
related. In this section we establish this relationship. As an application we strengthen and generalize a theorem of Farber \cite{Farber00AT} about
the relationship between the Turaev torsion and the $\eta$-invariant.

%----------------------------------------
\subsection{Notation}\Label{SS:notationlast}
Let $M$ be a closed oriented odd-dimensional manifold. In this section we view the refined analytic torsion as a function of a representation of
$\pi_1(M)$. Let $V\subset \Repo$ be an open set consisting of representations $\alp$ such that $\na\in \Flat'(E,g^M)$. Let $V'\subset V$ be the open
subset of $V$ such that, for all $\alp\in V'$, the connection $\na$ belongs to the open set $U'$ defined in \refss{RaySinger}.  For every $\alp\in V$
we set $T_\alp:= T(\na)$, $\TRS_\alp:= \TRS(\na)$, $\eta_\alp:= \eta(\na,g^M)$, etc.
%------------------------------------------
\subsection{Comparison between the Turaev and the Ray-Singer Torsion}\Label{SS:Turaevtorsion}
Theorem~10.2 of \cite{FarberTuraev00} establishes a relationship between the Turaev and the Ray-Singer torsions for real representations $\alp$. The
following result is an immediate extension of this theorem to complex representations.
\th{TTur-TRS}
Suppose $M$ is a closed oriented odd-dimensional manifold. Let $c(\eps)\in H_1(M,\ZZ)$ denote the characteristic class of the Euler structure $\eps$,
cf. \cite{Turaev90} or Section~5.2 of \cite{FarberTuraev00}. Then, for every $\alp\in \Repo$,
\eq{Tcomb-TRS}
    \log\,\frac{|\TTur_\alp(\eps,\gro)|}{\TRS_\alp} \ = \ -\,\pi\,\<\IM\Arg_\alp,c(\eps)\>,
\end{equation}
where the cohomology class\, $\Arg_{\alp}:= \Arg_{\na}\in H^1(M,\CC/\ZZ)$ is defined in \refss{Argalp} and $\<\cdot,\cdot\>$ denotes the natural
pairing
\[
    H^1(M,\CC/\ZZ)\times H_1(M,\ZZ) \ \longrightarrow \ \CC/\ZZ.
\]
In particular, if $\alp\in \Repho$ then
\eq{Tcomb-TRSun}
    \big|\,\TTur_\alp(\eps,\gro)\,\big| \ = \ \TRS_{\alp}.
\end{equation}
\eth
Note that, though $\<\Arg_\alp,c(\eps)\>$ is defined only modulo $\ZZ$, its imaginary part $\IM\<\Arg_\alp,c(\eps)\>$ is a well defined complex
number.
\prf
Let $\alp^{\RR}$ denote the representation $\alp$ considered as a real representation. Then, for every closed curve $\gam$ in $M$, we have
\[
    \det\Mon_{\alp^\RR}(\gam) \ = \ \big|\det\Mon_\alp(\gam)\,\big|^2.
\]
Define $\Arg_{\alp^\RR}\in H^1(M,\RR/\ZZ)$ by
\[
    \det\,\Mon_{\alp^\RR}(\gam) \ = \ \exp\,\big(\,2\pi i\<\Arg_{\alp^\RR},[\gam]\>\big).
\]
Then, from \refe{Arg}, we obtain
\[
    \exp\,\big(\,2\pi i\<\Arg_{\alp^\RR},[\gam]\>\big) \ = \ \exp\,\big(\,2\pi i\<2i\IM\Arg_{\alp},[\gam]\>\big).
\]
Hence,
\eq{Arg-ArgR}
    \<\Arg_{\alp^\RR},[\gam]\> \ \equiv \ 2i\,\IM\,\<\Arg_\alp,[\gam]\> \quad \MOD \ \ZZ.
\end{equation}

Let $\TTur_{\alp^\RR}(\eps,\gro)$ and $\TRS_{\alp^\RR}$ denote the Turaev and the Ray-Singer torsions associated to the representation $\alp^\RR$. Then
\eq{TC-TR}
    |\TTur_{\alp}(\eps,\gro)|^2 \ = \ \TTur_{\alp^\RR}(\eps,\gro), \qquad (\TRS_\alp)^2 \ = \ \TRS_{\alp^\RR}.
\end{equation}

By formula (10.3) of \cite{FarberTuraev00}, we have
\[
    \left(\,\frac{\TTur_{\alp^\RR}(\eps,\gro)}{\TRS_{\alp^\RR}}\,\right)^2 \ = \ \left|\, \exp\big(2\pi
    i\<\Arg_{\alp^\RR},c(\eps)\>\,\big)\,\right|.
\]
Combining this equality with \refe{Arg-ArgR} and \refe{TC-TR}, we obtain \refe{Tcomb-TRS}.

If $\alp$ is unitary, then $\IM\Arg_{\alp}=0$ and \refe{Tcomb-TRSun} follows.
\eprf

%-------------------------------------------
\subsection{Comparison between the Turaev and the refined analytic torsions}\Label{SS:T-TTur}
To simplify the notation let us denote by $\hatL(p)\in H_\b(M,\ZZ)$ the Poincar\'e dual of the cohomology class $[L(p)]$. Let $\hatL_1\in H_1(M,\ZZ)$
denote the component of $\hatL(p)$ in $H_1(M,\ZZ)$. Then
\[
    \<\, [L(p)]\cup \Arg_\alp, [M]\,\> \ = \ \<\Arg_\alp,\hatL_1\> \ \in \ \CC/\ZZ.
\]

Recall that the neighborhood $V'$ of $\Repho$ was defined in \refss{notationlast}. If $\alp\in V'$ then by Theorems~\ref{T:RaySinger} and
\refe{Tcomb-TRS}
\eq{RelogTTTur}
    \RE\,\log\frac{T_\alp}{\TTur_\alp(\eps,\gro)} \ = \ \log\frac{|T_\alp|}{|\TTur_\alp(\eps,\gro)|}
    \ = \
    \pi\,  \big\<\IM\Arg_\alp,c(\eps)+\hatL_1\big\>.
\end{equation}

Let $\Sig$ denote the set of singular points of the complex analytic set $\Rep$. By \refc{Tanalytic}, the refined analytic torsion $T_\alp$ is a
non-vanishing holomorphic function of $\alp\in V\backslash\Sig$. By the very construction \cite{Turaev86,Turaev90,FarberTuraev00} the Turaev
torsion is a non-vanishing holomorphic function of $\alp\in \Repo$. Hence,
\[
    \frac{T_\alp}{\TTur_\alp(\eps,\gro)}
\]
is a holomorphic function on $V'\backslash\Sig$.

By construction of the cohomology class $\Arg_\alp$, for every homology class $z\in H_1(M,\ZZ)$, the expression
 \(
   e^{2\pi\,i\,\<\Arg_\alp,z\>}
 \)
is a holomorphic function on $\Rep$.

Now the expression \refe{RelogTTTur} can be rewritten as
\[
    \left|\,\frac{T_\alp^2}{\TTur_\alp(\eps,\gro)^2}\,\right| \ = \ \left|\, e^{-2\pi i\<\Arg_\alp,c(\eps)+\hatL_1\>}\,\right|, \qquad \alp\in V'.
\]
If the absolute values of two non-vanishing holomorphic functions are equal on a connected open set then the functions must be equal up to a
factor $\phi\in \CC$ with $|\phi|=1$. Hence, on each connected component $\ccomp\subset V'\backslash\Sig$, there exists a constant
$\phi_\ccomp(\eps,\gro)\in \RR$, depending on $\eps$ and $\gro$, so that
\eq{Talp-Tturalp}
    \left(\, \frac{T_\alp}{\TTur_\alp(\eps,\gro)}\,e^{-i\phi_\ccomp(\eps,\gro)}\,\right)^2
    \ = \ e^{-2\pi i\<\Arg_\alp,c(\eps)+\hatL_1\>}, \qquad \alp\in \ccomp.
\end{equation}
Note that the constants $\phi_\ccomp(\eps,\gro)$ are defined up to an additive multiple of $\pi$. Since $T_\alp$, $\TTur_\alp(\eps,\gro)$, and
$e^{-2\pi i\<\Arg_\alp,c(\eps)+\hatL_1\>}$ depend continuously on $\alp\in V'$, we can choose these constants so that
\[
    \phi_{\ccomp_1}(\eps,\gro)\ = \ \phi_{\ccomp_2}(\eps,\gro)
\]
whenever $\ccomp_1$ and $\ccomp_2$ are contained in the same connected component of $V'$. Thus \refe{Talp-Tturalp} remains valid if $\ccomp$ is
a connected component of $V'$.

Let now $\ccomp$ be a connected component of\, $V'$. \refe{Talp-Tturalp} implies that on $\ccomp$ the function $\alp\mapsto e^{-2\pi
i\<\Arg_\alp,c(\eps)+\hatL_1\>}$ admits a square root
\eq{sqrtArg}
    {}_\ccomp\sqrt{e^{-2\pi i\<\Arg_\alp,c(\eps)+\hatL_1\>}} \ = \ \frac{T_\alp}{\TTur_\alp(\eps,\gro)}\,e^{-i\phi_\ccomp(\eps,\gro)},
\end{equation}
which is a holomorphic function on each connected component of $\ccomp\backslash\Sig$.

Thus we have proven the following extension of the Cheeger-M\"uller theorem about the equality between the Reidemeister and the Ray-Singer torsions.
\th{TTTur}
Suppose $M$ is a closed oriented odd dimensional manifold. Let $\eps$ be an Euler structure on $M$ and let $\gro$ be a cohomological orientation of
$M$. Let $V'\subset \Repo$ be as in \refss{notationlast}. Then, for each connected component $\ccomp$ of \/ $V'$, there exists a constant
$\phi_\ccomp= \phi_\ccomp(\eps,\gro)\in \RR$, depending on $\eps$ and $\gro$, such that
\eq{logTTTur}
    \frac{T_\alp}{\TTur_\alp(\eps,\gro)} \ = \ e^{i\phi_\ccomp}\,{}_{{}_\ccomp}\sqrt{e^{-2\pi i\<\Arg_\alp,c(\eps)+\hatL_1\>}},
\end{equation}
where ${{}_{{}_\ccomp}}\sqrt{\ }$ is the analytic square root defined in \refe{sqrtArg}.
\eth
\rem{Phi(alp)}
1.\ In general, the number
 \(
    \pi i\,\<\Arg_\alp,c(\eps)+\hatL_1\>
 \)
is defined only modulo $\pi{i}$. However, in some interesting examples it is defined modulo $2\pi{i}$. Suppose that the homology class
$c(\eps)+\hatL_1\in H_1(M,\ZZ)$ is divisible by $2$, i.e., there exists an integer homology class $z\in H_1(M,\ZZ)$ such that
\eq{4z}
    2z \ = \  \hatL_1 \ + \ c(\eps).
\end{equation}
Then
\[
     \pi i\,\<\Arg_\alp,c(\eps)+\hatL_1\> \ = \ 2\pi i\, \<\Arg_\alp,z\>
\]
is defined modulo $2\pi{i}$. In this situation we can rewrite \refe{logTTTur} as
\eq{TTTur2}
    \frac{T_\alp}{\TTur_\alp(\eps,\gro)} \ = \  e^{i\phi_\ccomp}\cdot e^{-i\pi \<\Arg_\alp,c(\eps)+\hatL_1\>}.
\end{equation}

2. \ It would be very interesting to calculate the constant $\phi_\ccomp(\eps,\gro)$. In particular, it would be interesting to know whether it
actually depends on the connected component $\ccomp$ of \/ $V'$. Another interesting question is for which representations $\alp$ one can find an
Euler structure $\eps$ and the cohomological orientation $\gro$ such that $T_\alp= \TTur_\alp(\eps,\gro)$.
\erem

The following simplest case of \reft{TTTur} is still very interesting.
\cor{TTTur}
Under the assumptions of \reft{TTTur} suppose that $c(\eps)=0$.
\begin{enumerate}
\item
If $\dim{M}\equiv 3$ ($\MOD\,4$), then the ratio
 \(
     {T_\alp}/{\TTur_\alp(\eps,\gro)}
 \)
is locally constant on $V'$.
\item
If $\dim{M}\equiv 1$ ($\MOD\,4$) and $\hatL_1\in H_1(M,\ZZ)$ is divisible by 2 (i.e., there exists a homology class $z\in H_1(M,\ZZ)$ such that
$2z=\hatL_1$), then
\eq{TTTurcor}\notag
    \frac{T_\alp}{\TTur_\alp(\eps,\gro)} \cdot e^{-i\pi\,\<[L(p)]\cup \Arg_{\alp},[M]\>}
\end{equation}
is locally constant on $V'$.
\end{enumerate}
\ecor
\rem{TTTurcor}
The condition $c(\eps)=0$ holds for the Euler structure $\eps$ which is used in the definition of the  {\em absolute torsion} introduced by Farber and Turaev
\cite{FarberTuraev99,FarberTuraev00}. Further, $\hatL_1$ is divisible by 2 if, for example, the Stiefel-Whitney class $w_{d-1}(M)\in H^{d-1}(M,\ZZ_2)$ vanishes and
$H_1(M)$ has no 2-torsion, cf. \S3.1 of \cite{Farber00AT}.
\erem

%-------------------------
\subsection{Phase of the Turaev torsion of a unitary representation}\Label{SS:argTTur}
As an application of our study of the refined analytic torsion we obtain a result about the phase of the Turaev torsion which improves and
generalizes a theorem of Farber \cite{Farber00AT}, cf. \refr{Farber} below.

We denote the phase of a complex number $z$ by $\Ph(z)\in [0,2\pi)$ so that $z= |z|e^{i\Ph(z)}$.

Suppose  $\alp\in \Repho$ is a unitary representation. Then the number $\xi_\alp= \xi(\na,g^M,\tet)$, defined in \refe{calB}, is real (in fact,
in this case, $\xi_\alp$ coincides with  $\log\TRS_\alp$,\/ cf. \refe{RaySingertor3}). Moreover, the $\eta$-invariant $\eta_\alp$ is real, cf.
\refss{symmetric}. Thus, \refe{DetB-eta2} and the definition of the refined analytic torsion (\refd{refantor2}) imply
\eq{PhT}
    \Ph(T_\alp) \ = \ -\,\pi\,\eta_\alp \ + \ \pi\,\frac{\rank \alp}2\, \int_{N}\,L(p)
    \qquad \MOD\, 2\pi\,\ZZ,
\end{equation}
where $N$ is an oriented manifold whose oriented boundary is the disjoint union of two copies of $M$. The second term on the right hand side of
\refe{PhT} vanishes if $\dim M\equiv 1\ (\MOD 4)$.

Combining \refe{PhT} with \reft{TTTur} we obtain the following
\th{extofFarber}
Under the assumptions of \reft{TTTur} suppose that $\alp_1,\, \alp_2\in \Repho$ are unitary representations which lie in the same connected
component of\, $V'$. In particular, they have the same rank. Then, modulo $\pi\,\ZZ$,
\meq{extofFarber2}
    \Ph(\TTur_{\alp_1}(\eps,\gro))  + \pi\, \eta_{\alp_1}  -
    \pi\,\big\<\Arg_{\alp_1},c(\eps)+\hatL_1\big\>
    \\ \equiv \ \Ph(\TTur_{\alp_2}(\eps,\gro))  + \pi\, \eta_{\alp_2} -
    \pi\, \big\<\Arg_{\alp_2},c(\eps)+\hatL_1\big\>.
\end{multline}
If, moreover, the condition \refe{4z} is satisfied, then \refe{extofFarber2} holds modulo $2\pi\ZZ$.
\eth

%-------------------------
\subsection{Sign of the absolute torsion}\Label{SS:signofAT}
Suppose that the Stiefel-Whitney class
\[
    w_{d-1}(M)\ \in \ H^{d-1}(M,\ZZ_2)
\]
vanishes (this is always the case when $\dim{M}\equiv 3\ (\MOD\,4)$, cf. \cite{Massey60}). Then one can choose an Euler structure $\eps$ such that
$c(\eps)=0$, cf. \cite[\S3.2]{FarberTuraev99}. Assume, in addition, that the first Stiefel-Whitney class $w_1(E_\alp)$, viewed as a homomorphism
$H_1(M,\ZZ)\to \ZZ_2$, vanishes on the 2-torsion subgroup of $H_1(M,\ZZ)$. In this case there is also a canonical choice of the cohomological
orientation $\gro$, cf.  \cite[\S3.3]{FarberTuraev99}. Then the Turaev torsion $\TTur_\alp(\eps,\gro)$ corresponding to any  $\eps$ with $c(\eps)=0$
and the canonically chosen $\gro$ will be the same.

If the above assumptions on $w_{d-1}(M)$ and $w_1(E_\alp)$ are satisfied, then the number
\[
    \Tabs_\alp \ := \ \TTur_\alp(\eps,\gro)\ \in\ \CC, \qquad (\,c(\eps)\,=\,0\,),
\]
is canonically defined, i.e., is independent of any choices. It was introduced by Farber and Turaev, \cite{FarberTuraev99}, who called it the {\em
absolute torsion}. If $\alp\in \Repho$, then $\Tabs_\alp\in \RR$, cf. Theorem~3.8 of \cite{FarberTuraev99} and, hence,
\[
    e^{i\Ph(\Tabs_\alp)} \ = \ \sign(\Tabs_\alp).
\]

From \refc{TTTur}, \refr{TTTurcor}, and \reft{extofFarber} we now obtain the following
\th{Farber}
Under the assumptions of \reft{TTTur} suppose that $\alp_1,\, \alp_2\in \Repho$ are unitary representations which lie in the same connected component
of\, $V'$.
\begin{enumerate}
\item
Let $\dim{M}\equiv 3$ ($\MOD\,4$). Assume that the first  Stiefel-Whitney class $w_1(E_{\alp_1})= w_1(E_{\alp_2})$ vanishes on the 2-torsion subgroup
of $H_1(M,\ZZ)$. Then
\[
    \sign\big(\,\Tabs_{\alp_1})\cdot e^{i\pi\eta_{\alp_1}} \ = \ \sign\big(\,\Tabs_{\alp_2})\cdot e^{i\pi\eta_{\alp_2}}.
\]
\item
Let\/ $\dim{M}\equiv 1$ ($\MOD\,4$). Assume that $w_{d-1}(M)=0$ and the group $H_1(M,\ZZ)$ has no 2-torsion. Then
\[
    \sign\big(\,\Tabs_{\alp_1})\cdot e^{i\pi\big(\, \eta_{\alp_1}- \< [L(p)]\cup \Arg_{\alp_1},[M]\>\,\big)}
    \ = \ \sign\big(\,\Tabs_{\alp_2})\cdot e^{i\pi\big(\, \eta_{\alp_2}- \< [L(p)]\cup \Arg_{\alp_2},[M]\>\,\big)}.
\]
\end{enumerate}
\eth

\rem{Farber}
For the special case when there is a real analytic path $\alp_t$ of unitary representations connecting $\alp_1$ and $\alp_2$ such that the
twisted deRham complex \refe{deRham} is acyclic for all but finitely many values of $t$, \reft{Farber} was established by Farber, using a completely different method,%
\footnote{Note that Farber's definition of the $\eta$-invariant differs from ours by a factor of 2 and also that the sign in front of
$\< [L(p)]\cup \Arg_{\alp_1},[M]\>$ is wrong in \cite{Farber00AT}.} see \cite{Farber00AT}, Theorems~2.1 and 3.1.
\erem

%-----------------------------------------------------------
%-----------------------------------------------------------
\appendix
\section{Determinant of an Operator with the Spectrum Symmetric about the Real Axis}\Label{S:det-eta-sa}

In this appendix we show that for a wide and important class of differential operators, including the self-adjoint ones, formula \refe{det-eta}
represents $\LD_{\tet}(D)$ as a sum of its real and imaginary parts.

%----------
\defe{realcoef}
The spectrum of $D$ is {\em symmetric  with respect to the real axis} if the following condition holds: if $\lam$ is an eigenvalue of $D$, then
$\olam$ also is an eigenvalue of\/ $D$ and has the same algebraic multiplicity as $\lam$.
\edefe
Note that every operator with {\em real coefficients} has this property. See \cite{BrAbanov} for examples of other
interesting operators with symmetric spectrum.%
\footnote{All the operators considered in \cite{BrAbanov} have spectrum symmetric about the {\em imaginary axis}. However,  the spectrum of the
operator considered in Section~5 of \cite{BrAbanov} is also symmetric about the real axis. Further, the spectrum of the operator
$\Gam_L\cdot{}D_{mn}$, discussed at the end of Section~6.8 of \cite{BrAbanov}, is symmetric about the real axes.}

%-----------------------
\th{realcoef}
Let $D:C^\infty(M,E)\to C^\infty(M,E)$ be an injective elliptic differential operator of order $m$ with self-adjoint leading symbol, whose spectrum
is symmetric about the real axis. Let $\tet\in (-\pi/2,0)$ be an Agmon angle for $D$. Then the numbers $\zet_{2\tet}(0,D^2)$, $\eta(D)$, and
$\Det_{2\tet}(D^2)= e^{-\zet_{2\tet}'(0,D^2)}$ are real. In particular, the following analogue of \refe{det-eta-sa} and \refe{det-eta-sa2} holds:
\eq{realcoef}
    \Det_\tet(D) \ = \ (-1)^{m_-}\cdot\sqrt{\big|\Det_{2\tet}(D^2)\big|}\,\cdot e^{-i\pi\big(\eta(D)-\frac12\zet_{2\tet(0,D^2)}\big)},
\end{equation}
where $m_-=\rank P_-$ is the number of the eigenvalues of $D$ (counted with their algebraic multiplicities) on the negative part of the imaginary
axis, cf. \refss{etainv}.
\eth
%%%
\cor{detsa}
If, in addition to the assumptions of \reft{realcoef}, $D$ is self-adjoint, then $m_-=0$ and $\Det_{2\tet}(D^2)$ is real and positive. Hence, as expected, formulas
\refe{det-eta-sa} and \refe{det-eta-sa2} hold.
\ecor

\noindent{\em Proof of \refc{detsa}.} \ If $D$ is self-adjoint, the spectrum of $D$ lies on the real line. Hence, in particular, $m_-=0$. It follows
from \refe{zet-Pizet} together with \refe{Imzet} and \refe{zet=ozet} below that $\IM\big(\zet'_{2\tet}(0,D^2)\big)= 0$. Hence, $\Det_{2\tet}(D^2)>0$.
\hfill$\square$
\rem{realcoef}
It is interesting to compare \refe{realcoef} with Theorem~3.2 of \cite{BrAbanov}. Suppose that the spectrum of $D$ is also symmetric about the imaginary axis. Then
$\eta(D)=0$. If, in addition, $\dim{M}$ is odd, then $\zet_{2\tet}(0,D^2)=0$, cf. \refr{det-eta}.c. Hence, \refe{realcoef} imply that, in this case, $\Det_\tet(D)$ is
real and its sign is equal to $(-1)^{m_-}$. Theorem~3.2 of \cite{BrAbanov} states that this is true
without the assumption that the spectrum of $D$ is symmetric about the real axis%
\footnote{In \cite{BrAbanov} the spectral cut was taken in the {\em upper half-plane}. Consequently, $m_-$ is replaced their by $m_+$.}
(i.e., for every invertible elliptic operator with self-adjoint leading symbol, whose spectrum is symmetric about the imaginary axis).
\erem
\noindent{\em Proof of \reft{realcoef}.} In view of \refe{zet-zet}, it is enough to consider the case when $\tet$ is sufficiently close to
$-\pi/2$ so that there are no eigenvalues of $D$ in the solid angles $L_{(-\pi/2,\tet]}$ and $L_{(\pi/2,\tet+\pi]}$, which we will henceforth
assume. By \refe{zet-Pizet}, \refe{eta-Pizet}, and \refe{zet'=2} it suffices to show that the numbers
\begin{multline}\notag
    \zet_{\tet}(0,\tilPi_+,D) \ \pm \ \zet_{\tet}(0,\tilPi_-,-D) \\ = \
    \big(\,\zet_{\tet}(0,\Pi_+,D) \ \pm \ \zet_{\tet}(0,\Pi_-,-D)\,\big) \ + \
    \big(\,\zet_{\tet}(0,P_+,D) \ \pm \ \zet_{\tet}(0,P_-,-D)\,\big)
\end{multline}
are real and that the imaginary part of the number
\begin{multline}\notag
    \zet_{\tet}'(0,\tilPi_+,D) \ + \ \zet_{\tet}'(0,\tilPi_-,-D) \\ = \
    \big(\,\zet_{\tet}'(0,\Pi_+,D) \ + \ \zet_{\tet}'(0,\Pi_-,-D)\,\big) \ + \
    \big(\,\zet_{\tet}'(0,P_+,D) \ + \ \zet_{\tet}'(0,P_-,-D)\,\big)
\end{multline}
is equal to \/ $-\pi{}m_-$.

Since the projections $P_\pm$ have finite rank, one has
\[
    \zet_{\tet}(0,P_\pm,\pm{}D) \ =\ \RANK{}P_\pm.
\]
Thus these numbers are real.

Because the spectrum of $D$ is symmetric about the real axis, $\RANK{}P_+= \RANK{}P_-$. As, for every $r>0$ one has
\[
    \frac{d}{ds}|_{s=0}\, (ir)^{-s}_{\tet} \ = \ -\log r \ - \ i\frac\pi2,
\]
we conclude that
\[
    \IM\,\zet_{\tet}'(0,P_+,D) \ = \ \IM\,{\zet_{\tet}'(0,P_-,-D)} \ = \ -\frac{\pi}2\,\RANK{}P_-.
\]
Hence,
\eq{Imzet}
    \IM\,\Big(\,\zet_{\tet}'(0,P_+,D) \ + \  \zet_{\tet}'(0,P_-,-D)\,\Big) \ = \ -\pi\,\RANK{}P_- \ \in \ \pi\,\ZZ.
\end{equation}
It remains to show that
\[
    \zet_{\tet}(0,\Pi_\pm,\pm D), \ \zet_{\tet}'(0,\Pi_\pm,\pm D)\ \in \ \RR.
\]
We will show that the numbers $\zet_{\tet}(0,\Pi_+,D)$ and $\zet_{\tet}'(0,\Pi_+,D)$ are real. The fact that the other two numbers are real as
well follows then by replacing $D$ with $-D$.

Let
\[
    \lam_j>0, \qquad j\in I_1\subset \NN
\]
be all the positive real eigenvalues of $D$ and let
\[
    \lam_j\ =\ \rho_j{}e^{i\alp_j}, \qquad j\in I_2\subset \NN
\]
be all the eigenvalues of $D$ which lie in the solid angle $L_{(0,\pi/2)}$. Let $m_j$ denote the algebraic multiplicity of $\lam_j$, cf.
\refss{det-sa}. Since the spectrum of $D$ is symmetric about the real axis,
\[
    \rho_j{}e^{-i\alp_j}, \qquad j\in I_2,
\]
are all the eigenvalues of $D$ in the solid angle $L_{(-\pi/2,0)}$ and
\begin{multline}\Label{E:zetPi-lam}\notag
    \zet_{\tet}(s,\Pi_+,D) \ = \ \sum_{j\in I_1}\,  m_j\,\lam_j^{-s} \ + \ \sum_{j\in I_2}\, m_j\,\rho_j^{-s}\, (e^{-is\alp_j}+e^{is\alp_j})
    \\ = \ \sum_{j\in I_1}\,  m_j\,\lam_j^{-s} \ + \ 2\sum_{j\in I_2}\, m_j\,\rho_j^{-s}\,\cos(s\alp_j), \qquad\qquad \RE s >\frac{\dim M}m.
\end{multline}
Hence,
\eq{zets-zetos}
    \zet_{\tet}(s,\Pi_+,D) \ = \ \overline{\zet_{\tet}(\os,\Pi_+,D)},  \qquad\qquad \RE s >\frac{\dim M}m.
\end{equation}
Since both sides of \refe{zets-zetos} are holomorphic functions of $s$, the equality \refe{zets-zetos} holds for all regular points of
$\zet_{\tet}(s,\Pi_+,D)$. In particular, $\zet_{\tet}(s,\Pi_+,D)$ is real for all real regular points. Hence, $\zet_{\tet}(0,\Pi_+,D)\in \RR$.
Since \refe{zets-zetos} implies
\eq{zet=ozet}
        \zet_{\tet}'(s,\Pi_+,D) \ = \ \overline{\zet_{\tet}'(\os,\Pi_+,D)},
\end{equation}
we conclude that the number $\zet_{\tet}'(0,\Pi_+,D)$ is also real. \hfill$\square$

%-------------------------------------------------------------------
%------------------------------------------------------------------
%\section{The logarithm of the determinant as a holomorphic function}\Label{S:prlogdethol}

%-------------------------------------------------------------------
%------------------------------------------------------------------
\section{Families of flat connections}\Label{S:famconGM}

In this appendix we review some of the results of \cite{GoldmanMillson88} and reformulate them in a form convenient for our purposes. These results
are used in \refs{analytic}.

%----------------
\subsection{Connections flat modulo lower order terms}\Label{SS:flatupto}
First, we introduce some definitions from \cite{GoldmanMillson88}, but we formulate them in a slightly different form which is more convenient for
our purposes.

Let $k[t]= k[t_1\nek t_r]$ denote the polynomial ring in $r$ variables  over a field $k$. Let $\grm\subset k[t]$ denote the unique maximal ideal of
$k[t]$ (the {\em augmentation ideal}), i.e., the ideal generated by $t_1\nek t_r$. Let $A_m=k[t]/\grm^{m+1}$. We denote by $G_m=\GL(n,A_m)$ the group
of matrices with entries in $A_m$.

Let $M$ be a manifold and let $E$ be a complex vector bundle over $M$. Suppose $\n$ is a flat connection on $E$. Let
\eq{connection}
    \n(t)\ = \n+ \sum_{0<|\alp|\le m}\, \ome_\alp{}t^\alp, \qquad t\in k^r,
\end{equation}
be a family of connections. Here $\alp\in (\ZZ_\ge)^r$ is a multi-index, $|\alp|=\alp_1+\cdots+\alp_r$, $t^\alp=t_1^{\alp_1}t_2^{\alp_2}\cdots
t_r^{\alp_r}$, and $\ome_\alp$ are smooth 1-forms with values in $\End E$. We say that the family $\n(t)$ is {\em flat modulo $t^{m+1}$} if
$\n(t)^2\in \grm^{m+1}$.

Fix a base point $x_*\in M$.  Given a continuous path $\phi:[0,1]\to M, \ \phi(0)=\phi(1)=x_*$, for any $t\in k^r$,  we denote by
$\Mon_{\n(t)}(\phi)$ the monodromy of $\n(t)$ along $\phi$, cf. \refe{ddxPhi}. If the family $\n(t)$ is flat modulo $t^{m+1}$ then, for any
pathes $\phi_i : [0,1] \to M,\ \phi_i(0) = \phi_i (1) = x_\ast\   (i=1,2)$,
\[
    \Mon_{\n(t)}(\phi_1) \ \equiv \ \Mon_{\n(t)}(\phi_2) \qquad \MOD \ \grm^{m+1}.
\]
Hence, we have a well defined representation
\eq{Mon}
    \Mon_{\n(t)}: \pi_1(M,x_*) \ \longrightarrow \ G_m.
\end{equation}

One says that two families of connections $\n_1(t)$ and $\n_2(t)$, which are flat modulo $t^{m+1}$  are {\em $A_m$-gauge equivalent} if there
exists a family of gauge transformations
\eq{gauge}
    \grg(t) \ = \ \grg_0 \ + \ \sum_{0<|\alp|\le m}\, \grg_\alp{}t^\alp
\end{equation}
where each $\grg_\alp$ is a gauge transformation of $E$, such that
\[
   \n_2(t) \ \equiv \grg(t)\cdot \n_1(t)\cdot \grg(t)^{-1} \qquad \MOD \ \grm^{m+1}.
\]

%----------------------
\subsection{Relationship between families of connections and families of representations of the fundamental group in $G_m$}\Label{SS:repr-conn}
Proposition~6.3 of \cite{GoldmanMillson88} states that {there is a one-to-one correspondence between the $A_m$-gauge equivalence classes of connections $\n(t)$ and
the isomorphism classes of representations $\gam(t)$ of $\pi_1(M)$ in $G_m$ given by the monodromy representation \refe{Mon}}. In other words, we have the following
\lem{GM}
(i)  \ \  For every family of representations $\gam(t):\pi_1(M,x_*)\to G_m$, there exists a flat modulo $t^{m+1}$ family of connections $\n(t)$
such that
\eq{Mon2}
    \Mon_{\n(t)} \ \equiv \ \gam(t) \qquad \mod \grm^{m+1}.
\end{equation}

(ii)\  \ Every two connections $\n_1(t)$ and $\n_2(t)$ which are of the form \refe{connection}, are flat modulo $t^{m+1}$, and satisfy
\refe{Mon2} are $A_m$-gauge equivalent, i.e., there exists a family of gauge transformations \refe{gauge} such that
\eq{gauge2}
   \n_2(t) \ \equiv \grg(t)\cdot \n_1(t)\cdot \grg(t)^{-1} \qquad \MOD \ \grm^{m+1}.
\end{equation}
Moreover, if\/ $\n_1(0)= \n_2(0)$ then one can choose $\grg(t)=\grg_0 + \sum_{0<|\alp|\le m}\, \grg_\alp{}t^\alp$ so that $\grg_0=\Id$.
\end{lemma}

%-------------------------------
\subsection{The case when $k=\CC$ or $\RR$}\Label{SS:CorR}
Suppose now that $k=\CC$ or $\RR$  and $r \in \ZZ _{\ge 1}$. Let $\calO\subset k^r$ be an open set and let $\n_\mu$\ ($\mu\in \calO$) be a
family of connections such that for some $\lam\in \calO$ we have
\eq{connection2}
    \n_\mu\ =\ \n_\lam\ +\ \sum_{0<|\alp|\le m}\, \ome_\alp{}(\mu-\lam)^\alp \ + \ o(|\mu-\lam|^{m}), \qquad \mu\in \calO,
\end{equation}
where $o(|\mu-\lam|^m)$ is understood in the sense of the Fr\'echet topology on $\calC(E)$ introduced in \refss{analofconnections}.

\renewcommand{\on}{\overline{\n}}

Denote $t=\mu-\lam$ and set
\eq{on(t)}
    \on(t)\ =\ \n_\lam+ \sum_{0<|\alp|\le m} \ome_\alp{}t^\alp.
\end{equation}
Then $\n_\mu \ = \ \on(\mu-\lam) + o(|\mu-\lam|^m)$. Hence, for every closed path $\phi:[0,1]\to M, \ \phi(0)=\phi(1)= x_*$ we have
\eq{Monn-Monon}
  \Mon_{\n_\mu}(\phi) \ = \ \Mon_{\on(\mu-\lam)}(\phi) \  + \ o(|\mu-\lam|^{m}),
\end{equation}
where $o(|\mu-\lam|^m)$ is understood in the sense of the Fr\'echet topology introduced in \refss{difopFr}. If the family $\on(t)$ is flat
modulo $t^{m+1}$, we will view $\Mon_{\n_\mu}$ as a map $\pi_1(M,x_*)\to G_m$ by identifying it with $\Mon_{\on(\mu-\lam)}$.

%------------------------------
\subsection{Application to real differentiable families of flat connections}\Label{SS:applrealdif}
Let $\calO\subset \CC$ be an open set. A family $\n_\mu$ ($\mu\in \calO$) of flat connections on $E$ is called {\em real differentiable} at
$\lam\in \calO$ if there exist $\ome_1,\ome_2\in \Ome^1(M,\End E)$ with
\eq{nmu-nlam1}
    \n_\mu \ = \ \n_\lam \ + \ \RE(\mu-\lam)\cdot\ome_1 \ + \ \IM(\mu-\lam)\cdot\ome_2 \ + \ o(\mu-\lam).
\end{equation}
(Again, $o(\mu-\lam)$ is understood in the sense of the Fr\'echet topology on $\calC(E)$ introduced in \refss{analofconnections}.)

%----------
\lem{realanal}
Let $\lam \in \CC$ and let $\calO\subset \CC$ be an open neighborhood of $\lam$ in $\CC$. Suppose that $\n_\mu$ ($\mu\in \calO$) is a family of
flat connections which is real differentiable at $\lam$, cf. \refe{nmu-nlam1}. Assume that the map
\[
    \calO\ \to\ \Rep, \qquad \mu\ \mapsto\  \Mon_{\n_\mu}
\]
is a holomorphic curve in $\Rep$. Then the following statements hold:

(i) \ There exists a smooth form $\ome\in \Ome^1(M,\End E)$ such that $\n_\lam\ome=0$ and
\eq{n+ome1}
  \Mon_{\n_\lam+(\mu-\lam)\ome}(\phi) \ = \ \Mon_{\n_\mu}(\phi) \ + \ o(\mu-\lam),
\end{equation}
for every closed path $\phi:[0,1]\to M,\ \phi(0)=\phi(1)=x_*$.

(ii) \ There exists a family of gauge transformations $G(\mu)\in \End{E}$ ($\mu\in \calO$) such that $G(\lam)=\Id$ and
\eq{guage1}
   \n_\lam \ + \ (\mu-\lam)\,\ome \ = \ G(\mu)\cdot\n_\mu\cdot G(\mu)^{-1} \ + \ o(\mu-\lam).
\end{equation}
\elem
\prf
To prove part (i) of the lemma we apply \refl{GM} with $k=\CC$, $t=t_1$ (i.e., $r=1$), $m=1$, and $t=\mu-\lam$. Since $\mu\to \Mon_{\n_\mu}$ is
a holomorphic curve, its Taylor expansion at $\lam$ up to first order, $\gamma(t)$, defines a map $\calO \to G_1$. Then
\eq{Mon=gam}
    \Mon_{\n_{\lam +t\ome}}\ =\ \gam(t) \ + \ o(t).
\end{equation}
By \refl{GM}(i), there exists a flat modulo $t^2$ family of connections $\tiln(t)=\tiln(0)+t\tilome$ such that
\eq{Montiln=Mon}
    \Mon_{\tiln(t)} \ \equiv \ \gam(t) \qquad \MOD\ t^2.
\end{equation}
Since $\Mon_{\tiln(0)}= \gam(0)$ there exists a gauge transformation $\grg\in \End E$ such that its restriction to the fiber of $E$ over the
base point $x_*$ is the identity map and
\[
  \n_\lam \ = \ \grg^{-1}\cdot \tiln(0)\cdot \grg.
\]
Then $\ome:= \grg^{-1}\,\tilome\,\grg$ is a smooth $\End E$-valued 1-form and \refe{Montiln=Mon} takes the form
\eq{Monnlam+t=gam(t)}
    \Mon_{\n_\lam+t\ome} \ \equiv \ \gam(t) \qquad \MOD\ t^2
\end{equation}
which together with \refe{Mon=gam} implies \refe{n+ome1}. Note that $\n_{\lam} +t\ome$ is a flat modulo $t^2$ connection and, hence,
$\n_{\lam}\ome  = 0$.

For part (ii) let us set $k=\RR$, $t=(t_1,t_2)$ (i.e., $r=2$), and $m=1$. Denote $t_1:=\RE(\mu-\lam), \ t_2:=\IM(\mu-\lam)$. Then, by the
assumption of real differentiability,  $\n_\mu $ is of the form
\[
       \n_{\lam} \ +\ t_1\, \ome_1  \ +\ t_2\,\ome_2\  =\ \n_{\mu} \ +\ o(\mu -\lam).
\]
Note that both, $\n_{\lam} +t_1 \ome_1  +t_2\ome_2$ and $\n_{\lam} +(t_1 + i t_2)\ome$, are flat modulo $t^2$ connections which, by
\refe{Mon=gam}, induce the same monodromy representation $\gam(t):\pi_1(M,x_*)\to G_1$. Hence, \refe{guage1} follows from \refl{GM}(ii).
\eprf

%-------------------------------------------------------------------------------------
%-------------------------------------------------------------------------------------
\providecommand{\bysame}{\leavevmode\hbox to3em{\hrulefill}\thinspace} \providecommand{\MR}{\relax\ifhmode\unskip\space\fi MR }
\providecommand{\MRhref}[2]{%
  \href{http://www.ams.org/mathscinet-getitem?mr=#1}{#2}
} \providecommand{\href}[2]{#2}


\begin{thebibliography}{10}

\bibitem{BrAbanov}
A.~Abanov and M.~Braverman, \emph{{Topological calculation of the phase of the
  determinant of a non self-adjoint elliptic operator}},
   \texttt{arXiv:math-ph/0401037}, to appear in Comm. of Math. Physics.

\bibitem{APS1}
M.~F. Atiyah, V.~K. Patodi, and I.~M. Singer, \emph{Spectral asymmetry and
  {R}iemannian geometry. {I}}, Math. Proc. Cambridge Philos. Soc. \textbf{77}
  (1975), 43--69.

\bibitem{APS2}
\bysame, \emph{Spectral asymmetry and {R}iemannian geometry. {II}}, Math. Proc.
  Cambridge Philos. Soc. \textbf{78} (1975), no.~3, 405--432.

\bibitem{BeGeVe}
N.~Berline, E.~Getzler, and M.~Vergne, \emph{Heat kernels and {Dirac}
  operators}, Springer-Verlag, 1992.

\bibitem{BisZh92}
J.-M. Bismut and W.~Zhang, \emph{An extension of a theorem by {Cheeger} and
  {M\"uller}}, Ast\'erisque \textbf{205} (1992).

\bibitem{BrKappelerRATshort}
M.~Braverman and T.~Kappeler,  \emph{A refinement of the {Ray}-{Singer} torsion}, C.R. Acad. Sci.
  Paris \textbf{341} (2005), 497--502.


\bibitem{BrKappelerRATdetline}
\bysame, \emph{{Refined Analytic Torsion as an Element of the Determinant
  Line}}, IHES preprint M/05/49, \texttt{math.GT/0510532}.

\bibitem{BruningLesch99}
J.~Br{\"u}ning and M.~Lesch, \emph{On the {$\eta$}-invariant of certain
  nonlocal boundary value problems}, Duke Math. J. \textbf{96} (1999), no.~2,
  425--468.

\bibitem{BFK91}
D.~Burghelea, L.~Friedlander, and T.~Kappeler, \emph{On the determinant of
  elliptic differential and finite difference operators in vector bundles over
  {$S\sp 1$}}, Comm. Math. Phys. \textbf{138} (1991), no.~1, 1--18.

\bibitem{BFK92}
D.~Burghelea, L.~Friedlander, and T.~Kappeler, \emph{Mayer-{Vietoris} type
  formula for determinants of elliptic differential operators}, Journal of
  Funct. Anal. \textbf{107} (1992), 34--65.

\bibitem{BFK3}
\bysame, \emph{Asymptotic expansion of the {Witten} deformation of the analytic
  torsion}, Journal of Funct. Anal. \textbf{137} (1996), 320--363.

\bibitem{BurgheleaHaller_Euler}
D.~Burghelea and S.~Haller, \emph{{{Euler} Structures, the Variety of
  Representations and the {Milnor}-{Turaev} Torsion}},
  \texttt{arXiv:math.DG/0310154}.

\bibitem{BurgheleaHaller_function}
D.~Burghelea and S.~Haller, \emph{Torsion, as a function on the space of
  representations}, \texttt{arXiv:math.DG/0507587}.
\bibitem{Cheeger79}
J.~Cheeger, \emph{Analytic torsion and the heat equation}, Ann. of Math.
  \textbf{109} (1979), 259--300.

\bibitem{Dineen99}
S.~Dineen, \emph{Complex analysis on infinite-dimensional spaces}, Springer
  Monographs in Mathematics, Springer-Verlag London Ltd., London, 1999.

\bibitem{Farber00AT}
M.~Farber, \emph{Absolute torsion and eta-invariant}, Math. Z. \textbf{234}
  (2000), no.~2, 339--349.

\bibitem{FarberLevine96}
M.~Farber and J.~Levine, \emph{Jumps of the eta-invariant}, Math. Z.
  \textbf{223} (1996), 197--246, With an appendix by S.~Weinberger:
  Rationality of $\rho$-invariants.

\bibitem{FarberTuraev99}
M.~Farber and V.~Turaev, \emph{Absolute torsion}, Tel Aviv Topology Conference:
  Rothenberg Festschrift (1998), Contemp. Math., vol. 231, Amer. Math. Soc.,
  Providence, RI, 1999, pp.~73--85.

\bibitem{FarberTuraev00}
\bysame, \emph{Poincar\'e-{R}eidemeister metric, {E}uler structures, and
  torsion}, J. Reine Angew. Math. \textbf{520} (2000), 195--225.

\bibitem{Gilkey84}
P.~B. Gilkey, \emph{The eta invariant and secondary characteristic classes of
  locally flat bundles}, Algebraic and differential topology---global
  differential geometry, Teubner-Texte Math., vol.~70, Teubner, Leipzig, 1984,
  pp.~49--87.

\bibitem{GohbergKrein_book69}
I.~C. Gohberg and M.~G. Kre\u{\i}n, \emph{Introduction to the theory of linear
  nonselfadjoint operators}, Translated from the Russian by A. Feinstein.
  Translations of Mathematical Monographs, Vol. 18, American Mathematical
  Society, Providence, R.I., 1969.

\bibitem{GoldmanMillson88}
W.~Goldman and J.~Millson, \emph{The deformation theory of representations of
  fundamental groups of compact {K}\"ahler manifolds}, Inst. Hautes \'Etudes
  Sci. Publ. Math. (1988), no.~67, 43--96.

\bibitem{Grubbseeley95}
G.~Grubb and R.~Seeley, \emph{Weakly parametric pseudodifferential operators
  and {A}tiyah-{P}atodi-{S}inger boundary problems}, Invent. Math. \textbf{121}
  (1995), no.~3, 481--529.

\bibitem{Guillemin85}
V.~Guillemin, \emph{A new proof of {W}eyl's formula on the asymptotic
  distribution of eigenvalues}, Adv. in Math. \textbf{55} (1985), 131--160.

\bibitem{HormanderSCV}
L.~H{\"o}rmander, \emph{An introduction to complex analysis in several
  variables}, third ed., North-Holland Mathematical Library, vol.~7,
  North-Holland Publishing Co., Amsterdam, 1990.

\bibitem{KontsevichVishik_long}
M.~Kontsevich and S.~Vishik, \emph{{Determinants of elliptic
  pseudo-differential operators}}, Preprint MPI /94-30,
  \texttt{arXiv:hep-th/9404046}.

\bibitem{KontsevichVishik_short}
\bysame, \emph{Geometry of determinants of elliptic operators}, Functional
  analysis on the eve of the 21st century, Vol.\ 1 (New Brunswick, NJ, 1993),
  Progr. Math., vol. 131, Birkh\"auser Boston, Boston, MA, 1995, pp.~173--197.

\bibitem{Lidskii59}
V.~B. Lidski\u{\i}, \emph{Non-selfadjoint operators with a trace}, Dokl. Akad.
  Nauk SSSR \textbf{125} (1959), 485--487.

\bibitem{Markus88}
A.~S. Markus, \emph{Introduction to the spectral theory of polynomial operator
  pencils}, Translations of Mathematical Monographs, vol.~71.

\bibitem{Massey60}
W.~S. Massey, \emph{On the {S}tiefel-{W}hitney classes of a manifold}, Amer. J.
  Math. \textbf{82} (1960), 92--102.

\bibitem{Mu1ller78}
W.~M\"uller, \emph{Analytic torsion and {R}-torsion on {Riemannian} manifolds},
  Adv. in Math. \textbf{28} (1978), 233--305.

\bibitem{Ponge-asymetry}
R.~Ponge, \emph{{Spectral Asymmetry, Zeta Functions and the Noncommutative
  Residue}}, \texttt{arXiv:math.DG/0310102}.

\bibitem{RaSi1}
D.~B. Ray and I.~M. Singer, \emph{{R}-torsion and the {Laplacian} on
  {Riemannian} manifolds}, Adv. in Math. \textbf{7} (1971), 145--210.

\bibitem{Retherford93}
J.~R. Retherford, \emph{Hilbert space: compact operators and the trace
  theorem}, London Mathematical Society Student Texts, vol.~27, Cambridge
  University Press, Cambridge, 1993.

\bibitem{Rudin_FA}
W.~Rudin, \emph{Functional analysis}, second ed., International Series in Pure
  and Applied Mathematics, McGraw-Hill Inc., New York, 1991.

\bibitem{Rudyak_book}
Y.~B. Rudyak, \emph{On {T}hom spectra, orientability, and cobordism}, Springer
  Monographs in Mathematics, Springer-Verlag, Berlin, 1998, With a foreword by
  Haynes Miller.

\bibitem{Seeley67}
R.~Seeley, \emph{Complex powers of elliptic operators}, Proc. Symp. Pure and
  Appl. Math. AMS \textbf{10} (1967), 288--307.

\bibitem{ShubinPDObook}
M.~A. Shubin, \emph{Pseudodifferential operators and spectral theory}, Springer
  Verlag, Berlin, New York, 1980.

\bibitem{Singer85Dirac}
I.~M. Singer, \emph{Families of {D}irac operators with applications to
  physics}, Ast\'erisque (1985), no.~Numero Hors Serie, 323--340, The
  mathematical heritage of \'Elie Cartan (Lyon, 1984).

\bibitem{Turaev86}
V.~G. Turaev, \emph{{R}eidemeister torsion in knot theory}, Russian Math.
  Survey \textbf{41} (1986), 119--182.

\bibitem{Turaev90}
\bysame, \emph{Euler structures, nonsingular vector fields, and
  {R}eidemeister-type torsions}, Math. USSR Izvestia \textbf{34} (1990),
  627--662.

\bibitem{Wall60}
C.~T.~C. Wall, \emph{Determination of the cobordism ring}, Ann. of Math. (2)
  \textbf{72} (1960), 292--311.

\bibitem{Warner}
F.~W. Warner, \emph{Foundations of differentiable manifolds and {Lie} groups},
  Graduate Texts in Mathematics, Springer-Verlag, New-York, Berlin, Heidelberg,
  Tokyo, 1983.

\bibitem{Wodzicki84}
M.~Wodzicki, \emph{Local invariants of spectral asymmetry}, Invent. Math.
  \textbf{75} (1984), no.~1, 143--177.

\bibitem{Wodzicki87}
\bysame, \emph{Noncommutative residue. {I}. {F}undamentals}, $K$-theory,
  arithmetic and geometry (Moscow, 1984--1986), Lecture Notes in Math., vol.
  1289, Springer, Berlin, 1987, pp.~320--399.

\bibitem{Wojciechowski99}
K.~P. Wojciechowski, \emph{Heat equation and spectral geometry. {I}ntroduction
  for beginners}, Geometric methods for quantum field theory (Villa de Leyva,
  1999), World Sci. Publishing, River Edge, NJ, 2001, pp.~238--292.

\bibitem{ZimmerFunctAn}
R.~J. Zimmer, \emph{Essential results of functional analysis}, Chicago Lectures
  in Mathematics, University of Chicago Press, Chicago, IL, 1990.

\end{thebibliography}
\end{document}